\ifx \labelsONmargin\undefined
  \ifx\RequirePackage\undefined
    \expandafter\ifx\csname amsart.sty\endcsname\relax
        \documentstyle[12pt,amssymb,amscd]{amsart}
    \fi

    \def\atSign{@@}

    \def\mathbb{\Bbb}
    \def\mathfrak{\frak}
    \def\mathbf{\bold}
    \ifx\boldsymbol\undefined	
      \def\boldsymbol#1{{\bold #1}}
    \fi

    \def\mathbit{\boldsymbol}

    \makeatletter
    \newenvironment{proof}{%
         \@ifnextchar[{%
                       \expandafter\let\expandafter\end@proof
                         \csname endpf*\endcsname
                         \my@proof
                      }{\let\end@proof\endpf\pf}%
        }{\end@proof}
    \def\my@proof[#1]{\@nameuse{pf*}{#1}}
    \makeatother

    \def\xrightarrow[#1]#2{@>{#2}>{#1}>}
    \def\xleftarrow[#1]#2{@<{#2}<{#1}<}

    \def\providecommand#1{\def#1}
    \def\emph#1{{\em #1}}
    \def\textbf#1{{\bf #1}}
  \else
      \ifx\usepackage\RequirePackage
      \else
        \documentclass[12pt]{amsart}
        \usepackage{amssymb,amscd}
      \fi

      \ifx\mathbit\undefined	
        \DeclareMathAlphabet{\mathbit}{OML}{cmm}{b}{it}
      \fi

      \def\atSign{@}

      \advance\oddsidemargin -0.1\textwidth
      \evensidemargin\oddsidemargin
      \def\Sb#1\endSb{_{\substack{#1}}}
      \def\Sp#1\endSp{^{\substack{#1}}}
  \fi
\makeatletter
\@ifundefined{DeclareFontShape} 
     {\@ifundefined{selectfont}
             {
                \message{We need AmS-LaTeX, this version of LaTeX does not 
                        support it}\@@end
             }{
                \def\mathcal{\cal}
                
                \def\pcyr{%
                        \def\default@family{UWCyr}%
                        \let\oldSl@\sl
                        \def\sl{\def\default@shape{it}\oldSl@}%
                        \cyracc
                        \language\Russian\family{UWCyr}\selectfont
                }
             }}{
                \DeclareFontEncoding{OT2}{}{\noaccents@}
                \DeclareFontSubstitution{OT2}{cmr}{m}{n}
                \DeclareFontFamily{OT2}{cmr}{\hyphenchar\font45 }
                \DeclareFontShape{OT2}{cmr}{m}{n}{%
                     <5><6><7><8><9><10>gen*wncyr %
                     <10.95><12><14.4><17.28><20.74><24.88> wncyr10 %
                }{}
                \DeclareFontShape{OT2}{cmr}{m}{it}{%
                     <5><6><7><8><9><10> gen * wncyi%
                     <10.95><12><14.4><17.28><20.74><24.88> wncyi10%
                }{}
                \DeclareFontShape{OT2}{cmr}{bx}{n}{%
                     <5><6><7><8><9><10> gen * wncyb%
                     <10.95><12><14.4><17.28><20.74><24.88> wncyb10%
                }{}
                \DeclareFontShape{OT2}{cmr}{m}{sl}{%
                     <-> ssub * cmr/m/it%
                }{}
                \DeclareFontShape{OT2}{cmr}{m}{sc}{%
                     <5><6><7><8><9><10>%
                     <10.95><12><14.4><17.28><20.74><24.88> wncysc10%
                }{}
                \DeclareFontFamily{OT2}{cmss}{\hyphenchar\font45 }
                \DeclareFontShape{OT2}{cmss}{m}{n}{%
                     <8><9><10> gen * wncyss%
                     <10.95><12><14.4><17.28><20.74><24.88> wncyss10%
                }{}
                
                \def\cyrencodingdefault{OT2}
                \def\pcyr{%
                        \cyracc
                        \let\encodingdefault\cyrencodingdefault
                        \language\Russian\fontencoding{OT2}\selectfont
                }
             }   

\@ifundefined{theorembodyfont} 
     {
        \def\theorembodyfont#1{\relax}
        \makeatletter
          \let\@@th@plain\th@plain
          \def\th@plain{ \@@th@plain \slshape }
        \makeatother
        \let\normalshape\relax
     }{}   

\ifx\cprime\undefined
     \def\cprime{$'$}
\fi
\makeatletter
  \def\@sect@my#1#2#3#4#5#6[#7]#8{%
\ifnum #2>\c@secnumdepth
   \let\@svsec\@empty
 \else
   \refstepcounter{#1}%
\edef\@svsec{\ifnum#2<\@m
             \@ifundefined{#1name}{}{\csname #1name\endcsname\ }\fi
\noexpand\rom{\csname the#1\endcsname.}\enspace}\fi
 \@tempskipa #5\relax
 \ifdim \@tempskipa>\z@ 
   \begingroup #6\relax
   \@hangfrom{\hskip #3\relax\@svsec}{\interlinepenalty\@M #8\par}%
   \endgroup
   \if@article\else\csname #1mark\endcsname{%
        \ifnum \c@secnumdepth >#2\relax\csname the#1\endcsname. \fi#7}\fi
\ifnum#2>\@m \else
       \let\@tempf\\ \def\\{\protect\\}\addcontentsline{toc}{#1}%
{\ifnum #2>\c@secnumdepth \else
             \protect\numberline{%
               \ifnum#2<\@m
               \@ifundefined{#1name}{}{\csname #1name\endcsname\ }\fi
               \csname the#1\endcsname.}\fi
           #8}\let\\\@tempf
     \fi
 \else
  \def\@svsechd{#6\hskip #3\@svsec
    \@ifnotempty{#8}{\ignorespaces#8\unskip
       \ifnum\spacefactor<1001.\fi}%
        \ifnum#2>\@m \else
          \let\@tempf\\ \def\\{\protect\\}\addcontentsline{toc}{#1}%
            {\ifnum #2>\c@secnumdepth \else
              \protect\numberline{%
                \ifnum#2<\@m
                \@ifundefined{#1name}{}{\csname #1name\endcsname\ }\fi
                \csname the#1\endcsname.}\fi
             #8}\let\\\@tempf\fi}%
 \fi
\@xsect{#5}}
  \ifx\RequirePackage\undefined
  \let\@sect\@sect@my             
  \fi
  \def\th@remark@my{\theorempreskipamount6\p@\@plus6\p@
    \theorempostskipamount\theorempreskipamount
    \def\theorem@headerfont{\it}\normalshape}
  \ifx\RequirePackage\undefined
  \let\th@remark\th@remark@my
  \fi
\let\myLabel\@gobble
\def\labelsONmargin{\@mparswitchfalse\def\myLabel##1{\@bsphack\marginpar
                                  {\normalshape\tiny\rm Label ##1}\@esphack}}
\ifx \url\undefined
  \def\url#1{{\tt #1}}%
\fi
\def\cyracc{\def\u##1{
                \if \i##1\char"1A%
                \else \if I##1\char"12%
                \else \accent"24 ##1\fi\fi }%
\def\"##1{\if e##1{\char"1B}%
                \else \if E##1{\char"13}%
                \else \accent"7F ##1\fi\fi }%
\def\9##1{\if##1z\char"19 
\else\if##1Z\char"11 
\else\if##1E\char"03 
\else\if##1e\char"0B 
\else\if##1u\char"18 
\else\if##1U\char"10 
\else\if##1A\char"17 
\else\if##1a\char"1F 
\else\if##1p\char"7E 
\else\if##1P\char"5E 
\else\if##1Q\char"5F 
\else\if##1q\char"7F 
\else\if##1i\char"1A 
\else\if##1I\char"12 
\else\if##1N\char"7D 
\fi
\fi
\fi
\fi
\fi
\fi
\fi
\fi
\fi
\fi
\fi
\fi
\fi
\fi
\fi
}%
\def\cydot{{\kern0pt}}}%
\def\cydot{$\cdot$}

\ifx\Russian\undefined
        \def\Russian{0\relax
    \message{Don't know the hyphenation rules for Russian^^J
                        Please do INITeX with `input  russhyph' in the 
                        command line}%
                \gdef\Russian{0\relax}%
        }
\fi
\def\@putname#1#2#3#4{\def\@@ref{#3}\let\old@bf\bf	
	\let\old@reset@font\reset@font			
        \def\bf##1{\old@bf\if?\noexpand##1?{#4}\else##1\fi}%
	\def\reset@font##1##2{\old@reset@font##1\if?\noexpand##2?{#4}\else##2\fi}#1{#2}%
        \let\bf\old@bf\let\reset@font\old@reset@font}
\let\my@ref=\ref
\def\ref#1{\@putname\my@ref{#1}{#1}{\tiny\rm\@@ref}}
\let\my@pageref=\pageref
\def\pageref#1{\@putname\my@pageref{#1}{#1}{\tiny\rm\@@ref}}
\let\my@cite=\cite
\def\cite#1{\@putname\my@cite{#1}{\@citeb}{\tiny\rm\@@ref}}
\makeatother
\fi
\relax 
\theoremstyle{plain} 
\theorembodyfont{\sl}

\newtheorem{nwthrmi}{Question } 
\newtheorem{nwthrmii}{Question }

\numberwithin{equation}{section}
\theorembodyfont{\rm}
\theoremstyle{definition}

\newtheorem{definition}{Definition}[section]
\newtheorem{conjecture}[definition]{Conjecture}
\newtheorem{example}[definition]{Example}

\theoremstyle{remark}

\newtheorem{remark}[definition]{Remark} 

\theoremstyle{plain} 
\theorembodyfont{\sl}

\newtheorem{theorem}[definition]{Theorem}
\newtheorem{lemma}[definition]{Lemma}
\newtheorem{corollary}[definition]{Corollary}
\newtheorem{proposition}[definition]{Proposition}
\newtheorem{amplification}[definition]{Amplification}

\begin{document}
\bibliographystyle{amsplain}
\relax 

\title[ Local geometry of bihamiltonian structures]{ Webs, Lenard schemes, and
the local geometry of bihamiltonian Toda and Lax structures }

\author{ Israel~M.~Gelfand }

\author{ Ilya Zakharevich }

\address{ Dept. of Mathematics, Rutgers University, Hill Center, New
Brunswick, NJ, 08903}

\email{igelfand\atSign{}math.rutgers.edu}

\address{ Department of Mathematics, Ohio State University, 231 W.~18~Ave,
Columbus, OH, 43210 }

\email {ilya\atSign{}math.ohio-state.edu}

\date{ March 1999 (Revision III: March 2000) Archived as
\url{math.DG/9903080} Printed: \today }

\setcounter{section}{-1}

\maketitle
\begin{abstract}
We introduce a criterion that a given bihamiltonian structure
admits a local coordinate system where both brackets have constant
coefficients. This criterion is applied to the bihamiltonian open Toda
lattice in a generic point, which is shown to be locally isomorphic to a
Kronecker odd-dimensional pair of brackets with constant coefficients.
This shows that the open Toda lattice cannot be locally represented as a
product of two bihamiltonian structures.

In a generic point the bihamiltonian periodic Toda lattice is shown
to be isomorphic to a product of two open Toda lattices (one of which is
a (trivial) structure of dimension 1).

While the above results might be obtained by more traditional
methods, we use an approach based on general results on geometry of webs.
This demonstrates a possibility to apply a geometric language to problems
on bihamiltonian integrable systems, such a possibility may be no less
important than the particular results proven in this paper.

Based on these geometric approaches, we conjecture that
decompositions similar to the decomposition of the periodic Toda lattice
exist in local geometry of the Volterra system, the complete Toda
lattice, the multidimensional Euler top, and a regular bihamiltonian Lie
coalgebra. We also state general conjectures about geometry of more
general ``homogeneous'' finite-dimensional bihamiltonian structures.

The class of homogeneous structures is shown to coincide with the class
of system integrable by Lenard scheme. The bihamiltonian structures which
admit a non-degenerate Lax structure are shown to be locally isomorphic to
the open Toda lattice.

\end{abstract}
\tableofcontents

\section{Introduction }\label{h01}\myLabel{h01}\relax 

A {\em local-geometric approach\/} consists of considering a geometric
structure (for the purpose of our discussion this is a collection of
tensor fields) up to a local diffeomorphism, studying its local
automorphisms, invariant tensor fields for these automorphisms, and a
possibility to decompose the structure into direct products. When applied
to integrable systems, this accounts to forgetting all the information
related to the given coordinate system (say, whether the structure is
polynomial {\em in this system\/}).

This approach cannot explain the phenomenon of integrability of a
Hamiltonian system, when the initial geometric structure is a Poisson
bracket and a function on a manifold. This local geometric structure has
too large group of automorphism, and there is no additional invariant
functions one could have used to integrate the system. One needs global
(or non-invariant) data to integrate a Hamiltonian system.

There is an alternative {\em bihamiltonian\/} approach to dynamic
systems in which integrability becomes meaningful on the local level
already. In this approach one starts with two {\em compatible\/}\footnote{Two Poisson brackets $ \left\{,\right\}_{1} $ and $ \left\{,\right\}_{2} $ on $ M $ are {\em compatible\/} if
the bracket $ \lambda_{1}\left\{,\right\}_{1} +\lambda_{2}\left\{,\right\}_{2} $ is Poisson for any $ \lambda_{1} $, $ \lambda_{2} $.
} Poisson
brackets $ \left\{,\right\}_{1} $ and $ \left\{,\right\}_{2} $ on $ M $. Basing on these brackets one constructs a
dynamical system which is Hamiltonian with respect to any one of these
brackets (and in fact to any linear combination of the brackets). The
construction of the dynamical system basing on the brackets is called
{\em Lenard scheme}. It provides a family of functions in involution (w.r.t.~
any linear combination of the brackets). Considering any function of this
family as a Hamiltonian w.r.t.~any bracket of two one obtains many
Hamiltonian flows. In most cases which appear in practice the above
family of functions is large enough to make these dynamics integrable
(compare with examples in Section~\ref{h48} and statements of Section~\ref{h55}).

Lenard scheme was formalized in \cite{Lax76Alm,Mag78Sim,%
GelDor79Ham,FokFuch80Str}, see also \cite{KosMag96Lax}. Most of these
formalizations assume that at least one of the brackets is symplectic\footnote{Any symplectic structure carries a Poisson bracket. We call such Poisson
brackets {\em symplectic}.}
(thus $ M $ is even-dimensional). That time it was not realized how these
formalizations relate to known applications of Lenard scheme, which
consist of a recurrence relation, and of initial data for these
relations. The above formalizations of \cite{Lax76Alm,Mag78Sim,%
GelDor79Ham,FokFuch80Str,KosMag96Lax} studied the recurrence
relations only, ignoring the initial data.

When even-dimensional bihamiltonian structures were classified in
\cite{Tur89Cla,Mag88Geo,Mag95Geo,McKeanPC,GelZakh93}, it became
clear that there is exactly one case where the above ``symplectic''
formalizations are compatible with the initial data for recurrence. This
case is in no way analogous to known examples (see Remark~\ref{rem55.80}).

Later, when the analysis of \cite{GelZakhFAN,GelZakh94Spe} had shown
that the periodic KdV system should be considered as an odd-dimensional
(though infinite-dimensional) bihamiltonian structure, an alternative
approach to the Lenard scheme became necessary. The philosophy of
\cite{GelZakhWeb} and \cite{GelZakh93} is that such a substitute is given by the
local classification of bihamiltonian structures.

By this philosophy the mentioned above ``symplectic'' formalizations
of Lenard scheme are substituted by the local descriptions of generic
even-dimensional bihamiltonian structures in \cite{Tur89Cla,Mag88Geo,%
Mag95Geo,McKeanPC,GelZakh93}. Indeed, these descriptions provide
all the information contained in \cite{Mag78Sim} and \cite{GelDor79Ham}, and
demystify the assumptions of the former papers.

From the classification of even-dimensional bihamiltonian structures
in general position, it turns out that this geometry is pretty rigid: on
an open subset the structure may be canonically decomposed into a direct
product of two-dimensional components, with one distinguished canonically
defined coordinate on each of these components. (It is this rigidity
which allows a local construction of a big family of commuting
Hamiltonians.) However, as in the case of a Hamiltonian system, locally
it has discrete parameters only (up to minor details the only parameter
is dimension). The morale of this classification is that only
$ 2 $-dimensional geometry is important, anything else can be combined from
$ 2 $-dimensional building blocks.

The situation becomes very different in an odd-dimensional case: the
structures in general position are indecomposable. In fact such
structures are even {\em micro-indecomposable}, i.e., one cannot represent them
as a product of two structures of smaller dimension---even if one
restricts attention to one tangent space to a point of the manifold. For
{\em analytical\/} structures in general position a local classification is also
possible (\cite{GelZakhWeb,GelZakh93}), but it is equivalent to a (local)
classification of {\em non-linear\/} $ 1 $-dimensional bundles over a rational curve,
i.e., analytical surfaces which have a submanifold isomorphic to $ {\mathbb P}^{1} $ and a
fixed projection onto this curve\footnote{Dimension of the initial bihamiltonian structure depends on the
degree of the normal (line) bundle to this curve.}. This classification involves
functional parameters (several functions of two complex variables).

The geometry of such bihamiltonian structures is also very rigid,
thus basing on local geometric data one can canonically construct enough
functions in involution, thus produce integrable systems. Out of this
huge pool of micro-indecomposable integrable systems of the given odd
dimension one can single out one particular {\em flat\/} structure, with {\em both\/}
Poisson structures having constant coefficients in the same coordinate
system (any two flat odd-dimensional indecomposable structures are
locally isomorphic, compare \cite{GelZakhFAN}).

However, after the heuristic of \cite{GelZakhFAN} that the KdV system is
in fact an infinite-dimensional analogue of an odd-dimensional
bihamiltonian structure, no other bihamiltonian structure were
(explicitly) considered from the point of view of classification up to a
diffeomorphism\footnote{Since the geometry of many ``classical'' bihamiltonian structures is
investigated up to minor details, a specialist could easily concoct an
answer to such a question from the known results. The conjectural reason
why this was not done before is that the answer would not fit into the
fixed mindset of ``everything is a product of $ 2 $-dimensional components'',
compare with discussion in Section~\ref{h005}.}. One of the targets of this paper is to investigate from
this point of view the simplest classical bihamiltonian structures: the
open and the periodic finite-dimensional Toda lattices.

While we proceed to this goal, we also provide generally-useful
easy-to-check criteria of flatness, investigate Lenard scheme in context
of odd-dimensional bihamiltonian geometry, and provide geometric
description of systems which admit a Lax representation.

For a detailed overview of the presented results in Section~\ref{h003} we
need to introduce some notions which are going to be used throughout the
paper. We do this in Section~\ref{h002}. Here we only list the principal steps
of our presentation:
\begin{enumerate}
\item
criteria of being homogeneous and being Kronecker of corank 1;
\item
introduction of webs as a way to encode mutual positions of Casimir
functions;
\item
proof of the criteria;
\item
examples of bihamiltonian structures which demonstrate purposes of
different conditions of the criteria;
\item
relation of Lenard integrability and homogeneous structures;
\item
relation of Lax structures and flatness;
\item
application of criteria to Toda lattices.
\end{enumerate}
We also discuss geometric conjecture which might provide geometric
description of many other finite-dimensional bihamiltonian structures.

Authors are indebted to A.~S.~Fokas, A.~Givental,
Y.~Kosmann-Schwarzbach, F.~Magri, H.~McKean, T.~Ratiu, N.~Reshetikhin for
fruitful discussions, and to A.~Gorokhovsky, M.~Braverman, B.~Khesin,
A.~Panasyuk, and V.~Serganova for the remarks which lead to improvements
of this paper. Special thanks go to A.~Panasyuk for letting us see the
preprint of \cite{Pan99Ver} before it went to print, and to M.~Gekhtman for
his suggestions on using known B\"acklund--Darboux transformation for
Volterra systems.

{\bf Revisions: }The revision II of this paper (January 2000) introduced
references to new papers \cite{Tur99Equi} and \cite{Zakh99Kro}, expanded
bibliography on ``classical'' bi-Hamiltonian systems, and minor stylistic
corrections. The revision III (March 2000) added Remark~\ref{rem48.91}.
Numbering of statements did not change. The archive name of this paper is
\url{math.DG/9903080} at \url{http://arXiv.org/math/abs}.

\section{Basic notions }\label{h002}\myLabel{h002}\relax 

All the geometric definitions which follow are applicable in $ C^{\infty} $ and
analytic geometry. We state only the $ C^{\infty} $-variant, the analytic one can be
obtained by substituting $ {\mathbb R} $ by $ {\mathbb C} $.

In what follows if $ f $ is a function or a tensor field on $ M $, $ f|_{m} $
denotes the value of $ f $ at $ m\in M $.

\begin{definition} A {\em bracket\/} on a manifold $ M $ is a $ {\mathbb R} $-bilinear skewsymmetric
mapping $ f,g \mapsto \left\{f,g\right\} $ from pairs of smooth functions on $ M $ to smooth
functions on $ M $. This mapping should satisfy the Leibniz identity
$ \left\{f,gh\right\}=g\left\{f,h\right\}+h\left\{f,g\right\} $. A bracket is {\em Poisson\/} if it satisfies Jacobi
identity too (thus defines a structure of a Lie algebra on functions on
$ M $).

A {\em Poisson structure\/} is a manifold $ M $ equipped with a Poisson
bracket. \end{definition}

\begin{remark} \label{rem002.20}\myLabel{rem002.20}\relax  Leibniz identity implies $ \left\{f,g\right\}|_{m}=0 $ if $ f $ has a zero of
second order at $ m\in M $. Thus a bracket is uniquely determined by describing
functions $ \left\{f_{i},f_{j}\right\} $, here $ \left\{f_{i}\right\}_{i\in I} $ is an arbitrary collection of smooth
functions on $ M $ such that for any $ m\in M $ the collection $ \left\{df_{i}|_{m}\right\}_{i\in I} $ of vectors
in $ {\mathcal T}_{m}^{*}M $ generates $ {\mathcal T}_{m}^{*}M $ as a vector space. \end{remark}

\begin{definition} Call two Poisson brackets $ \left\{,\right\}_{1} $ and $ \left\{,\right\}_{2} $ on $ M $ {\em compatible\/} if
the bracket $ \lambda_{1}\left\{,\right\}_{1} +\lambda_{2}\left\{,\right\}_{2} $ is Poisson for any $ \lambda_{1} $, $ \lambda_{2} $.

A {\em bihamiltonian structure\/} is a manifold $ M $ with a pair of compatible
Poisson brackets. \end{definition}

In fact it is possible to show that if {\em one\/} linear combination $ \lambda_{1}\left\{,\right\}_{1}
+\lambda_{2}\left\{,\right\}_{2} $ of two Poisson brackets is Poisson and $ \lambda_{1}\not=0 $, $ \lambda_{2}\not=0 $, then {\em any\/}
linear combination $ \lambda_{1}\left\{,\right\}_{1} +\lambda_{2}\left\{,\right\}_{2} $ is Poisson. In the analytic situation
the coefficients $ \lambda_{1} $, $ \lambda_{2} $ may be taken to be complex numbers.

If $ M $ is a $ C^{\infty} $-manifold with a bracket, we may consider the extension
of the bracket to the $ {\mathbb C} $-vector space of complex-valued functions on $ M $. In
this case $ \lambda_{1}\left\{,\right\}_{1} +\lambda_{2}\left\{,\right\}_{2} $ is well-defined even for complex values of
$ \lambda_{1},\lambda_{2} $. By the above remarks, complex linear combinations of brackets of a
bihamiltonian structure are also Poisson. In what follows we always
consider brackets as acting on the spaces of complex-valued functions.

\begin{definition} Given two brackets, $ \left\{\right\}_{M} $ on $ M $ and $ \left\{\right\}_{N} $ on $ N $, the {\em direct
product\/} of brackets $ \left\{\right\}_{M} $ and $ \left\{\right\}_{N} $ is the bracket on $ M\times N $ defined by
\begin{equation}
\left\{f_{M}\times f_{N},g_{M}\times g_{N}\right\}_{M\times N} \buildrel{\text{def}}\over{=} \left\{f_{M},g_{M}\right\}_{M}\times\left(f_{N}g_{N}\right)+\left(f_{M}g_{M}\right)\times\left\{f_{N},g_{N}\right\}_{N}.
\notag\end{equation}
Call a bihamiltonian structure {\em decomposable\/} if it isomorphic to a direct
product of two bihamiltonian structures of positive dimension. \end{definition}

Obviously, a direct product of two Poisson structures is a Poisson
structure, and a direct product of two bihamiltonian structures is a
bihamiltonian structure.

\begin{definition} \label{def002.40}\myLabel{def002.40}\relax  Consider a bihamiltonian structure $ \left(V,\left\{,\right\}_{1},\left\{,\right\}_{2}\right) $,
here $ V $ is a vector space. The bihamiltonian structure is
{\em translation-invariant\/} if $ \left\{{\mathfrak T}f,{\mathfrak T}g\right\}_{a}={\mathfrak T}\left\{f,g\right\}_{a} $, $ a=1,2 $, for any parallel
translation $ {\mathfrak T} $ on $ V $, any $ f $, and any $ g $. \end{definition}

\begin{definition} \label{def002.43}\myLabel{def002.43}\relax  A bihamiltonian structure on $ M $ is {\em flat\/} if it is
locally isomorphic to a translation-invariant bihamiltonian structure,
i.e., there is a collection of open subsets $ M_{i}\subset M $ such that $ M=\bigcup_{i\in I}M_{i} $, and
for any $ i\in I $ the restriction of the bihamiltonian structure on $ M $ to $ M_{i} $ is
isomorphic to an open subset $ \widetilde{M}_{i}\subset V_{i} $, here $ V_{i} $ is a vector space with a
translation-invariant bihamiltonian structure.

A bihamiltonian structure on $ M $ is {\em generically flat\/} if it is flat on
a dense open subset $ U\subset M $. \end{definition}

\begin{remark} Throughout the paper the phrase ``{\em at generic points\/}'' means ``at
points of an appropriate open dense subset''. Similarly, a ``{\em small open
subset\/}'' is used instead of ``an appropriate neighborhood of any given
point''. \end{remark}

\begin{remark} It is possible to give a complete classification of
translation-invariant bihamiltonian structures and a complete local
classification of flat bihamiltonian structures. (See Remark~\ref{rem6.13}.)
Classification of generically flat bihamiltonian structures is an
interesting unsolved problem which we do not consider in this paper. \end{remark}

\begin{remark} Any flat structure is generically flat, and any
translation-invariant structure is flat, but the opposite is not true. To
construct an example of non-translation-invariant flat structure one can
take a quotient of a translation-invariant structure on $ V $ by an arbitrary
discrete subgroup of $ V $. Later we will construct many generically flat
structures which are not flat. One of the simplest possible cases will be
provided in Example~\ref{ex002.45}, see also Theorems~\ref{th01.60},~\ref{th01.70}.

Not every bihamiltonian structure is generically flat. Important
examples of non-generically-flat structures will be constructed in
Section~\ref{h47}. \end{remark}

\begin{remark} The classification of Remark~\ref{rem6.13} shows that {\em indecomposable\/}
flat bihamiltonian structures break into two types with principally
different geometries: even-dimensional structures are modeled by Jordan
blocks, and odd-dimensional ones are modeled by Kronecker blocks. \end{remark}

Consider an interesting example of a translation-invariant
bihamiltonian structure. In fact it is going to be a key example of this
paper: we are going to show that this example is a ``building block'' in
decomposition of many ``classical'' examples of bihamiltonian structures.

\begin{example} Consider a vector space $ V $ with coordinates $ x_{0},\dots ,x_{2k-2} $ and
the Poisson brackets of coordinates
\begin{equation}
\left\{x_{2l},x_{2l+1}\right\}_{1}=1,\qquad \left\{x_{2l+1},x_{2l+2}\right\}_{2}=1,\qquad 0\leq l\leq k-2,
\label{equ45.20}\end{equation}\myLabel{equ45.20,}\relax 
any other brackets of coordinate functions $ x_{0},\dots ,x_{2k-2} $ vanishing. This
pair of brackets is in fact a translation-invariant bihamiltonian
structure. \end{example}

The following example is the simplest of classical examples of
bihamiltonian structures arising in theory of integrable systems.

\begin{example} \label{ex002.45}\myLabel{ex002.45}\relax  Given a Lie algebra $ {\mathfrak g} $ and an element $ \alpha\in{\mathfrak g}^{*} $, define a
bihamiltonian structure on $ {\mathfrak g}^{*} $ as in \cite{Bol91Com}. An element $ X\in{\mathfrak g} $ defines a
linear function $ f_{X} $ on $ {\mathfrak g}^{*} $. Due to Remark~\ref{rem002.20}, to define a
bihamiltonian structure on $ {\mathfrak g}^{*} $ it is enough to describe brackets $ \left\{f_{X},f_{Y}\right\}_{a} $,
$ a=1,2 $, $ X,Y\in{\mathfrak g} $.

Let $ \left\{f_{X},f_{Y}\right\}_{1} $ be a constant function on $ {\mathfrak g}^{*} $ and $ \left\{f_{X},f_{Y}\right\}_{2} $ be a linear
function on $ {\mathfrak g}^{*} $ given by the formulae
\begin{equation}
\left\{f_{X},f_{Y}\right\}_{1}\equiv c\left(X,Y\right)\buildrel{\text{def}}\over{=}f_{\left[X,Y\right]}\left(\alpha\right),\qquad \left\{f_{X},f_{Y}\right\}_{2}=f_{\left[X,Y\right]}.
\notag\end{equation}
The bracket $ \left\{,\right\}_{2} $ is the natural Lie--Kirillov--Kostant--Souriau Poisson
bracket on $ {\mathfrak g}^{*} $. The bracket $ \left\{,\right\}_{1} $ is translation-invariant. The bracket
$ \left\{,\right\}_{2} $ is translation-invariant only if $ {\mathfrak g} $ is abelian.

Call this bihamiltonian structure {\em regular\/} if $ {\mathfrak g} $ is semisimple and $ \alpha $
is regular semisimple. In such a case Conjecture~\ref{con01.100} states that
this structure is in fact generically flat (compare with \cite{Pan99Ver},
where a weaker property is proven\footnote{Paper \cite{Zakh99Kro} contains a proof of generic flatness of this
structure.}). In the case $ {\mathfrak g}={\mathfrak s}{\mathfrak l}_{2} $ the conjecture
follows from Theorem~\ref{th1.10}. This provides an example of generically
flat, but not flat and not translation-invariant structure.

In the case $ {\mathfrak g}={\mathfrak s}{\mathfrak l}_{2} $ it is easy to see that this structure is not flat.
Indeed, $ \left\{f,g\right\}_{2}|_{0}=0 $ for any $ f,g $. If the structure were flat, this would
imply $ \left\{f,g\right\}_{2}=0 $ for any $ f,g $, which is obviously false. \end{example}

By its definition, any flat bihamiltonian structure is locally
isomorphic to a direct product of several translation-invariant
indecomposable bihamiltonian structures. Introduce a special class of
bihamiltonian structures by allowing only special class of factors in the
above direct product.

\begin{definition} \label{def01.103}\myLabel{def01.103}\relax  A bihamiltonian structure is a {\em Kronecker\/} structure
if it is locally isomorphic to a direct product of several
translation-invariant odd-dimensional indecomposable structures. A {\em type\/}
of a Kronecker structure is the sequence of dimensions of factors in the
above direct product. The Kronecker structure is {\em indecomposable\/} if the
above product consists of one factor only.

A structure is {\em generically Kronecker\/} if it is Kronecker on an open
dense subset. \end{definition}

Note that a direct product of translation-invariant structures is
translation-invariant. In Section~\ref{h25} we will see that components of a
product of translation-invariant structures are uniquely determined by
the product. Thus Kroneker structures are flat structures open subsets of
which have {\em no even-dimensional\/} indecomposable components.

\begin{remark} The restriction of having no even-dimensional factors looks
very artificial. Moreover, one may think that bihamiltonian structures
which have {\em only\/} Jordan blocks should be the common case. Say, the
classification of even-dimensional bihamiltonian structures in general
position (\cite{Tur89Cla,Mag88Geo,Mag95Geo,McKeanPC,GelZakh93})
shows that on an open dense subset such pairs are isomorphic to direct
product of $ 2 $-dimensional bihamiltonian factors (thus have Jordan blocks
only in their decompositions). However, as we show later, some
``classical'' bihamiltonian systems are in fact generically Kronecker, and
we conjecture that many more such examples exist.

The condition of having no Jordan blocks is equivalent to the
condition of {\em completeness\/} of \cite{Bol91Com}. Note that the idea of the last
condition is to be one of possible {\em integrability criteria\/}: bihamiltonian
structures which are complete deserve to be called integrable. \end{remark}

By Remark~\ref{rem6.13}, flat bihamiltonian structures are essentially
pairs of skewsymmetric pairings on vector spaces, thus objects of linear
algebra. These objects of linear algebra have a classification, but the
building blocks of this classification are not only Jordan blocks, but
also some new blocks, constructed by Kronecker one year after Jordan.
This was the reason for our choice of the name.

\begin{remark} As Remark~\ref{rem6.13} will show, indecomposable odd-dimensional
flat bihamiltonian structures are locally isomorphic to the structure
given by~\eqref{equ45.20}. Thus the local geometry of a Kronecker structure is
uniquely determined by its type. \end{remark}

\begin{definition} \label{def01.120}\myLabel{def01.120}\relax  Consider a bracket $ \left\{,\right\} $ on a manifold $ M $. The
{\em associated bivector\/}\footnote{A {\em bivector field\/} is a skewsymmetric contravariant tensor of valence 2.} {\em field\/} $ \eta $ is the section of $ \Lambda^{2}{\mathcal T}M $ given by $ \left\{f,g\right\}|_{m}=\left<
\eta|_{m},df\wedge dg|_{m} \right> $, $ m\in M $, here $ \left<, \right> $ denotes the canonical pairing between
$ \Lambda^{2}{\mathcal T}_{m}M $ and $ \Omega_{m}^{2}M $. \end{definition}

\begin{definition} Consider a bracket $ \left\{,\right\} $ on $ M $ and $ m_{0}\in M $. The {\em associated pairing\/}
(,) in $ {\mathcal T}_{m_{0}}^{*}M $ is defined as $ \left(\alpha,\beta\right)=\left\{f,g\right\}|_{m_{0}} $ if $ \alpha=df|_{m_{0}} $, $ \beta=dg|_{m_{0}} $. \end{definition}

Obviously, the associated bivector field uniquely determines the
bracket and visa versa. The associated pairing is a skewsymmetric
bilinear pairing.

Given a pair of brackets $ \left\{,\right\}_{1} $ and $ \left\{,\right\}_{2} $, one obtains two bivector
fields $ \eta_{1} $, $ \eta_{2} $. Analogously, one obtains two skewsymmetric bilinear
pairings $ \left(,\right)_{1} $, $ \left(,\right)_{2} $ on $ {\mathcal T}_{m}^{*}M $, so that $ \left(\alpha,\beta\right)_{a}=\left\{f,g\right\}_{a}|_{m} $ if $ \alpha=df|_{m} $, $ \beta=dg|_{m} $,
$ a=1,2 $.

\begin{definition} The {\em rank\/} of the bracket $ \left\{,\right\} $ at $ m\in M $ is $ r $ if the associated
skewsymmetric bilinear pairing on $ {\mathcal T}_{m}^{*}M $ has rank $ r $. In this case the
{\em corank\/} of the bracket is $ \dim  M-r $.

A bracket has a {\em constant (co)rank\/} if its rank does not depend on the
point $ m\in M $. A bracket is {\em symplectic\/} if the corank is constant and equal to
0. \end{definition}

\begin{definition} Given a pair of vector spaces $ V^{\alpha} $ and $ V^{\beta} $, each equipped with
a pair of skewsymmetric bilinear pairings, equip $ V^{\alpha}\oplus V^{\beta} $ with two pairings
$ \left(,\right)_{a}\buildrel{\text{def}}\over{=}\left(,\right)_{a}^{\alpha}\oplus\left(,\right)_{a}^{\beta} $, $ a=1,2 $. If a pair is isomorphic to such a direct sum
with $ \dim  V^{i}\not=0 $, $ i=\alpha,\beta $, it is {\em decomposable}. \end{definition}

It is possible to provide a complete description of indecomposable
pairs of skewsymmetric pairings (we will do it in Theorem~\ref{th6.10}).

\begin{definition} \label{def01.105}\myLabel{def01.105}\relax  A bihamiltonian structure $ \left(M,\left\{\right\}_{1},\left\{\right\}_{2}\right) $ is {\em homogeneous\/}\footnote{A similar definition appears in \cite{Pan99Ver}.}
of type $ \left(2k_{1}-1,2k_{2}-1,\dots ,2k_{l}-1\right) $ if for any $ m\in M $ the pair of bilinear
pairings on $ {\mathcal T}_{m}^{*}M $ decomposes into a direct sum of indecomposable blocks of
dimensions $ 2k_{1}-1 $, $ 2k_{2}-1 $, \dots , $ 2k_{l}-1 $.

Such homogeneous system is {\em micro-indecomposable\/} if $ l=1 $. \end{definition}

By uniqueness of decomposition into indecomposable blocks (Theorem
~\ref{th6.10}), Kronecker structures are those bihamiltonian structures which
are simultaneously homogeneous and flat. There exist important examples
of homogeneous structures which are not flat (see Section~\ref{h47}).

What makes homogeneous structures important is the fact that the
standard algorithm of ``complete integration'' (so-called {\em anchored Lenard
scheme\/}) is applicable to these structures, and this algorithm provides
enough functions in involution for these structures only. (See Section
~\ref{h55} for details.)

In fact Kronecker structures are a {\em very special\/} case of homogeneous
structures:

\begin{conjecture} Given a sequence $ \left(2k_{1}-1,2k_{2}-1,\dots ,2k_{l}-1\right) $ there exist $ N>0 $
and a natural ways to assign tensor fields $ K_{1},\dots ,K_{N} $ to a homogeneous
bihamiltonian structure such that the structure is Kronecker iff $ K_{i}=0 $,
$ 1\leq i\leq N $. \end{conjecture}

In \cite{GelZakhWeb} we proved this conjecture in the case of
micro-indecomposable structures of dimension 3. This generalized to the
case of a general micro-indecomposable structure. In these cases $ N=1 $, and
the tensor field $ K_{1} $ is in fact a $ 2 $-form of curvature of a connection on
an appropriate line bundle (compare with \cite{Rig98Sys}). This $ 2 $-form plays
the same r\^ole for bihamiltonian structures as tensor of curvature plays
for Riemannian structures.

In what follows we provide criteria of homogeneity and of being an
indecomposable Kronecker structure. All these criteria are going to be
expressed in the following terms:

\begin{definition} Call a smooth function $ F $ on a manifold $ M $ with a Poisson
bracket $ \left\{,\right\} $ a {\em Casimir\/} function if $ \left\{F,f\right\}=0 $ for any smooth function $ f $ on $ M $.
\end{definition}

Obviously, any function $ \varphi\left(F_{1},F_{2},\dots ,F_{k}\right) $ of several Casimir functions
is again Casimir.

\begin{definition} A collection of smooth functions $ F_{1},\dots ,F_{r} $ on $ M $ is {\em dependent\/}
if $ \varphi\left(F_{1},\dots F_{r}\right)\equiv 0 $ for an appropriate smooth function $ \varphi\not\equiv 0 $. \end{definition}

We will use this definition when we want to pick up a small
independent collection of Casimir function out of the set of all Casimir
functions (possibly Casimir functions for several different brackets).

\section{Overview }\label{h003}\myLabel{h003}\relax 

One of the principal targets of this paper is to state three
criteria which for a given bihamiltonian structure determine whether it
is
\begin{enumerate}
\item
homogeneous micro-indecomposable structure (Theorem~\ref{th1.07});
\item
indecomposable Kronecker structure (Theorem~\ref{th1.10});
\item
homogeneous structure (Amplification~\ref{amp1.07}).
\end{enumerate}
We will use the criterion of Theorem~\ref{th1.10} to prove that open and
periodic {\em Toda lattices\/} are generically Kronecker (in Theorems~\ref{th01.60}
and~\ref{th01.70}), and to show that so-called {\em Lax structures\/} are
indecomposable Kronecker structures provided some conditions of general
position hold (in Theorem~\ref{th60.30}.)

The most interesting feature of all these criteria is that they are
stated in terms of {\em mutual position\/}\footnote{Given several functions $ \left\{F_{i}\right\}_{i\in I} $ on a manifold $ M $ and a point $ m_{0}\in M $,
consider the directions of differentials $ dF_{i}|_{m_{0}} $ of these functions at $ m_{0} $.
These directions can be considered as points of the projectivization
$ {\mathcal P}\left({\mathcal T}_{m_{0}}^{*}M\right) $ of the vector space $ {\mathcal T}_{m_{0}}^{*}M $. Thus we obtain a configuration of $ |I| $
points in a projective space, and this configuration depends on $ m_{0}\in M $. The
term ``mutual position'' refers to studying these configurations of points.} of Casimir functions for different
linear combinations $ \lambda_{1}\left\{,\right\}_{1} +\lambda_{2}\left\{,\right\}_{2} $ of Poisson brackets of the
bihamiltonian structure. We propose a way to encode these mutual
positions in a geometric structure of a new type, which we call a {\em web}.

Recall that the traditional Liouville approach to complete
integration of a dynamical system is to provide a system of so-called
{\em action-angle variables}. It so happens that in typical examples the
Casimir functions depend on action variables only. Moreover, the action
variables are typically much easier to find than the angle variable. This
indicates a fundamental asymmetry between action variables and angle
variables.

The notion of web (Definition~\ref{def02.20}) amplifies this asymmetry by
providing a way to remove angle variables from consideration whatsoever.
Since the Casimir functions do not depend on angle variables, it is
possible to study the mutual position of Casimir functions in terms of
the geometry of the web which corresponds to the given bihamiltonian
structure. Thus the conditions of the above criteria (of
being homogeneous or Kronecker structures) may be reformulated in terms
of webs.

The webs for micro-indecomposable bihamiltonian structures coincide
with {\em Veronese webs\/} which were studied\footnote{A beginning of a similar study in the case of general homogeneous
structures is done in \cite{Pan99Ver}.} in \cite{GelZakhWeb} and \cite{GelZakh93}.
After the criterion of being a Kronecker structure is reformulated as a
statement about webs, it becomes a direct corollary of results of
\cite{GelZakhWeb}. The results of \cite{GelZakhWeb} we use here only scratch the
surface of the beautiful theories of \cite{GelZakhWeb,GelZakh93,%
Pan99Ver}, in Section~\ref{h2} we provide an independent formulation of
these results, and prove the simplest of them. In Section~\ref{h45} we deduce
from these results the criterion~\ref{th1.10} of being an indecomposable
Kronecker structure.

Though the criteria~\ref{th1.07} and~\ref{amp1.07} of being a homogeneous
system may be formulated in terms of webs, in fact both the hypotheses
and the conclusions of these statements may be stated in terms of
individual cotangent spaces $ {\mathcal T}_{m}^{*}M $ to the bihamiltonian structure $ M $. Thus
these statements may be reduced to appropriate statements of linear
algebra. We do this reduction in Section~\ref{h25}.

The criterion~\ref{th1.10} of being an indecomposable Kronecker structure
is expressed in terms of several inequalities. In Sections~\ref{h47} and~\ref{h62}
we provide examples of bihamiltonian structures which show that no
inequality may be weakened without breaking the criterion. These examples
are homogeneous bihamiltonian structures which are not flat. One of these
examples shows that even a presence of a family of Casimir functions
which depend {\em polynomially\/} on a parameter does not guarantee flatness.

Note that all the examples of Section~\ref{h47} are completely
integrable. Here we use this vague term in the following sense: the
``anchored'' Lenard scheme works for these examples, and provides enough
functions in involution to construct action-angle variables. In Section
~\ref{h48} we describe the anchored Lenard scheme, and show its relations with
Casimir functions (thus with webs).

In Section~\ref{h55} we show that any homogeneous structure is completely
integrable via the anchored Lenard scheme. Theorem~\ref{th55.50} shows that in
fact the class of bihamiltonian structures which may be completely
integrated via the anchored Lenard scheme {\em coincides\/} with the class of
homogeneous structures. This answers a long-standing question in the
theory of integrable systems.

We finish the paper with applications of the criterion of flatness
to classical examples of integrable systems. After recalling (in Section
~\ref{h0}) definitions of {\em Toda lattices}, we show that the open and the periodic
Toda lattices are in fact generically flat (Theorems~\ref{th01.60} and
~\ref{th01.70}).

In Section~\ref{h60} we introduce a notion of a {\em Lax structure}. It is a
natural modification of the notion of Lax operator from \cite{KosMag96Lax}. We
show that under appropriate non-degeneracy conditions all the Lax
structures (in generic points) are indecomposable Kronecker structures.
In particular, two non-degenerate Lax structures of the same dimension
become isomorphic when restricted to appropriate open subsets.

Section~\ref{h005} contains conjectures which extend results of this
paper to the case of homogeneous systems which are not
micro-indecomposable.

\section{The principal criteria }\label{h1}\myLabel{h1}\relax 

One of the key ideas of this paper (compare with Conjecture
~\ref{con01.100}) is that many integrable systems admit a decomposition into a
product of ``simple'' bihamiltonian structures given by~\eqref{equ45.20}. Theorem
~\ref{th1.10} will provide an easy-to-check criterion when an open subset of a
given bihamiltonian structure is {\em isomorphic\/} to one given by~\eqref{equ45.20}.
Note that to check the criterion all one needs to know are Casimir
functions.

Note that a structure is locally isomorphic to one given by~\eqref{equ45.20}
iff it is an indecomposable Kronecker structure. In other words, it is
simultaneously a micro-indecomposable homogeneous structure, and a flat
structure. The following statement provides a criterion for the first
part, being a micro-indecomposable homogeneous structure.

\begin{theorem} \label{th1.07}\myLabel{th1.07}\relax  Consider a manifold $ M $, $ \dim  M\not=0 $, with two compatible
Poisson structures $ \left\{,\right\}_{1} $ and $ \left\{,\right\}_{2} $. Consider an open subset $ {\mathcal U}\subset{\mathbb R} $ and a
family of smooth functions $ F_{\lambda} $, $ \lambda\in{\mathcal U} $, on $ M $. Suppose that for any $ \lambda\in{\mathcal U} $ the
function $ F_{\lambda} $ is Casimir w.r.t.~the Poisson bracket $ \lambda\left\{,\right\}_{1}+\left\{,\right\}_{2} $, and that
$ dF_{\lambda}|_{m}\in{\mathcal T}_{m}^{*}M $ depends continuously on $ \lambda $ for any $ m\in M $. For $ m\in M $ denote by
$ W_{1}\left(m\right)\subset{\mathcal T}_{m}^{*}M $ the vector subspace spanned by the the differentials $ dF_{\lambda}|_{m} $ for
all possible $ \lambda\in{\mathcal U} $. If
\begin{enumerate}
\item
for one particular value $ m_{0}\in M $ one has $ \dim  W_{1}\left(m_{0}\right)\geq\frac{\dim  M}{2} $;
\item
for one particular value of $ \lambda_{1},\lambda_{2}\in{\mathbb R}^{2} $ the Poisson structure
$ \lambda_{1}\left\{,\right\}_{1}+\lambda_{2}\left\{,\right\}_{2} $ has at most one independent Casimir function on any open
subset of $ M $ near $ m_{0} $;
\end{enumerate}
then $ \dim  M $ is odd, and the bihamiltonian structure on $ M $ is
homogeneous of type $ \left(\dim  M\right) $ on an open subset $ U\subset M $ such that $ m_{0} $ is in
the closure of $ U $. \end{theorem}

The proof of this theorem is finished with the proof of Corollary
~\ref{cor25.40} in Section~\ref{h25}. Note that this proof implies also that $ \dim 
W_{1}\left(m_{0}\right)=\frac{\dim  M+1}{2} $. In fact the proof will show that if the Poisson
bracket $ \lambda_{1}\left\{,\right\}_{1}+\lambda_{2}\left\{,\right\}_{2} $ is of constant corank 1, then one may require that
$ m_{0}\in U $.

Amplification~\ref{amp1.07} provides a similar criterion of homogeneity
with an arbitrary type.

The following statement shows what one needs to know about
a micro-indecomposable homogeneous structure to ensure its flatness (thus
it being Kronecker):

\begin{theorem} \label{th1.10}\myLabel{th1.10}\relax  In addition to the conditions of Theorem~\ref{th1.07} suppose
that $ M $ is analytic, and $ F_{\lambda}\left(m\right) $ depends polynomially on $ \lambda $:
\begin{equation}
F_{\lambda}\left(m\right)=\sum_{k=0}^{d}f_{k}\left(m\right)\lambda^{k},
\notag\end{equation}
with analytic coefficients $ f_{k}\left(m\right) $ and the degree $ d $ satisfying $ d<\frac{\dim  M}{2} $.
Then the bihamiltonian structure on $ M $ is flat indecomposable of odd
dimension on an open subset $ U $ the closure of which contains $ m_{0} $. \end{theorem}

The proof of this theorem takes up to Section~\ref{h45}. Note that this
proof implies also that $ d=\frac{\dim  M -1}{2} $. Note that Conjecture~\ref{con01.120}
may provide a similar criterion applicable to arbitrary (i.e., not
necessarily indecomposable) Kronecker structures. The proof will actually
show the following statement (which cannot be expressed in terms of
Casimir functions only):

\begin{amplification} \label{amp1.12}\myLabel{amp1.12}\relax  In the case when in addition to conditions of
Theorem~\ref{th1.10} the Poisson structure $ \lambda_{1}\left\{,\right\}_{1}+\lambda_{2}\left\{,\right\}_{2} $ is of constant corank
1, the open subset $ U $ is in fact a neighborhood of $ m_{0} $. \end{amplification}

Remark~\ref{rem6.13} will show that all flat indecomposable structures of
dimension $ 2k-1 $ are locally isomorphic to each other, thus to the
structure given by~\eqref{equ45.20}. It is easy to see that for the structure of
~\eqref{equ45.20} one has $ \dim  M=2k-1 $, the vector space $ W_{1}\left(m\right) $ is spanned by $ dx_{0} $,
$ dx_{2},\dots $, $ dx_{2k-2} $, and the family $ F_{\lambda}\left(x\right) $ of degree $ k-1 $ is given by
~\eqref{equ45.25}.

\begin{remark} Not all homogeneous bihamiltonian structures of type $ \left(2k-1\right) $ are
flat, as the examples of Section~\ref{h47} show (already in the case $ k=2 $).

The example of $ \left\{,\right\}_{1}=\left\{,\right\}_{2}\equiv 0 $ shows that in Theorem~\ref{th1.07} one cannot
drop the restriction on the number of independent Casimir functions.
Considering a direct product of $ M $ with any bihamiltonian structure shows
the significance of the bound on $ \dim  W_{1} $. Moreover, Proposition~\ref{prop47.40}
implies that one cannot weaken the bound $ d<\frac{\dim  M}{2} $ of Theorem~\ref{th1.10}.
\end{remark}

\begin{remark} As Theorem~\ref{th01.60} will show, one can also consider Theorem
~\ref{th1.10} as a criterion that a given bihamiltonian structure is locally
isomorphic to an open subset of the open Toda lattice. \end{remark}

\begin{remark} Theorems~\ref{th1.07} and~\ref{th1.10} are almost immediate corollaries
of results of \cite{GelZakhWeb} and \cite{GelZakh93}. However, since we will need
many results of these papers anyway, the following three sections provide
almost self-contained proof of these theorems. The only component of the
proof which requires a reference to \cite{GelZakhWeb} is the last statement of
Theorem~\ref{th2.07}. The proof of this statement is outside of the scope of
this paper (compare with Remark~\ref{rem2.017}). \end{remark}

\section{Linear case and criterion of homogeneity }\label{h25}\myLabel{h25}\relax 

Recall the classification of pairs of skewsymmetric bilinear
pairings from \cite{GelZakhFAN} (see also \cite{GelZakhWeb,GelZakh93}). For $ k\in{\mathbb N} $
consider the identity $ k\times k $ matrix $ I_{k} $. For $ \mu\in{\mathbb C} $ consider the Jordan block
$ J_{k,\mu} $ of size $ k $ and eigenvalue $ \mu $. The pair of matrices
\begin{equation}
{\text H}_{1}^{\left(\mu\right)}= \left( 
\begin{matrix}
0 & J_{k,\mu}
\\
-J_{k,\mu}^{t} & 0
\end{matrix}
\right),\qquad {\text H}_{2}^{\left(\mu\right)}=\left( 
\begin{matrix}
0 & I_{k}
\\
-I_{k} & 0
\end{matrix}
\right)
\notag\end{equation}
defines a pair of skewsymmetric bilinear pairings on vector space $ {\mathbb C}^{2k} $. The
limit case of $ \mu \to \infty $ may be deformed to
\begin{equation}
{\text H}_{1}^{\left(\infty\right)}= \left( 
\begin{matrix}
0 & I_{k}
\\
-I_{k} & 0
\end{matrix}
\right),\qquad {\text H}_{2}^{\left(\infty\right)}= \left( 
\begin{matrix}
0 & J_{k,0}
\\
-J_{k,0}^{t} & 0
\end{matrix}
\right).
\notag\end{equation}
Denote the pair $ \left({\text H}_{1}^{\left(\mu\right)},{\text H}_{2}^{\left(\mu\right)}\right) $ of skewsymmetric bilinear pairings by
$ {\mathcal J}_{2k,\mu} $, $ k\in{\mathbb N} $, $ \mu\in{\mathbb C}{\mathbb P}^{1} $.

Add to this list the so-called Kroneker pair $ {\mathcal K}_{2k-1} $. This is a pair
in a vector space $ {\mathbb C}^{2k-1} $ with a basis $ \left({\mathbit w}_{0},{\mathbit w}_{1},\dots ,{\mathbit w}_{2k-2}\right) $. The only non-zero
pairings are
\begin{equation}
\left({\mathbit w}_{2l},{\mathbit w}_{2l+1}\right)_{1}=1,\qquad \left({\mathbit w}_{2l+1},{\mathbit w}_{2l+2}\right)_{2}=1,
\label{equ2.10}\end{equation}\myLabel{equ2.10,}\relax 
for $ 0\leq l\leq k-2 $. Obviously, different pairs from this list are not isomorphic.

\begin{theorem} \label{th6.10}\myLabel{th6.10}\relax  (\cite{GelZakhFAN,Thom91Pen}) Any pair of skewsymmetric
bilinear pairings on a finite-dimensional complex vector space can be
decomposed into a direct sum of pairs of the pairings isomorphic to
$ {\mathcal J}_{2k,\mu} $, $ k\in{\mathbb N} $, $ \mu\in{\mathbb P}^{1} $, and $ {\mathcal K}_{2k-1} $, $ k\in{\mathbb N} $. The types of the components of this
decomposition are uniquely determined. \end{theorem}

Though this simple statement was known for a long time (say, the
preprint of \cite{Thom91Pen} existed in 1973), we do not know whether it was
published before it was used in \cite{GelZakhFAN}. The discussions in
\cite{Gan59The} and \cite{TurAith61Int} come very close, but do not state this
result.

\begin{remark} The papers \cite{GelZakhFAN,GelZakh94Spe} described significance
of Kronecker blocks in the spectral theory of pencils $ A_{\lambda}=A+\lambda B $, $ \lambda\in{\mathbb C} $, of
differential operators. Though it is not used in this paper, let us
highlight the details of this description.

The Jordan blocks which appear in spectral theory of pencils
correspond to values of $ \lambda $ where the dimension of $ \operatorname{Ker} A_{\lambda} $ jumps up. It so
happens that due to special properties of the pencil $ A_{\lambda} $ (say, skew
symmetry of operators) it may happen that $ \operatorname{Ker} A_{\lambda}\not=0 $ for any $ \lambda $ (this is
what actually happens in the pencil related to the periodic case of KdV
equation). In such a case the direct sum of Jordan blocks has a
non-trivial complement in the vector space where the pencil acts.

For so-called {\em finite gap potentials\/} this {\em defect space\/} happens to be
exactly the Kronecker block $ {\mathcal K}_{2k-1} $ (here $ k $ is the number of gaps), thus the
situation is absolutely parallel to the finite-dimensional case discussed
above. In the case of infinitely many gaps an appropriate
infinite-dimensional analogue of Kronecker blocks may be described.

Note, however, that it is absolutely unclear how to translate this
description of the linear situation (which is associated to one cotangent
space to the phase space of KdV) to the nonlinear bihamiltonian geometry
of KdV. While results and conjectures of this paper illuminate the
bihamiltonian geometry of finite-dimensional systems in many details,
they do not look applicable in infinite-dimensional situation.

The main obstruction is that while all the Kronecker blocks of the
same dimension are isomorphic, infinite-dimensional Kronecker blocks
acquire new invariants---{\em fuzzy eigenvalues}. Though fuzzy, these data in
fact completely disambiguate points which may be distinguished by Casimir
functions (at least for real-analytic potentials, for details see
\cite{GelZakhFAN}).

One can see that the linearized geometry of periodic KdV is very
similar to geometry on odd-dimensional manifolds---there is exactly one
Kronecker block, the rest is Jordan blocks with $ k=1 $, and in generic
points there is no Jordan block. But the non-linear geometry of KdV is in
some regards also similar to even-dimensional geometry in the sense that
the points $ m_{1},m_{2}\in M $ which are separated by Casimir functions also have
non-isomorphic pairings in $ {\mathcal T}_{m_{1}}^{*}M $, $ {\mathcal T}_{m_{2}}^{*}M $. \end{remark}

\begin{remark} \label{rem6.13}\myLabel{rem6.13}\relax  Given a skewsymmetric bilinear pairing (,) on a vector
space $ V^{*} $, consider the bracket $ \left\{,\right\} $ on the vector space $ V $ described by
$ \left\{f,g\right\}|_{m}=\left(df|_{m},dg|_{m}\right) $. As it is easy to check, this bracket is
translation-invariant and Poisson. Given a pair of such pairings $ \left(,\right)_{1} $,
$ \left(,\right)_{2} $ on $ V^{*} $ one obtains a translation-invariant bihamiltonian structure on
$ V $. Obviously, any translation-invariant bihamiltonian structure may be
obtained this way.

Similarly, any decomposable flat bihamiltonian structure is locally
isomorphic to a product of two flat bihamiltonian structures. Indeed, it
is enough to show that if an open subset $ U $ of the above bihamiltonian
structure on $ V $ is decomposable, then the pair of pairings on $ V^{*} $ is
decomposable, which is obvious.

Thus Theorem~\ref{th6.10} gives also a complete classification of
translation-invariant bihamiltonian structures, a complete local
classification of flat bihamiltonian structures, and a description of
indecomposable flat structures. \end{remark}

For the topics we discuss here it is not necessary to answer the
following question, but it is interesting nevertheless:

\begin{conjecture} Consider two bihamiltonian structures on $ M_{1} $ and $ M_{2} $.
Suppose that $ M_{1}\times M_{2} $ is flat. Then $ M_{1} $ and $ M_{2} $ are flat. \end{conjecture}

The first step in the proof of Theorem~\ref{th1.07} is the following

\begin{proposition} \label{prop6.15}\myLabel{prop6.15}\relax  Consider a pair of skewsymmetric bilinear pairings
$ \left(,\right)_{1} $, $ \left(,\right)_{2} $ on a finite-dimensional complex vector space $ W $. Suppose there
is a finite set $ L $ and there are families of vectors $ w_{l,\lambda}\in W $, $ l\in L $,
polynomially depending on $ \lambda $ such that $ \lambda\left(w_{l,\lambda},w\right)_{1}+\left(w_{l,\lambda},w\right)_{2}=0 $ for any $ w\in W $,
$ l\in L $, and $ \lambda\in{\mathbb C} $. Denote by $ W_{1} $ the vector subspace spanned by $ w_{l,\lambda} $, $ l\in L $, $ \lambda\in{\mathbb C} $.
Suppose that for one particular value of $ \lambda_{1} $, $ \lambda_{2} $ the corank of the
bilinear pairing $ \lambda_{1}\left(,\right)_{1}+\lambda_{2}\left(,\right)_{2} $ is $ r $. If $ \dim  W_{1} \geq \frac{\dim  W+r-1}{2} $, then the
pair $ \left(,\right)_{1} $, $ \left(,\right)_{2} $ is isomorphic to $ \oplus_{t=1}^{r}{\mathcal K}_{2k_{t}-1} $ with $ \sum_{t}k_{t}=\dim  W_{1} $. In
particular, $ \dim  W_{1} = \frac{\dim  W+r}{2} $. \end{proposition}

\begin{proof} We may assume that the pair $ \left(,\right)_{1} $, $ \left(,\right)_{2} $ is a direct sum of
several blocks of the form $ {\mathcal J}_{2k,\mu} $ and $ {\mathcal K}_{2k-1} $, and that for any $ l\in L $ the
family $ w_{l,\lambda}\not\equiv 0 $. We suppose that $ \left(\lambda_{1},\lambda_{2}\right)\not=\left(0,0\right) $, it is easy to consider the
remaining case separately.

Start with supposing that there are only blocks of the form $ {\mathcal K}_{2k_{t}-1} $,
$ t=1,\dots ,T $. Then the only things we need to prove is that $ T=r $, and $ \dim  W_{1}
\leq\sum_{t}k_{t} $. The first statement is obvious.

The following lemma follows immediately from the explicit
description of the pair $ {\mathcal K}_{2k-1} $:

\begin{lemma} \label{lm25.20}\myLabel{lm25.20}\relax  For the pair $ {\mathcal K}_{2k-1} $ of skewsymmetric pairings there exists a
family of vectors $ \widetilde{w}_{\lambda}\in W $ polynomially depending on $ \lambda $ such that
$ \lambda\left(\widetilde{w}_{\lambda},w\right)_{1}+\left(\widetilde{w}_{\lambda},w\right)_{2}=0 $ for any $ w\in W $ and $ \lambda\in{\mathbb C} $, and the degree of $ \widetilde{w}_{\lambda} $ in $ \lambda $ is $ k-1 $.
This family is defined uniquely up to multiplication by a constant, and
it spans a $ k $-dimensional vector subspace. Any other polynomial family $ w_{\lambda} $
such that $ \lambda\left(w_{\lambda},w\right)_{1}+\left(w_{\lambda},w\right)_{2}=0 $ for any $ w\in W $ and $ \lambda\in{\mathbb C} $ may be written as $ p\left(\lambda\right)\widetilde{w}_{\lambda} $
for an appropriate scalar polynomial $ p $. \end{lemma}

Denote the family $ \widetilde{w}_{\lambda} $ for the Kronecker block $ {\mathcal J}_{2k_{t}-1} $ by $ \widetilde{w}_{\lambda}^{\left(t\right)} $. Due to
this lemma one can write $ w_{l,\lambda}=\sum_{t=1}^{T}p_{lt}\left(\lambda\right)\widetilde{w}_{\lambda}^{\left(t\right)} $, thus $ \dim  W_{1}\leq\sum_{t=1}^{r}k_{t}=\frac{\dim 
W +r}{2} $. Since $ \dim  W + r $ is even, this shows that $ \dim  W_{1} = \frac{\dim  W+r}{2} $,
thus finishes proof of the proposition in the case when there are no
Jordan blocks.

Consider now the general case. First of all, $ w_{l,\lambda}\not=0 $ for a generic $ \lambda $,
thus $ w_{l,\lambda} $ (for a generic $ \lambda $) is in the null-space of the linear
combination $ \lambda\left(,\right)_{1}+\left(,\right)_{2} $. Since for a block of the form $ {\mathcal J}_{2k,\mu} $ and generic $ \lambda $
this combination has no null-space, it is obvious that $ w_{l,\lambda} $ is in the sum
of components of the form $ {\mathcal K}_{2k-1} $. Since removing a component of the form
$ {\mathcal J}_{2k,\mu} $ decreases $ \dim  W $ by $ 2k $, does not change $ \dim  W_{1} $, and may only
decrease $ r $, one can see that conditions of the proposition are applicable
to the sum of components of the form $ {\mathcal K}_{2k-1} $, but the equality on $ \dim  W_{1} $ is
sharpened by at least $ k $. However, we have seen that it is not possible to
sharpen this inequality more than by $ \frac{1}{2} $, which proves that $ W $
contains no Jordan components. \end{proof}

\begin{amplification} In Lemma~\ref{lm25.20} and Proposition~\ref{prop6.15} one may
suppose (without changing the conclusions\footnote{With an obvious exception that $ p $ in Lemma~\ref{lm25.20} becomes a continuous
function.} of these statements) that
families $ w_{l,\lambda} $ are continuous functions of $ \lambda $ defined on a given open
subset $ {\mathcal U}\subset{\mathbb C} $ or $ {\mathcal U}\subset{\mathbb R} $. \end{amplification}

\begin{corollary} \label{cor25.40}\myLabel{cor25.40}\relax  In conditions of Theorem~\ref{th1.07} the dimension of $ M $
is odd. There is a point $ m_{1} $ of $ M $ such that the pair of skewsymmetric
bilinear pairings in $ {\mathcal T}_{m_{1}}^{*}M $ is isomorphic to $ {\mathcal K}_{2k-1} $ with $ \dim  M=2k-1 $. \end{corollary}

\begin{proof} In this prove we assume that $ M $ is a complex manifold, so that
$ {\mathcal T}_{m}^{*}M $ is a complex vector space for any $ m\in M $. If $ M $ is a $ C^{\infty} $-manifold, one
should substitute $ {\mathcal T}_{m}^{*}M\otimes{\mathbb C} $ instead of $ {\mathcal T}_{m}^{*}M $ in the arguments below.

In conditions of Theorem~\ref{th1.07} if $ m_{1} $ in a neighborhood $ \widetilde{U} $ of
the point $ m_{0}\in M $, then vectors $ dF_{\lambda}|_{m_{1}}\in{\mathcal T}_{m_{1}}^{*}M $ span a vector subspace $ W_{1}\left(m_{1}\right) $
satisfying $ \dim  W_{1}\left(m_{1}\right)>\frac{\dim  M}{2} $. There is an open subset $ U_{r}\subset\widetilde{U} $ where
$ \lambda_{1}\left\{,\right\}_{1}+\lambda_{2}\left\{,\right\}_{2} $ has a constant corank $ r $. Obviously, there is $ r\in{\mathbb Z} $ such that
the point $ m_{0} $ is in the closure of $ U_{r} $. Restrict our attention to this
value of $ r $. Let $ m_{1} $ be in $ U_{r} $, and $ W={\mathcal T}_{m_{1}}^{*}M $, $ L=\left\{\bullet\right\} $, and $ w_{\bullet,\lambda}=dF_{\lambda}|_{m_{1}} $. Then
the span $ W_{1} $ of vectors $ w_{\bullet,\lambda} $ considered for all possible $ \lambda\in{\mathcal U} $ satisfies $ \dim 
W_{1}>\frac{\dim  W}{2} $, thus $ \dim  W_{1}\geq\frac{\dim  W+1}{2} $.

The vector space $ W $ is equipped with two skewsymmetric bilinear
pairings $ \left(,\right)_{1} $, $ \left(,\right)_{2} $ given by values of $ \eta_{1} $, $ \eta_{2} $ (see Definition
~\ref{def01.120}) at $ m_{1} $. Obviously, $ w_{\bullet,\lambda} $ is in the kernel of $ \lambda\left(,\right)_{1}+\left(,\right)_{2} $.

By the conditions of Theorem~\ref{th1.07}, there is at most one
independent Casimir function near $ m_{1} $, thus $ r\leq1 $. Obviously, this is the
same $ r $ as in Proposition~\ref{prop6.15}, thus $ \dim  M\not=0 $ implies $ r\not=0 $. Hence the
pair $ \left(,\right)_{1} $, $ \left(,\right)_{2} $ is isomorphic to $ {\mathcal K}_{2k-1} $ for an appropriate $ k $, thus $ \dim  M $
is odd. \end{proof}

This proves Theorem~\ref{th1.07}. In Section~\ref{h2} we show that it also
allows one to apply the results of \cite{GelZakhWeb,GelZakh93} to prove
Theorem~\ref{th1.10} as well.

Corollary~\ref{cor25.40} uses a particular case of Proposition~\ref{prop6.15}
with $ r=1 $. While we will not need it in this paper, it is possible to
strengthen Corollary~\ref{cor25.40} so that it uses the full power of
Proposition~\ref{prop6.15}. This result would move us one step in the
direction of Conjecture~\ref{con01.120}.

\begin{amplification} \label{amp1.07}\myLabel{amp1.07}\relax  Consider a manifold $ M $ with two compatible Poisson
structures $ \left\{,\right\}_{1} $ and $ \left\{,\right\}_{2} $. Consider a finite set $ L $, open subsets $ {\mathcal U}_{l}\subset{\mathbb C} $,
$ l\in L $, and families of smooth functions $ F_{l,\lambda} $, $ l\in L $, $ \lambda\in{\mathcal U}_{l} $, on $ M $. Suppose that
for any $ l\in L $ and any $ \lambda\in{\mathcal U}_{l} $ the function $ F_{l,\lambda} $ is Casimir w.r.t.~the Poisson
bracket $ \lambda\left\{,\right\}_{1}+\left\{,\right\}_{2} $, and that $ dF_{l,\lambda}|_{m}\in{\mathcal T}_{m}^{*}M $ depends continuously on $ \lambda $ for
any $ l\in L $ and $ m\in M $. For $ m\in M $ denote by $ W_{1}\left(m\right)\subset{\mathcal T}_{m}^{*}M $ the vector subspace spanned
by the the differentials $ dF_{l,\lambda}|_{m} $ for all possible $ l $ and $ \lambda\in{\mathcal U}_{l} $. If for an
appropriate $ R\in{\mathbb Z}_{\geq0} $
\begin{enumerate}
\item
for one particular value $ m_{0}\in M $ one has $ \dim  W_{1}\left(m_{0}\right)\geq\frac{\dim  M+R}{2} $;
\item
for one particular value of $ \lambda_{1},\lambda_{2}\in{\mathbb C}^{2} $ the Poisson structure
$ \lambda_{1}\left\{,\right\}_{1}+\lambda_{2}\left\{,\right\}_{2} $ has at most $ R $ independent Casimir functions on any open
subset of $ M $ near $ m_{0} $;
\end{enumerate}
then $ \dim  M-R $ is even, $ \dim  W_{1}\left(m_{0}\right)=\frac{\dim  M+R}{2} $, and the bihamiltonian
structure on $ M $ is homogeneous of type $ \left(t_{1},\dots ,t_{R}\right) $ on an open subset $ U\subset M $
such that $ m_{0} $ is in the closure of $ U $. Here $ t_{k}\in{\mathbb Z}_{>0} $ are appropriate numbers
with $ \sum_{k}t_{k}=\dim  M $. \end{amplification}

\begin{proof} First of all, one can proceed as in Corollary~\ref{cor25.40} up to
the moment we concluded $ r\leq1 $. Under the conditions of the amplification we
conclude that $ r\leq R $, thus $ \dim  W_{1}\left(m_{1}\right)\geq\frac{\dim  W+r}{2} $. Proposition~\ref{prop6.15}
implies that $ \dim  W_{1}\left(m_{1}\right)=\frac{\dim  M+r}{2} $, thus $ r\geq R $. This shows that in fact
$ r=R $.

We can conclude that for $ m $ in an appropriate open subset $ U\subset M $ the pair
of bilinear pairings on the vector space $ {\mathcal T}_{m}^{*}M $ is isomorphic to a direct sum
of $ R $ Kronecker blocks. What remains to prove is that the dimensions of
these blocks do not depend on $ m $ in an appropriate open subset of $ U $.

Fix a vector space $ V $. For a sequence $ T=\left(t_{1}\leq\dots \leq t_{R}\right) $ denote by
$ {\mathfrak F}_{T}\subset\Lambda^{2}V^{*}\times\Lambda^{2}V^{*} $ the set of pairs of skewsymmetric bilinear pairings which are
isomorphic to $ \bigoplus_{a=k}^{R}{\mathcal K}_{t_{k}} $. In particular, $ {\mathfrak F}_{T} $ is not empty iff all $ t_{k} $ are
odd and $ \sum t_{k}=\dim  V $. Moreover, $ {\mathfrak F}_{T} $ is a $ \operatorname{GL}\left(V\right) $-orbit.

It follows that if $ {\mathfrak F}_{T'} $ intersects the closure of $ {\mathfrak F}_{T} $, then $ {\mathfrak F}_{T'} $ is
contained in this closure. Fix a neighborhood $ U_{1} $ of $ m_{0} $, let $ T^{\left(1\right)},\dots ,T^{\left(N\right)} $
be such sequences that there are points $ m $ in $ U\cap U_{1} $ where the pair of
pairings is in each of $ {\mathfrak F}_{T^{\left(k\right)}} $, $ 1\leq k\leq N $. Suppose that $ {\mathfrak F}_{T^{\left(1\right)}},\dots ,{\mathfrak F}_{T^{\left(M\right)}} $ are of
maximal possible dimension among $ {\mathfrak F}_{T^{\left(1\right)}},\dots ,{\mathfrak F}_{T^{\left(N\right)}} $, then the points $ m $ in
$ U\cap U_{1} $ where the pair of pairings is in any one of $ {\mathfrak F}_{T^{\left(k\right)}} $, $ 1\leq k\leq M $, form an open
subset. Obviously, at least one of these subsets has $ m_{0} $ in its closure. \end{proof}

\begin{remark} It is not clear whether one can improve the statement of
Amplification~\ref{amp1.07} provided that the rank of $ \lambda_{1}\left\{,\right\}_{1}+\lambda_{2}\left\{,\right\}_{2} $ is
constant near $ m_{0} $. Recall that in Theorem~\ref{th1.07} one {\em could\/} conclude that
the structure is homogeneous in a {\em neighborhood\/} of $ m_{0} $. However, under the
condition of constant rank one can weaken the condition on dimension to
become $ \dim  W_{1}\left(m_{0}\right)\geq\frac{\dim  M+R-1}{2} $.

To recognize a possibility of a jump of the type of decomposition of
$ {\mathcal T}_{m}^{*}M $, consider the vector space with a basis $ {\mathbit w}_{0},\dots ,{\mathbit w}_{4},{\mathbit W} $ with the only
non-zero pairings being
\begin{equation}
\left({\mathbit w}_{2l},{\mathbit w}_{2l+1}\right)_{1}=1,\qquad \left({\mathbit w}_{2l+1},{\mathbit w}_{2l+2}\right)_{2}=1,
\notag\end{equation}
for $ 0\leq l\leq1 $, and $ \left({\mathbit W},{\mathbit w}_{1}\right)=\left({\mathbit W},{\mathbit w}_{3}\right)=\varepsilon $. If $ \varepsilon\not=0 $, then this pair is of the type
$ {\mathcal K}_{3}\oplus{\mathcal K}_{3} $, if $ \varepsilon=0 $, it is of the type $ {\mathcal K}_{5}\oplus{\mathcal K}_{1} $. Thus different orbits $ {\mathfrak F}_{T} $ may be
adjacent\footnote{The recent preprint \cite{Pan99Ver} contains example of a {\em bihamiltonian
structure\/} where such an adjacency takes place. Thus the statement of
Amplification~\ref{amp1.07} cannot be improved.} indeed. \end{remark}

\section{Bihamiltonian structures and webs }\label{h02}\myLabel{h02}\relax 

Consider a manifold $ M $ with a Poisson bracket $ \left\{,\right\} $. To define the
notion of a symplectic leaf on $ M $, consider Casimir functions on $ M $. The
local classification of Poisson structures of {\em constant rank\/} \cite{Kir76Loc,%
Wei83Loc} shows that for an arbitrary Poisson bracket there is an open
(and in interesting cases dense) subset $ U\subset M $ and $ k\in{\mathbb Z}_{\geq0} $ such that on $ U $
there are $ k $ independent Casimir functions $ F_{1},\dots ,F_{k} $, and any Casimir
function on $ U $ may be written as a function of $ F_{1},\dots ,F_{k} $ (we do not
exclude the case $ k=0 $). The common level sets $ F_{1}=C_{1},\dots ,F_{k}=C_{k} $ form an
invariantly defined foliation on $ U $, which is called the {\em symplectic
foliation}. Note that one can define this foliation as an equivalence
relation given by $ m_{1}\sim m_{2} $ iff $ F\left(m_{1}\right)=F\left(m_{2}\right) $ for any Casimir function $ F $ on $ U $.

Consider now a pair $ \left\{,\right\}_{1} $, $ \left\{,\right\}_{2} $ of compatible Poisson structures on $ M $
(i.e., a bihamiltonian structure). Proceed as with the above construction
of leaves, and consider

\begin{definition} A smooth function $ F $ on $ M $ is {\em semi-Casimir\/} if there is
$ \left(\lambda_{1},\lambda_{2}\right)\not=\left(0,0\right) $ such that $ F $ is a Casimir function for $ \lambda_{1}\left\{,\right\}_{1}+\lambda_{2}\left\{,\right\}_{2} $. \end{definition}

For any open subset $ U\subset M $ define an equivalence relation on $ U $ by $ m_{1}\sim m_{2} $
iff $ F\left(m_{1}\right)=F\left(m_{2}\right) $ for any semi-Casimir function $ F $ on $ U $. Denote by $ {\mathcal B}_{U} $ the
topological space of equivalence classes. Then any semi-Casimir function
$ F $ on $ U $ induces a continuous function on $ {\mathcal B}_{U} $. Any function on $ {\mathcal B}_{U} $ induces a
pull-back function on $ U $.

As a result, to any local bihamiltonian structure $ \left(U,\left\{,\right\}_{1},\left\{,\right\}_{2}\right) $ we
associated a topological space $ {\mathcal B}_{U} $. Let $ \lambda=\left(\lambda_{1}:\lambda_{2}\right)\in{\mathbb C}{\mathbb P}^{1} $, let $ {\mathfrak C}_{\lambda} $ be the
vector space of functions on $ {\mathcal B}_{U} $ pull-backs of which are Casimir functions
for $ \lambda_{1}\left\{,\right\}_{1}+\lambda_{2}\left\{,\right\}_{2} $. Note that $ \varphi\left(F_{1},F_{2},\dots ,F_{k}\right)\in{\mathfrak C}_{\lambda} $ if $ \varphi $ is smooth and
$ F_{1},F_{2},\dots ,F_{k}\in{\mathfrak C}_{\lambda} $. This allows one to consider $ {\mathfrak C}_{\lambda} $ as a $ C^{0} $-analogue of a
set of local equations of a foliation.

Later we will see that in the cases we study here $ {\mathcal B}_{U} $ is a
manifold, and for any $ \lambda $ the space $ {\mathfrak C}_{\lambda} $ is the set of local equations of a
foliation on $ {\mathcal B}_{U} $. The codimension of this foliation is not going to depend
on $ \lambda\in{\mathbb C}{\mathbb P}^{1} $. Anyway, we come to

\begin{definition} \label{def02.20}\myLabel{def02.20}\relax  A {\em web\/}\footnote{The reason for this name is that $ {\mathcal B} $ is equipped with a huge family of
canonically defined subsets: for any $ \lambda $ one consider intersections of
level sets of functions from $ {\mathfrak C}_{\lambda} $. Moreover, one can consider intersections
of such subsets for different values of $ \lambda $. If one assumes that $ {\mathcal B} $ and
these intersections are manifolds, then one gets a delicate network of
submanifolds, with infinitely many of them passing through each given
point $ b\in{\mathcal B} $.} is a topological space $ {\mathcal B} $ with a given subset
$ {\mathfrak C}_{\lambda} $ of the set of continuous functions on $ {\mathcal B} $ for any $ \lambda\in{\mathbb C}{\mathbb P}^{1} $. We require that
$ \varphi\left(F_{1},F_{2},\dots ,F_{k}\right)\in{\mathfrak C}_{\lambda} $ if $ \varphi $ is smooth and $ F_{1},F_{2},\dots ,F_{k}\in{\mathfrak C}_{\lambda} $. \end{definition}

One can also introduce a notion of $ {\mathcal U} $-{\em web\/} for any subset $ {\mathcal U}\subset{\mathbb C}{\mathbb P}^{1} $, the
only change being that $ \lambda\in{\mathcal U} $ instead of $ \lambda\in{\mathbb C}{\mathbb P}^{1} $.

\begin{proposition} To any bihamiltonian structure $ \left(M,\left\{\right\}_{1},\left\{\right\}_{2}\right) $ one can
associate a structure of a web on $ {\mathcal B}_{M}=M/\sim $. \end{proposition}

In \cite{GelZakhWeb} and \cite{GelZakh93} it was shown that in some
particularly interesting types of bihamiltonian structures the class of
the web $ {\mathcal B}_{U} $ up to an isomorphism determines the class of bihamiltonian
structure on $ U $ up to an isomorphism (compare Theorem~\ref{th2.07}), at least
for small open subsets $ U\subset M $. This is going to be the main instrument used
in this paper: we show that the bihamiltonian structure from Theorem
~\ref{th1.10} and the structure given by~\eqref{equ45.20} are of the type mentioned
above, and show that the corresponding webs are locally isomorphic. This
will imply a local isomorphism of bihamiltonian structures.

To illustrate advantages of the approach of \cite{GelZakhWeb} and
\cite{GelZakh93} introduce

\begin{definition} A smooth function $ F $ on $ M $ is an {\em action\/} function if locally on
each small open subset $ U\subset M $ it is a pull-back from a function on $ {\mathcal B}_{U} $. \end{definition}

Obviously, any function of the form $ \varphi\left(F_{1},\dots ,F_{l}\right) $ with semi-Casimir
functions $ F_{1},\dots ,F_{l} $ (not necessarily corresponding to the same $ \lambda $) is an
action function. (The name is related to the fact that in bihamiltonian
geometry {\em action-\/} and {\em angle-variables\/} may be defined by local means.
Action function are functions of action variables.)

In these terms the approach of \cite{GelZakhWeb} and \cite{GelZakh93} states
that to construct an isomorphism of bihamiltonian structures $ M' $ and $ M'' $
it is enough to associate to each action function on $ M' $ an action
function on $ M'' $ (with appropriate compatibilities conditions this is
equivalent to constructing a diffeomorphism of the webs). One needs not
care about ``angle'' variables. Since explicit constructions of ``angle''
variables is the most complicated step of integration of a dynamical
system, this leads to very significant simplifications.

In particular, we are going to construct an isomorphism of manifolds
of (approximately) half the dimension of the initial manifolds. Moreover,
these smaller manifolds have a very {\em rigid\/} geometric structure\footnote{In particular, it has at most $ 1 $-dimensional group of automorphisms
preserving a given point, as opposed to the group of automorphisms of
bihamiltonian structures themselves. Recall that in \cite{GelZakhWeb} and
\cite{GelZakh93} it was shown that automorphisms of bihamiltonian structures
are enumerated by several functions of two variables.}, so it is
quite straightforward to construct an explicit diffeomorphism---the
moment one suspects that such a diffeomorphism exists.

\section{Webs for odd-dimensional bihamiltonian structures }\label{h2}\myLabel{h2}\relax 

In this section we suppose that $ \dim  M=2k-1 $.

\begin{definition} A pair of bilinear skewsymmetric pairings on a
finite-dimensional vector space $ V $ is {\em indecomposable\/} if the decomposition
of Theorem~\ref{th6.10} has only one component.

Call a pair of brackets $ \left\{,\right\}_{1} $ and $ \left\{,\right\}_{2} $ on $ M $ {\em micro-indecomposable\/} at
$ m\in M $ if the corresponding pair of bilinear pairings on $ {\mathcal T}_{m}^{*}M $ is
indecomposable. \end{definition}

\begin{definition} A pair of brackets $ \left\{,\right\}_{1} $ and $ \left\{,\right\}_{2} $ on $ M $ is {\em generic\/} at $ m\in M $, if
two corresponding bilinear pairings on $ {\mathcal T}_{m}^{*}M $ are in general position\footnote{For the purpose of this discussion, this means that $ \operatorname{GL}\left({\mathcal T}_{m}^{*}M\right) $-orbit of
the given pair of pairings is open.}. \end{definition}

Note that Theorem~\ref{th6.10} implies that an indecomposable pair of
parings on an odd-dimensional vector space $ W $ is isomorphic to $ {\mathcal K}_{2k-1} $, here
$ \dim  W=2k-1 $.

Now we can codify the program outlined in Section~\ref{h02}:

\begin{theorem} \label{th2.07}\myLabel{th2.07}\relax  (\cite{GelZakhWeb,GelZakh93}) Consider a pair of compatible
Poisson structures on an odd-dimensional manifold $ M $. This pair is generic
at $ m $ iff it is micro-indecomposable at $ m $. If it is micro-indecomposable
at $ m $, then it is micro-indecomposable at $ m' $ for any $ m' $ in a neighborhood
of $ m $.

If a pair is micro-indecomposable at $ m\in M $, then
\begin{enumerate}
\item
The web $ {\mathcal B}_{U} $ is a manifold for any small open neighborhood $ U $ of $ m $, in
other words, for any open $ U\ni m $ there is an open subset $ U' $, $ m\in U'\subset U $, such
that $ {\mathcal B}_{U'} $ is a manifold;
\item
The dimension of the manifold $ {\mathcal B}_{U} $ is $ \frac{\dim  M+1}{2} $;
\item
For any $ \lambda\in{\mathbb C}{\mathbb P}^{1} $ there is a foliation $ {\mathcal F}_{\lambda} $ on $ {\mathcal B}_{U} $ of codimension 1 such
that the subspace $ {\mathfrak C}_{\lambda} $ consists of smooth functions which are constant on
leaves of the foliation $ {\mathcal F}_{\lambda} $.
\item
Consider a micro-indecomposable pair of compatible Poisson
structures on a manifold $ M' $, and the corresponding manifold $ {\mathcal B}_{U'} $ with
foliations $ {\mathcal F}_{\lambda}' $. Suppose that both $ M $ and $ M' $ are analytic. If there is a
diffeomorphism $ \xi\colon {\mathcal B}_{U} \to {\mathcal B}_{U'} $ which sends the foliation $ {\mathcal F}_{\lambda} $ to the foliation
$ {\mathcal F}_{\lambda}' $ for any $ \lambda\in{\mathbb C}{\mathbb P}^{1} $, then the bihamiltonian structures on $ M $ and $ M' $ are
locally diffeomorphic. This local diffeomorphism is compatible with the
diffeomorphism $ \xi $.
\end{enumerate}
\end{theorem}

\begin{remark} Note that the conjecture of \cite{GelZakhWeb} implies that the last
statement of this theorem holds in the $ C^{\infty} $-case too. In \cite{Tur99Equi}
it is announced that this conjecture holds. \end{remark}

We are not going to repeat the proof of this theorem here, but we
sketch some arguments which should convince the reader that the first
several statements are true (it is the last one which is complicated).

Note that for the pair~\eqref{equ2.10} the pairing $ \lambda_{1}\left(,\right)_{1}+\lambda_{2}\left(,\right)_{2} $ is
degenerate (as any skewsymmetric pairing on an odd-dimensional vector
space) for any $ \left(\lambda_{1},\lambda_{2}\right)\in{\mathbb C}^{2} $, $ \left(\lambda_{1},\lambda_{2}\right)\not=\left(0,0\right) $, and has a $ 1 $-dimensional
null-space. (In other words, the dimension of the kernel does not jump up
for any $ \left(\lambda_{1},\lambda_{2}\right)\not=\left(0,0\right) $.) Moreover, Theorem~\ref{th6.10} momentarily implies
that a pair of pairings is indecomposable iff for any $ \left(\lambda_{1},\lambda_{2}\right)\not=\left(0,0\right) $ the
null-space of $ \lambda_{1}\left(,\right)_{1}+\lambda_{2}\left(,\right)_{2} $ is $ 1 $-dimensional.

This together with compactness of $ {\mathbb C}{\mathbb P}^{1} $ immediately implies that a
small deformation of $ {\mathcal K}_{2k-1} $ is indecomposable, thus isomorphic to $ {\mathcal K}_{2k-1} $.
In turn, this implies that a Zariski open (thus dense) subset of all
possible pairs consists of pairs isomorphic to $ {\mathcal K}_{2k-1} $. This shows that
the property of being generic coincides with indecomposability.

As a corollary, if a pair of brackets $ \left\{,\right\}_{1} $ and $ \left\{,\right\}_{2} $ is generic at
$ m\in M $, then in an appropriate neighborhood $ U $ of $ m $ the bracket $ \left\{,\right\}^{\left(\lambda_{1},\lambda_{2}\right)}
\buildrel{\text{def}}\over{=} \lambda_{1}\left\{,\right\}_{1}+\lambda_{2}\left\{,\right\}_{2} $ has a corank 1 in $ U $ for any $ \left(\lambda_{1},\lambda_{2}\right)\in{\mathbb C}^{2}\smallsetminus\left\{0\right\} $.

Now suppose that the bracket $ \left\{,\right\}^{\left(\lambda_{1},\lambda_{2}\right)} $ is Poisson. Since the rank
of the corresponding tensor field $ \eta $ is constant it is easy to see
(\cite{Kir76Loc,Wei83Loc}) that there is a (locally defined) Casimir
function $ F_{\lambda_{1},\lambda_{2}} $. Since the corank is 1, the level hypersurfaces of $ F_{\lambda_{1},\lambda_{2}} $
are canonically defined.

On the other hand, the normal direction $ {\mathbit n}_{\lambda_{1},\lambda_{2}} $ to the level
hypersurfaces of $ F_{\lambda_{1},\lambda_{2}} $ at $ m $ is the kernel of the corresponding
skewsymmetric pairing on $ {\mathcal T}_{m}^{*}M $. Let $ \lambda=\left(\lambda_{1}:\lambda_{2}\right)\in{\mathbb C}{\mathbb P}^{1} $, $ {\mathbit n}_{\lambda}\buildrel{\text{def}}\over{=}{\mathbit n}_{\lambda_{1},\lambda_{2}} $. Use the
isomorphism of $ {\mathcal T}_{m}^{*}M $ with the form~\eqref{equ2.10} to investigate how $ {\mathbit n}_{\lambda} $
depends on $ \lambda $. It is easy to see that the image of the vectors $ {\mathbit n}_{\lambda} $ in the
coordinate system of~\eqref{equ2.10} is proportional to
\begin{equation}
{\mathbit w}_{0}+\lambda{\mathbit w}_{2}+\dots +\lambda^{k}{\mathbit w}_{2k},
\label{equ2.20}\end{equation}\myLabel{equ2.20,}\relax 
thus taken for any $ k+1 $ distinct values $ \left\{\lambda_{i}\right\} $ of $ \lambda $ the vectors $ {\mathbit n}_{\lambda_{i}} $ span the
vector subspace $ \left< {\mathbit w}_{0},{\mathbit w}_{2},\dots ,{\mathbit w}_{2k} \right> $. Translating back to the language of
differential geometry, one obtains

\begin{corollary} Consider foliations $ {\mathcal F}_{\lambda} $ given by level sets of $ F_{\lambda_{1},\lambda_{2}} $, here
$ \lambda=\left(\lambda_{1}:\lambda_{2}\right) $. For any $ k+1 $ distinct values $ \lambda^{\left(0\right)},\lambda^{\left(1\right)},\dots ,\lambda^{\left(k\right)}\in{\mathbb C}{\mathbb P}^{1} $ the
foliations $ {\mathcal F}_{\lambda^{\left(0\right)}} $, $ {\mathcal F}_{\lambda^{\left(1\right)}} $, \dots , $ {\mathcal F}_{\lambda^{\left(k\right)}} $ intersect transversally, and the
intersection foliation $ {\mathcal F} $ does not depend on the choice of
$ \lambda^{\left(0\right)},\lambda^{\left(1\right)},\dots ,\lambda^{\left(k\right)} $. The foliation $ {\mathcal F} $ is a subfoliation of the foliation $ {\mathcal F}_{\lambda} $
for any $ \lambda\in{\mathbb C}{\mathbb P}^{1} $. \end{corollary}

Now one can momentarily see that the local base of the foliation $ {\mathcal F} $
coincides with the web $ {\mathcal B}_{U} $ of the bihamiltonian structure, and
the push-forward of $ {\mathcal F}_{\lambda} $ to the base of $ {\mathcal F} $ is a foliation on $ {\mathcal B}_{U} $ which
corresponds to the subspace $ {\mathfrak C}_{\lambda} $ of functions on $ {\mathcal B}_{U} $.

\begin{remark} \label{rem2.017}\myLabel{rem2.017}\relax  The proof of the last statement of Theorem~\ref{th2.07} might
be broken into two parts. The first one proves this statement under the
condition that both bihamiltonian structures admit an involution $ i\colon M \to
M $ such that $ \pi\circ i=i\circ\pi $, and such that $ i^{*}\left\{f,g\right\}_{a}=-\left\{i^{*}f,i^{*}g\right\}_{a} $ for any functions
$ f $ and $ g $ on $ M $ and any $ a=1,2 $ (here $ \pi $ is the projection from $ M $ to its web
$ {\mathcal B}_{M} $, and we suppose that $ M $ is small enough for the conditions of Theorem
~\ref{th2.07} to be applicable). Given such an involution, one can use the set
of fixed points of $ i $ as the common level set $ \left\{\varphi_{i}=0\right\} $ of would-be angle
variables $ \varphi_{i} $. After this choice it is possible to construct the angle
variables in purely geometric terms.

The second part of the proof consists of showing local existence of
such a section for any bihamiltonian structure of Theorem~\ref{th2.07}. This
part of the proof uses a hard cohomological statement related to
solvability of some overdetermined partial differential equations with
variable coefficients. The constant-coefficient variant of this
cohomological statement bears some similarity to the Dolbeault lemma.

In fact for the proof of Theorem~\ref{th1.10} only this constant
coefficients variant is needed, so it is possible that our proof of
Theorem~\ref{th1.10} may be significantly simplified. \end{remark}

\begin{remark} Note also that in many applications one may avoid using the
above cohomological statement, since one may be able to construct an
involution $ i $ explicitly. Say, for the structure~\eqref{equ45.20} the involution
is given by $ x_{j} \mapsto \left(-1\right)^{j}x_{j} $. \end{remark}

\section{Criterion of flatness }\label{h45}\myLabel{h45}\relax 

Here we prove Theorem~\ref{th1.10}. Suppose that conditions of Theorem
~\ref{th1.07} are satisfied. By Corollary~\ref{cor25.40} we know that the manifold
in question is odd-dimensional, and the pair of Poisson structures is
micro-indecomposable on an open subset $ U $. By Theorem~\ref{th2.07} the
corresponding web $ {\mathcal B}_{U} $ is a manifold with a family of foliations depending
on parameter $ \lambda\in{\mathbb C}{\mathbb P}^{1} $.

On the other hand, it is easy to describe this web explicitly. By
definition, for any given $ \lambda\in{\mathcal U} $ the level sets of the function $ F_{\lambda} $ are
unions of fibers of the foliation $ {\mathcal F}_{\lambda} $ on $ U $. Since the foliation $ {\mathcal F} $ is a
subfoliation of $ {\mathcal F}_{\lambda} $, the function $ F_{\lambda} $ is constant on leaves of $ {\mathcal F} $, thus
induces a function $ \widetilde{F}_{\lambda} $ on $ {\mathcal B}_{U} $. We obtain a mapping $ \widetilde{F}_{\bullet}\colon {\mathcal B}_{U}\times{\mathcal U} \to {\mathbb C}\colon \left(b,\lambda\right) \mapsto
\widetilde{F}_{\lambda}\left(b\right) $. Considered for variable $ \lambda\in{\mathcal U} $, the mapping $ \widetilde{F}_{\bullet} $ induces a mapping from
$ {\mathcal B}_{U} $ to the topological vector space $ C^{0}\left({\mathcal U}\right) $ of continuous functions on $ {\mathcal U} $.
This mapping sends a given point $ b\in{\mathcal B}_{U} $ to the function $ \widetilde{F}_{\lambda}\left(b\right) $ considered as
a function of $ \lambda $.

The topological space $ C^{0}\left({\mathcal U}\right) $ carries a canonical structure of a $ {\mathcal U} $-web
with $ \lambda\in{\mathcal U} $, with the subspace $ {\mathfrak C}_{\lambda} $ which consists of functions on $ C^{0}\left({\mathcal U}\right) $ of
the form $ f \mapsto \varphi\left(f|_{\lambda}\right) $ with an arbitrary smooth function $ \varphi $. The above
description of the mapping $ {\mathcal B}_{U} \mapsto C^{0}\left({\mathcal U}\right) $ shows that the $ {\mathcal U} $-web structure on
$ {\mathcal B}_{U} $ is induced from the $ {\mathcal U} $-web structure on $ C^{0}\left({\mathcal U}\right) $. One can also note that
by the condition of Theorem~\ref{th1.07} the mapping $ \widetilde{F}_{\bullet} $ is an immersion.
Indeed, the rank of $ d\widetilde{F}_{\bullet}|_{b} $ is exactly the dimension of the span of $ dF_{\lambda}|_{m} $
for $ \lambda\in{\mathcal U} $ (here $ b $ is the projection of $ m $ to $ {\mathcal B}_{U} $), which coincides with $ \dim 
{\mathcal B}_{U} $.

Now suppose that the conditions of Theorem~\ref{th1.10} are satisfied, so
$ F_{\lambda} $ depends polynomially on $ \lambda $. Denote the degree of $ F_{\lambda} $ in $ \lambda $ by $ d $. By
conditions of Theorem~\ref{th1.10} one has $ d<\frac{\dim  M}{2} $. On the other hand, by
Lemma~\ref{lm25.20} the degree in $ \lambda $ of $ dF_{\lambda}|_{m} $ cannot be less than $ \frac{\dim  M-1}{2} $.
Thus the degree of $ F_{\lambda}\left(m\right) $ is exactly $ \frac{\dim  M-1}{2} $ near $ m_{0} $.

In particular, for any $ b\in{\mathcal B}_{U} $ the image $ \widetilde{F}_{\bullet}\left(b\right) $ of $ b $ is in fact inside
the $ \left(d+1\right) $-dimensional vector space of polynomials of degree $ d $ in $ \lambda $.
Denote this vector space by $ {\mathcal P}_{d} $. Similarly to $ C^{0}\left({\mathcal U}\right) $, it carries a natural
structure of $ {\mathbb C} $-web, moreover, this structure may be extended to become
$ {\mathbb C}{\mathbb P}^{1} $-web by noting that $ {\mathcal P}_{d}=\Gamma\left({\mathbb C}{\mathbb P}^{1},{\mathcal O}\left(d\right)\right) $. Since $ \dim {\mathcal P}_{d}=d+1=\dim  {\mathcal B}_{U} $, and $ \widetilde{F} $ is
an immersion, one can see that

\begin{proposition} In conditions of Theorem~\ref{th1.10}
\begin{enumerate}
\item
The mapping $ \widetilde{F} $ is a local diffeomorphism compatible with structures
of webs on $ {\mathcal B}_{U} $ and $ {\mathcal P}_{d} $;
\item
The structure of the web on $ {\mathcal P}_{d} $ is invariant w.r.t.~parallel
translations on $ {\mathcal P}_{d} $;
\item
For any polynomial $ p\left(\lambda\right) $ of degree $ d $ the
family of functions $ G_{\lambda}\left(m\right)\buildrel{\text{def}}\over{=}F_{\lambda}\left(m\right)+p\left(\lambda\right) $ on $ M $ satisfies the conditions of
Theorem~\ref{th1.10};
\item
The mapping $ \widetilde{G}_{\bullet}\colon {\mathcal B}_{U} \to {\mathcal P}_{d} $ associated to $ G_{\lambda}\left(m\right) $ is a parallel translation
of the mapping $ \widetilde{F}_{\bullet} $.
\end{enumerate}
\end{proposition}

In other words, for any point $ p\in{\mathcal P}_{d} $ and any point $ b\in{\mathcal B}_{U} $ one can find a
local diffeomorphism of the webs $ {\mathcal B}_{U} $ and $ {\mathcal P}_{d} $ which sends $ b $ to $ p $. This shows
that if two bihamiltonian structures satisfy the conditions of Theorem
~\ref{th1.10}, then two corresponding webs are isomorphic.

\begin{lemma} There is a pair of compatible translation-invariant Poisson
brackets on $ {\mathbb C}^{2k-1} $ which may be equipped with a family of functions
satisfying Theorem~\ref{th1.10}. \end{lemma}

\begin{proof} For translation-invariant brackets the tensor fields $ \eta_{1} $ and $ \eta_{2} $
are constant, thus to describe the bracket we need to describe the
pairing on {\em one\/} cotangent space. On the other hand, we know that to
satisfy Theorem~\ref{th1.07} this pairing should be isomorphic to $ {\mathcal K}_{2k-1} $, thus
any pair which satisfies the lemma is isomorphic to one given by the
brackets~\eqref{equ45.20}.

Obviously,
\begin{equation}
F_{\lambda}=x_{0}+\lambda x_{2}+\lambda^{2}x_{4}+\dots +\lambda^{k-1}x_{2k-2}
\label{equ45.25}\end{equation}\myLabel{equ45.25,}\relax 
satisfies the conditions of Theorem~\ref{th1.10}, and the bracket $ \left\{,\right\}_{1} $ has
exactly one independent Casimir function $ x_{2k-2} $. \end{proof}

\begin{corollary} \label{cor45.40}\myLabel{cor45.40}\relax  In conditions of Theorem~\ref{th1.10} the web $ {\mathcal B}_{U} $ is
locally isomorphic to the web corresponding to the bihamiltonian structure
given by Equation~\eqref{equ45.20}, here $ \dim  M=2k-1 $. \end{corollary}

Indeed, both these webs are locally isomorphic to the web on $ {\mathcal P}_{k-1} $.

Now the last part of Theorem~\ref{th2.07} implies that the bihamiltonian
structure on $ M $ is isomorphic to the structure given by~\eqref{equ45.20}, which
finishes the proof of Theorem~\ref{th1.10}.

\section{Examples of non-flat structures }\label{h47}\myLabel{h47}\relax 

Here we show that Theorem~\ref{th1.10} is not a tautology. To do this, we
construct a huge pool of bihamiltonian structures which satisfy the
conditions of Theorem~\ref{th1.07}, but are not isomorphic to each other (in
particular, only one of them is flat). All these structures are
integrable by the anchored Lenard scheme (see Section~\ref{h48}, compare with
descriptions \cite{Mag78Sim,GelDor79Ham} in symplectic settings).

All the constructions below can be performed in $ C^{\infty} $-geometry and in
analytic geometry (unless explicitly specified). We state the $ C^{\infty} $-case
only.

Fix an open subset $ {\mathcal B}\subset{\mathbb R}^{2} $ and a smooth function $ f\left(x,y\right) $ of two
variables $ \left(x,y\right)\in{\mathcal B} $. Consider two brackets on $ {\mathcal B}\times{\mathbb R} $ defined by
\begin{equation}
\begin{gathered}
\left\{x,y\right\}_{1}=\left\{y,z\right\}_{1}=0,\qquad \left\{x,y\right\}_{2}=\left\{x,z\right\}_{2}=0,
\\
\left\{x,z\right\}_{1}=\frac{\partial f\left(x,y\right)}{\partial y},\qquad \left\{y,z\right\}_{2}=-\frac{\partial f\left(x,y\right)}{\partial x}.
\end{gathered}
\label{equ47.10}\end{equation}\myLabel{equ47.10,}\relax 
Obviously, these two brackets form a bihamiltonian structure on $ {\mathcal B}\times{\mathbb R} $.

\begin{definition} \label{def47.11}\myLabel{def47.11}\relax  Denote this bihamiltonian structure on $ {\mathcal B}\times{\mathbb R} $ by $ M_{f} $. \end{definition}

Assume that both $ \frac{\partial f}{\partial x} $ and $ \frac{\partial f}{\partial y} $ do not vanish in $ {\mathcal B} $. One can see
that
any non-zero linear combination $ \lambda_{1}\left\{,\right\}_{1}+\lambda_{2}\left\{,\right\}_{2} $ has rank 2, thus has
exactly one independent Casimir function near $ b\times z\in{\mathcal B}\times{\mathbb R} $. Moreover, it is
easy to construct a family of Casimir functions for different $ \lambda\buildrel{\text{def}}\over{=}\lambda_{1}:\lambda_{2} $
which depend smoothly on $ \lambda\in{\mathcal U}\subset{\mathbb R} $ (compare with Section~\ref{h48}). Thus the
structure~\eqref{equ47.10} satisfies conditions of Theorem~\ref{th1.07}, and is
homogeneous of type\footnote{Using Theorem~\ref{th2.07}, it is easy to show that any {\em analytic\/} homogeneous
bihamiltonian structure of type (3) is locally isomorphic to this
structure for an appropriate $ f $. For the discussion of global geometry for
such structures, see \cite{Rig95Tis}.} (3).

One can write explicitly Casimir functions for the values $ \lambda_{1}:\lambda_{2} $
being $ \infty $, 0, and 1, i.e., for $ \left\{,\right\}_{1} $, for $ \left\{,\right\}_{2} $, and for $ \left\{,\right\}_{1}+\left\{,\right\}_{2} $. They have
the form $ F_{\infty}\left(y\right) $, $ F_{0}\left(x\right) $, and $ F_{1}\left(f\left(x,y\right)\right) $ for arbitrary functions $ F_{0} $, $ F_{1} $,
$ F_{\infty} $. This implies

\begin{lemma} Consider a bihamiltonian structure $ \left(M,\left\{,\right\}_{a},\left\{,\right\}_{b}\right) $ and a mapping
$ p\colon M \to M_{f} $ which is a local isomorphism of bihamiltonian structures.
Consider $ x,y,f $ as functions on $ M_{f} $. Then $ x\circ p $ is a Casimir function for
$ \left\{,\right\}_{b} $, $ y\circ p $ is a Casimir function for $ \left\{,\right\}_{a} $, and $ f\circ p $ is a Casimir function
for $ \left\{,\right\}_{a}+\left\{,\right\}_{b} $. \end{lemma}

Since $ \left\{,\right\}_{1} $, $ \left\{,\right\}_{2} $, and $ \left\{,\right\}_{1}+\left\{,\right\}_{2} $ are of corank 1, the structure
$ \left(M,\left\{,\right\}_{a},\left\{,\right\}_{b}\right) $ of the lemma determines the functions $ x'=x\circ p $, $ y'=y\circ p $, and
$ f'=f\circ p $ uniquely up to transformations of the form $ \widetilde{A}=\alpha_{A}\left(A\right) $, $ A\in\left\{x',y',p'\right\} $.
Thus $ \left(M,\left\{,\right\}_{a},\left\{,\right\}_{b}\right) $ determines the function $ f\left(x,y\right) $ up to transformations
of the form $ \widetilde{f}=\gamma\left(f\left(\alpha\left(x\right),\beta\left(y\right)\right)\right) $. Indeed, the graph of $ f\left(x,y\right) $ coincides with
the image of $ M $ w.r.t.~the mapping $ x'\times y'\times f' $.

\begin{lemma} Consider a smooth function $ f\colon U_{1}\times U_{2} \to U_{3} $, $ U_{1,2,3}\subset{\mathbb R} $, and
diffeomorphisms $ \alpha\colon U_{1}' \to U_{1} $, $ \beta\colon U_{2}' \to U_{2} $, $ \gamma\colon U_{3} \to U_{3}' $. Let
\begin{equation}
g\left(x,y\right) \buildrel{\text{def}}\over{=} \gamma\left(\alpha\left(x\right),\beta\left(y\right)\right).
\notag\end{equation}
Then the bihamiltonian structures $ M_{f} $ and $ M_{g} $ are isomorphic. \end{lemma}

\begin{proof} It is easy to write the diffeomorphism explicitly as $ x'=\alpha\left(x\right) $,
$ y'=\beta\left(y\right) $, $ z'=\xi\left(x,y\right)z $ with an appropriate function $ \xi\left(x,y\right) $. \end{proof}

The next step is to introduce additional conditions on a function
$ g\left(x,y\right) $ which would define the diffeomorphisms $ \alpha $, $ \beta $, $ \gamma $ almost uniquely.
First, assume $ \left(0,0\right)\in{\mathcal B} $, $ f\left(0,0\right)=0 $. Consider $ O=\left(0,0,0\right) $ as a marked point on
$ {\mathcal B}\times{\mathbb R} $. Then the condition that $ x'- $ and $ y' $-coordinates of $ O $ remain 0 leads
to $ \alpha\left(0\right)=0 $ and $ \beta\left(0\right)=0 $.

Given a function $ f\left(x,y\right) $ such that $ f\left(0,0\right)=0 $, $ \partial f/\partial x\not=0 $ and $ \partial f/\partial y\not=0 $ near
$ x=y=0 $, one can find local coordinate changes $ x'=\alpha\left(x\right) $, $ y'=\beta\left(y\right) $, and
$ f'=\gamma\left(f\right) $, $ \alpha\left(0\right)=\beta\left(0\right)=\gamma\left(0\right)=0 $, such that $ \varphi\left(x',y'\right)\buildrel{\text{def}}\over{=}\gamma\left(f\left(\alpha^{-1}\left(x'\right),\beta^{-1}\left(y'\right)\right)\right) $
satisfies
\begin{equation}
\begin{gathered}
\frac{\partial\varphi}{\partial x'}|_{x'=0}=1\text{, }\frac{\partial\varphi}{\partial y'}|_{x'=0}=1,
\\
\frac{\partial\varphi}{\partial x'}|_{y'=0}=\frac{\partial\varphi}{\partial y'}|_{y'=0},
\end{gathered}
\qquad \varphi\left(0,0\right)=0.
\label{equ47.20}\end{equation}\myLabel{equ47.20,}\relax 
Moreover, such a coordinate change and the function $ \varphi $ are defined
uniquely up to simultaneous multiplication of $ x' $, $ y' $ and $ \varphi $ by the same
constant. If $ \frac{\partial^{2}\varphi}{\partial x'\partial y'}|_{\left(0,0\right)}\not=0 $, then this last degree of freedom may
be eliminated by a requirement that $ \frac{\partial^{2}\varphi}{\partial x'\partial y'}|_{\left(0,0\right)}=1 $. (In fact if
$ \varphi\left(x',y'\right)\not\equiv x'+y' $, then one can fix coordinates $ x' $ and $ y' $ by normalizing an
appropriate derivative of $ \varphi $ of higher order.)

These arguments lead to the following statement from (\cite{GelZakhWeb,%
GelZakh93}):

\begin{theorem} \label{th47.13}\myLabel{th47.13}\relax  Consider two functions $ \varphi\left(x,y\right) $ and $ \varphi'\left(x,y\right) $ defined in a
neighborhoods $ {\mathcal B} $, $ {\mathcal B}' $ of (0,0). Suppose that both $ \varphi $ and $ \varphi' $ satisfy
~\eqref{equ47.20}. There exists a local diffeomorphism between $ M_{\varphi} $ and $ M_{\varphi'} $ which
preserves the point (0,0,0) iff there exists $ C\in{\mathbb R} $, $ C\not=0 $, such that
$ \varphi\left(Cx,Cy\right)\equiv C\varphi'\left(x,y\right) $ near $ x=y=0 $. \end{theorem}

\begin{corollary} \label{cor47.30}\myLabel{cor47.30}\relax  If $ \varphi\left(x,y\right) $ satisfies~\eqref{equ47.20}, then $ M_{\varphi} $ is flat iff
$ \varphi\left(x,y\right)=x+y $. \end{corollary}

\begin{proof} Indeed, $ \varphi\left(x,y\right)=f\left(x,y\right)=x+y $ defines a structure with constant
coefficients (compare with~\eqref{equ45.20}), thus a flat one. Any other
function $ \varphi\left(x,y\right) $ which satisfies~\eqref{equ47.20} will define a non-isomorphic
bihamiltonian structure, thus a non-flat one. \end{proof}

As a corollary, one obtains a lot of structures which are not flat,
thus cannot satisfy the conditions of Theorem~\ref{th1.10}. The next logical
step is to check whether these bihamiltonian structures are ``integrable''.
To do so, we need a formalization of the notion of integrability. One of
the simplest such notions is integrability by the anchored Lenard scheme,
which is introduced in Section~\ref{h48}. Example~\ref{ex48.80} will demonstrate
that any homogeneous bihamiltonian structure $ M_{f} $ is integrable\footnote{In fact Proposition~\ref{prop55.60} will show that {\em any\/} homogeneous structure is
Lenard-integrable.} in this
sense.

\section{One counterexample }\label{h62}\myLabel{h62}\relax 

The examples of Section~\ref{h47} show that one cannot expect to prove
the conclusion of Theorem~\ref{th1.10} in conditions of Theorem~\ref{th1.07}, even
if one requires $ {\mathcal U} $ to become the whole complex plane $ {\mathbb C} $, and requires $ F_{\lambda} $ to
depend analytically on $ \lambda $. Recall that the notation $ M_{f} $ was introduced in
Definition~\ref{def47.11}.

To show that in fact even the restriction on the degree of the
polynomial in Theorem~\ref{th1.10} cannot be improved if $ k=2 $, consider

\begin{definition} Consider a bihamiltonian structure on $ M $ and a smooth
function on $ {\mathbb R}\times M $, $ \left(\lambda,m\right) \mapsto C_{\lambda}\left(m\right) $. Call this function a $ \left[d\right] $-{\em family\/} if
for any fixed $ m\in M $ the function $ C_{\lambda}\left(m\right) $ of $ \lambda $ depends on $ \lambda $ as a polynomial of
degree $ d $ or less, and for any fixed $ \lambda\in{\mathbb R} $ the function $ C_{\lambda} $ of $ m $ is a Casimir
function for $ \lambda\left\{,\right\}_{1}+\left\{,\right\}_{2} $. \end{definition}

\begin{proposition} \label{prop47.40}\myLabel{prop47.40}\relax  There exists a function $ \varphi\left(x,y\right)\not\equiv x+y $ which satisfies
~\eqref{equ47.20}, and such that the bihamiltonian structure $ M_{\varphi} $ admits a
$ \left[2\right] $-family $ F_{\lambda} $. \end{proposition}

\begin{proof} Actually it is possible to describe {\em all\/} analytic functions
$ f\left(x,y\right) $ such that $ M_{f} $ admits a $ \left[2\right] $-family $ F_{\lambda} $, at least if we are allowed to
restrict our attention to smaller open subsets.

If $ G_{\lambda} $ is a $ \left[1\right] $-family, then $ \left(d\lambda+e\right)G_{\lambda}+a\lambda^{2}+b\lambda+c $, $ a,b,c,d,e\in{\mathbb C} $, gives a
$ \left[2\right] $-family. Call such families {\em simple\/} families. By Theorem~\ref{th1.10} a
simple family may exist on a flat bihamiltonian structure only.

Given an open subset $ {\mathcal U}\subset{\mathbb C} $ with two analytic functions $ \eta,\zeta\colon {\mathcal U} \to {\mathbb C} $,
and an open subset $ {\mathfrak B}\subset{\mathbb C}\times{\mathbb C} $ with an analytic function $ \Lambda\colon {\mathfrak B} \to {\mathcal U} $, let
\begin{equation}
F_{\lambda}^{\left(\eta\zeta\Lambda\right)}\left(x,y\right)\buildrel{\text{def}}\over{=}\left(\Lambda\left(x,y\right)-\lambda\right)^{2}y+\zeta\left(\Lambda\left(x,y\right)\right)+\lambda\eta\left(\Lambda\left(x,y\right)\right).
\notag\end{equation}
\begin{lemma} \label{lm62.20}\myLabel{lm62.20}\relax  Consider an analytic bihamiltonian structure $ M_{f} $ with a
$ \left[2\right] $-family $ F_{\lambda} $ which is not simple. Then there exists an open subset $ {\mathcal U}\subset{\mathbb C} $
with two analytic functions $ \eta,\zeta\colon {\mathcal U} \to {\mathbb C} $, and an open subset $ {\mathfrak B}\subset{\mathbb C}\times{\mathbb C} $ with an
analytic function $ \Lambda\colon {\mathfrak B} \to {\mathcal U} $ such that $ F_{\lambda}\left(x,y,z\right)=F_{\lambda}^{\left(\eta\zeta\Lambda\right)}\left(x,y\right) $ and
\begin{equation}
\frac{d\zeta}{dt}+t\frac{d\eta}{dt}=0\quad \text{if }t\in{\mathcal U},\qquad x=\Lambda^{2}\left(x,y\right)y+\zeta\left(\Lambda\left(x,y\right)\right)\quad \text{if }\left(x,y\right)\in{\mathfrak B},
\notag\end{equation}
In particular, $ M_{f} $ is locally isomorphic to an open subset of
$ M_{F_{1}^{\left(\eta\zeta\Lambda\right)}} $. \end{lemma}

\begin{proof} Write $ F_{\lambda}=\lambda^{2}H_{0}+\lambda H_{1}+H_{2} $. Being a Casimir function for $ \left\{,\right\}_{1} $, $ H_{2} $
depends on $ x $ only, similarly $ H_{0} $ depends on $ y $ only. If $ H_{0} $ does not depend
on $ y $, then $ H_{0}=\operatorname{const} $, thus $ F_{\lambda} $ is simple. Similarly one can exclude the
case when $ H_{2} $ does not depend on $ x $. Restricting to a smaller open subset,
one can assume that $ x=\alpha\left(H_{2}\right) $, $ y=\beta\left(H_{0}\right) $, thus one can consider $ H_{0} $ and $ H_{2} $
instead of coordinates $ x $ and $ y $ on $ {\mathcal B} $. In particular, we may assume that
$ H_{2}=x $, $ H_{0}=y $.

By Lemma~\ref{lm25.20}, given a point $ \left(x,y\right)\in{\mathcal B} $, $ dF_{\lambda}|_{\left(x,y\right)} $ may be written
as $ p\left(\lambda\right)\widetilde{w}_{\lambda} $, here $ p\left(\lambda\right) $ is a scalar polynomial in $ \lambda $, and $ \widetilde{w}_{\lambda} $ is a
vector-valued polynomial of degree 1 in $ \lambda $. Thus $ \deg  p=1 $, denote the zero
of $ p $ by $ \Lambda $. We conclude that for each $ \left(x,y\right) $ there is $ \Lambda\left(x,y\right) $ such that
$ dF_{\lambda}|_{\left(x,y\right)}=0 $ if $ \lambda=\Lambda\left(x,y\right) $.

Restricting to an appropriate open subset of $ {\mathcal B} $, we may assume that $ \Lambda $
depends analytically on $ x $ and $ y $. If $ \Lambda $ is constant, then $ dF_{\Lambda}=0 $ implies
that $ \frac{F_{\lambda}\left(x,y\right)-F_{\lambda}\left(x_{0},y_{0}\right)}{\lambda-\Lambda} $ is linear, thus $ F_{\lambda} $ is simple. Thus,
decreasing $ {\mathcal B} $ again, we may assume that either $ \left(\Lambda,y\right) $ or $ \left(\Lambda,x\right) $ give a local
coordinate system on $ {\mathcal B} $. Assume that $ \left(\Lambda,y\right) $ is a local coordinate system.

The condition $ dF_{\Lambda}|_{\left(x,y\right)}=0 $ implies
\begin{equation}
dx=-\Lambda^{2}dy-\Lambda dH_{1},
\label{equ47.30}\end{equation}\myLabel{equ47.30,}\relax 
Thus $ 2 $-form $ d\left(-\Lambda^{2}dy-\Lambda dH_{1}\right) $ vanishes, in other words, $ 2\Lambda d\Lambda dy+d\Lambda dH_{1}=0 $. One
can conclude that in the coordinate system $ \left(\Lambda,y\right) $ one has
$ \frac{\partial H_{1}}{\partial y}|_{\Lambda=\operatorname{const}}=-2\Lambda $, or $ H_{1}=-2\Lambda y+\eta\left(\Lambda\right) $ with an unknown function
$ \eta\left(t\right) $. Equation~\eqref{equ47.30} leads to
\begin{equation}
x=\Lambda^{2}y+\zeta\left(\Lambda\right),\qquad \frac{d\zeta}{dt}=-t\frac{d\eta}{dt}.
\label{equ47.40}\end{equation}\myLabel{equ47.40,}\relax 
This leads to a formula for $ F_{\lambda} $ in coordinates $ y $ and $ \Lambda $:
\begin{equation}
F_{\lambda}=\left(\lambda-\Lambda\right)^{2}y+\zeta\left(\Lambda\right)+\lambda\eta\left(\Lambda\right),\qquad \frac{d\zeta}{dt}=-t\frac{d\eta}{dt}.
\label{equ47.45}\end{equation}\myLabel{equ47.45,}\relax 
By~\eqref{equ47.10}, $ F_{1}=\gamma\left(f\right) $ for an appropriate function $ \gamma $. If
$ F_{1}\left(x,y\right)\equiv \operatorname{const} $, then $ \frac{F_{\lambda}-F_{1}}{\lambda-1} $ is linear, thus $ F_{\lambda} $ is simple. Hence
decreasing $ {\mathcal B} $ we may assume that one can write $ f=\varepsilon\left(F_{1}\right) $ for an appropriate
function $ \varepsilon $. Thus $ M_{F_{1}} $ locally isomorphic to $ M_{f} $.

Moreover,~\eqref{equ47.40} implies that $ \left(\Lambda,x\right) $ is also a coordinate system
on an open subset of $ {\mathcal B} $. Exchanging $ x $ and $ y $, we see that our assumption
that $ \left(\Lambda,y\right) $ is a local coordinate system is {\em always\/} satisfied. \end{proof}

\begin{lemma} \label{lm62.30}\myLabel{lm62.30}\relax  There is a way to associate to an open subset $ {\mathcal U}\subset{\mathbb C} $ and a
function $ \eta\colon {\mathcal U} \to {\mathbb C} $ a homogeneous bihamiltonian structure $ M^{\left(\eta\right)} $ of type (3)
with a family of function $ F_{\lambda}^{\left(\eta\right)} $ which is quadratic in $ \lambda $. In conditions of
Lemma~\ref{lm62.20} there is a diffeomorphism of an open subset of $ M_{f} $ with an
open subset of $ M^{\left(\eta\right)} $ and $ C\in{\mathbb C} $ such that the diffeomorphism sends one
bihamiltonian structure to another, and the family $ F_{\lambda} $ to the family
$ F_{\lambda}^{\left(\eta\right)}+C $. A change of $ \eta $ of the form $ \eta'\left(t\right)=\eta\left(t\right)+at+b $ leads to an isomorphic
bihamiltonian structure with the isomorphism sending $ F_{\lambda}^{\left(\eta\right)} $ to
$ F'{}_{\lambda}^{\left(\eta\right)}+A\lambda^{2}+B\lambda+C $, $ A,B,C\in{\mathbb C} $. \end{lemma}

\begin{proof} Indeed, given the functions $ \zeta $ and $ \eta $, let $ \bar{\Sigma} = \left\{ \left(x,y,\Lambda\right)\in{\mathbb C}^{3} \mid x
=\Lambda^{2}y+\zeta\left(\Lambda\right)\right\} $. Let $ \Sigma=\left\{\left(x,y,\Lambda\right)\in\bar{\Sigma} \mid y\not=\frac{1}{2}\frac{d\eta}{d\Lambda},\Lambda\not=0\right\} $ be the subset of $ \bar{\Sigma} $
where $ x $ and $ y $ provide a local coordinate system. Plugging into~\eqref{equ47.45},
one obtains functions $ \Lambda $, $ x $, $ y $, $ F_{\lambda} $ defined on $ \Sigma $.

Functions $ x $ and $ y $ provide a local coordinate system near any point
of $ \Sigma $. On an open subset $ \Sigma_{1}\subset\Sigma $ one has $ \frac{\partial F_{1}}{\partial x}\not=0 $, $ \frac{\partial F_{1}}{\partial y}\not=0 $, thus putting
$ F_{1} $ into~\eqref{equ47.10} instead of $ f $ defines a homogeneous bihamiltonian
structure on $ M^{\left(\eta,\zeta\right)}\buildrel{\text{def}}\over{=}\Sigma_{1}\times{\mathbb C} $. As we have seen in the proof of Lemma
~\ref{lm62.20}, an open piece of $ M_{f} $ is isomorphic to an open piece of $ M^{\left(\eta,\zeta\right)} $,
moreover, the families $ F_{\lambda} $ are preserved by this isomorphism.

By~\eqref{equ47.40},~\eqref{equ47.45}, a change of the form $ \eta'\left(t\right)=\eta\left(t\right)+2a t+b $
together with the change $ \zeta'\left(t\right)= \zeta\left(t\right) $-at$ ^{2}+d $ would lead to a parallel
translation of the surface $ \bar{\Sigma} $, and to the required change of functions $ F_{\lambda} $.
In particular, a change in $ \zeta $ only will not change $ M^{\left(\eta,\zeta\right)} $, and will change
the family $ F_{\lambda} $ by an additive constant only. \end{proof}

\begin{lemma} In the conditions of Lemma~\ref{lm62.30} the family $ F_{\lambda}^{\left(\eta\right)} $ on $ M^{\left(\eta\right)} $
is a $ \left[2\right] $-family. \end{lemma}

\begin{proof} By construction of bihamiltonian structure on $ M^{\left(\eta\right)} $, the
functions $ F_{0}^{\left(\eta\right)}\equiv x $, $ F_{1}^{\left(\eta\right)} $, and the leading coefficient $ H_{0} $ of the quadratic
family $ F_{\lambda}^{\left(\eta\right)} $ are Casimir functions for $ \left\{,\right\}_{1} $, $ \left\{,\right\}_{1}+\left\{,\right\}_{2} $, and $ \left\{,\right\}_{2} $
correspondingly. Fix $ m_{0}\in M $. Then $ dF_{\lambda}^{\left(\eta\right)}|_{m_{0}} $ is a vector-function which is
quadratic in $ \lambda $. Moreover, at $ \lambda=\Lambda\left(m_{0}\right) $ this vector-function vanishes, thus
$ dF_{\lambda}^{\left(\eta\right)}|_{m_{0}}=\left(\lambda-\Lambda\left(m_{0}\right)\right)w\left(\lambda\right) $, here $ w\left(\lambda\right) $ is a vector-function of degree 1 in $ \lambda $.
Extend $ w $ to become a homogeneous function $ w\left(\lambda_{1},\lambda_{2}\right) $ of homogeneity degree
1, here $ \lambda=\lambda_{1}:\lambda_{2} $.

This function $ w\left(\lambda_{1},\lambda_{2}\right) $ is in the null-space of $ \lambda_{1}\left(,\right)_{1}+\lambda_{2}\left(,\right)_{2} $ for
three values 0, 1 and $ \infty $ of $ \lambda_{1}:\lambda_{2} $. However, the null-space for a pair of
pairing which is isomorphic to $ {\mathcal K}_{3} $ depends linearly on $ \lambda_{1}:\lambda_{2} $ (compare with
~\eqref{equ2.10}). We conclude that $ w\left(\lambda_{1},\lambda_{2}\right) $ is in the null-space for
$ \lambda_{1}\left(,\right)_{1}+\lambda_{2}\left(,\right)_{2} $ for any $ \lambda_{1},\lambda_{2} $, thus $ F_{\lambda}^{\left(\eta\right)} $ is a $ \left[2\right] $-family indeed.\end{proof}

\begin{lemma} The bihamiltonian structure $ M^{\left(\eta\right)} $ of Lemma~\ref{lm62.30} is flat on
any open subset iff $ \frac{d^{2}\eta}{dt^{2}}\equiv 0 $. \end{lemma}

\begin{proof} By arguments of Section~\ref{h47}, $ M_{F_{1}^{\left(\eta\right)}} $ is flat on an open subset
iff there is a local dependence between $ F_{1}^{\left(\eta\right)} $, $ x $ and $ y $ of the form
$ a\left(F_{1}^{\left(\eta\right)}\right)+b\left(x\right)+c\left(y\right)=0 $. If $ \eta=\zeta\equiv 0 $, then $ F_{\lambda}^{\left(\eta\right)}=\lambda^{2}y\pm2\lambda\sqrt{xy}+x $, thus
$ F_{1}^{\left(\eta\right)}=\left(\sqrt{x}\pm\sqrt{y}\right)^{2} $, thus $ M_{F_{1}^{\left(\eta\right)}} $ is flat. By Lemma~\ref{lm62.30} this proves the
``if'' part.

If a dependence $ a\left(F_{1}^{\left(\eta\right)}\right)+b\left(x\right)+c\left(y\right)=0 $ exists, then
$ \frac{\partial}{\partial\Lambda}|_{y=\operatorname{const}}\left(a\left(F_{1}^{\left(\eta\right)}\right)+b\left(x\right)\right)=0 $, or
\begin{equation}
\left(\Lambda-1\right)\alpha\left(F_{1}^{\left(\eta\right)}\right)=-\Lambda\beta\left(x\right),
\label{equ47.47}\end{equation}\myLabel{equ47.47,}\relax 
here $ \alpha\left(F_{1}^{\left(\eta\right)}\right)=da/dF_{1}^{\left(\eta\right)} $, $ \beta\left(x\right)=db/dx. $

Taking derivative $ \frac{\partial}{\partial y}|_{\Lambda=\operatorname{const}} $ of~\eqref{equ47.47}, and dividing by the
cube of~\eqref{equ47.47}, one obtains $ \left(\alpha^{-3}\frac{d\alpha}{dF_{1}^{\left(\eta\right)}}\right)\left(F_{1}^{\left(\eta\right)}\right)=\beta^{-3}\frac{d\beta}{dx}\left(x\right) $.
Since $ F_{1}^{\left(\eta\right)} $ and $ x $ are independent, and $ \alpha\not\equiv 0 $, $ \beta\not\equiv 0 $, we conclude that
\begin{equation}
\frac{d\alpha}{dF_{1}^{\left(\eta\right)}}=C\alpha^{3},\qquad \frac{d\beta}{dx}=C\beta^{3}.
\notag\end{equation}
Thus $ \alpha\left(F_{1}^{\left(\eta\right)}\right)=D/\sqrt{F_{1}^{\left(\eta\right)}-\varphi_{0}} $, $ \beta\left(x\right)=D/\sqrt{x-x_{0}} $, $ x_{0},\varphi_{0},D\in{\mathbb C} $, $ D\not=0 $. Hence
$ \left(\Lambda-1\right)^{2}\left(x-x_{0}\right)=\Lambda^{2}\left(F_{1}-\varphi_{0}\right) $ by~\eqref{equ47.47}, and $ \left(-2\Lambda+1\right)\left(\zeta\left(\Lambda\right)-\zeta_{0}\right)=\Lambda^{2}\left(\eta\left(\Lambda\right)-\eta_{0}\right) $ for
appropriate $ \eta_{0},\zeta_{0}\in{\mathbb C} $. Together with $ \zeta'\left(t\right)=-t\eta'\left(t\right) $ this shows that
$ \eta\left(\Lambda\right)=C\Lambda+E $. \end{proof}

This finishes the proof of Proposition~\ref{prop47.40}. \end{proof}

\section{Anchored Lenard scheme }\label{h48}\myLabel{h48}\relax 

Recall how Lenard scheme works\footnote{For history of Lenard scheme and of the term ``Lenard scheme'' see
\cite{KosMag96Lax}.}. Since descriptions of Lenard scheme
are in many cases based on an assumption that the Poisson bracket is
symplectic, here we supply as many details as possible to unravel the
relation of the anchored Lenard scheme with Casimir functions which
depend on parameters (such functions do not exists in symplectic
situation). In turn, such functions are directly related to webs.

\begin{remark} \label{rem48.02}\myLabel{rem48.02}\relax  Before we proceed with description of the problem which
the Lenard scheme solves, we need to resolve a possible ambiguity. The
Lenard scheme ``as a method of integration'' consists of recurrence
relations and initial data for these relations. However, the existing
{\em formalizations\/} of Lenard scheme (e.g., \cite{Mag78Sim,GelDor79Ham,%
KosMag96Lax}) consider the recurrence relations only, omitting the
initial data. The latter approach has an advantage of being more general,
in particular, it works in symplectic case too. However, this approach
does not address the question when recurrence relations have solutions
(these relations are overdetermined), in particular they do not specify
how to find the initial data which would make the Lenard scheme succeed.

Since in our settings symplectic Poisson structures have only a
tangential r\^ole, here we consider only the variant of Lenard scheme which
is important in applications, when both the initial data and the
recurrence relation are specified. To avoid any confusion, we call this
variant the {\em anchored Lenard scheme}. For this scheme we will not only be
able to describe when recurrence relations are solvable, we will also
describe which bihamiltonian systems are completely integrable by Lenard
scheme. As Theorem~\ref{th55.50} will show, such systems are never symplectic,
which may explain why the anchored Lenard scheme was not formalized
before. \end{remark}

The method of {\em first integrals\/} to ``integrate'' a system of ordinary
differential equations $ \frac{d}{dt}m\left(t\right)=v\left(m\left(t\right)\right) $, $ m\in M $, starts with writing the
{\em Hamiltonian representation\/} for this system, i.e.,
\begin{equation}
\frac{d}{dt}\left(f|_{m\left(t\right)}\right)=\left\{H,f\right\}|_{m\left(t\right)}\qquad \text{for any function }f\text{ on }M.
\notag\end{equation}
To do this one needs to find the Poisson bracket $ \left\{,\right\} $ and the function $ H $
(called the {\em Hamiltonian\/} of the equation). Note that $ H $ and $ \left\{,\right\} $ uniquely
determine the initial vector field $ v $.

Additionally, one needs to find a large enough independent
collection of functions $ H_{i} $ on $ M $ which all commute with each other w.r.t.~
$ \left\{,\right\} $ and such that $ H $ can be expressed as a function of $ H_{i} $. Alternatively,
one starts with a given bracket $ \left\{,\right\} $ and a function $ H $, then the problem is
to find the family $ H_{i} $. In fact, given the family $ H_{i} $, one can take as $ H $
any function of $ H_{i} $.

Thus to construct an integrable system a key problem is to find a
large family of independent functions $ H_{i} $ {\em in involution}, i.e., such that
$ \left\{H_{i},H_{j}\right\}=0 $. The anchored Lenard scheme is a particular algorithm to
construct such a family on a bihamiltonian manifold.

Start with a way to find many functions in involutions, not
necessarily independent. Most statements below are applicable both in
$ C^{\infty} $-geometry and in analytic geometry. In such cases we state the smooth
variant only, for the corresponding analytic statement one needs to
substitute $ {\mathbb R}{\mathbb P}^{1} $ by $ {\mathbb C}{\mathbb P}^{1} $.

\begin{definition} \label{def48.05}\myLabel{def48.05}\relax  Consider a bihamiltonian structure on $ M $, an open
subset $ {\mathcal U}\subset{\mathbb R}{\mathbb P}^{1} $ and a smooth function on $ {\mathcal U}\times M $, $ \left(\lambda,m\right) \mapsto C_{\lambda}\left(m\right) $. Consider this
function as a family $ C_{\lambda} $, $ \lambda\in{\mathcal U} $, of functions on $ M $. Call $ C_{\lambda} $ a
$ \lambda $-{\em Casimir family\/} on $ M $ if $ C_{\lambda_{1}:\lambda_{2}} $ is a Casimir function for $ \lambda_{1}\left\{,\right\}_{1}+\lambda_{2}\left\{,\right\}_{2} $
for $ \left(\lambda_{1}:\lambda_{2}\right)\in{\mathcal U} $. \end{definition}

\begin{proposition} \label{prop48.10}\myLabel{prop48.10}\relax  Consider a bihamiltonian structure $ M $ and a point
$ m_{0}\in M $. Fix $ \lambda^{0}\in{\mathbb R}{\mathbb P}^{1} $. Suppose that the corank of the bracket $ \lambda_{1}\left\{,\right\}_{1}+\lambda_{2}\left\{,\right\}_{2} $ at
$ m $ is $ r\in{\mathbb Z} $ for any $ m $ near $ m_{0} $ and any $ \lambda_{1}:\lambda_{2} $ near $ \lambda^{0} $. Then there is a
neighborhood $ U\times{\mathcal U} $ of $ \left(m_{0},\lambda^{0}\right)\in M\times{\mathbb R}{\mathbb P}^{1} $ and $ r $ families $ C_{t,\lambda} $, $ 1\leq t\leq r $, $ \lambda\in{\mathcal U} $, of
functions on $ U $ such that
\begin{enumerate}
\item
for any given $ t $, $ 1\leq t\leq r $, $ C_{t,\lambda} $ is a $ \lambda $-Casimir family on $ U $, and
\item
for any given $ \lambda\in{\mathcal U} $ the functions $ C_{t,\lambda} $, $ 1\leq t\leq r $, are independent.
\end{enumerate}
\end{proposition}

\begin{proof} Let $ \lambda^{0}=\left(\lambda_{1}^{0}:\lambda_{2}^{0}\right)\in{\mathbb R}{\mathbb P}^{1} $. Consider the symplectic leaf of
$ \lambda_{1}^{0}\left\{,\right\}_{1}+\lambda_{2}^{0}\left\{,\right\}_{2} $ passing through $ m_{0} $. The codimension of this leaf is $ r $,
fix a transversal manifold $ N $ of dimension $ r $, and coordinate functions $ c_{t} $,
$ 1\leq t\leq r $, on this manifold. Obviously, if $ \left(\lambda_{1}:\lambda_{2}\right) $ is close to $ \left(\lambda_{1}^{0}:\lambda_{2}^{0}\right) $ and
$ m $ is close to $ m_{0} $, then the symplectic leaf of $ \lambda_{1}\left\{,\right\}_{1}+\lambda_{2}\left\{,\right\}_{2} $ passing
through $ m $ intersects $ N $ in exactly one point, and this point depends
smoothly on $ \lambda=\lambda_{1}:\lambda_{2} $ and $ m $.

Thus there is exactly one Casimir function $ C_{t,\lambda} $ for $ \lambda_{1}\left\{,\right\}_{1}+\lambda_{2}\left\{,\right\}_{2} $
which coincides with $ c_{t} $ when restricted to $ N $. Obviously, it satisfies the
conditions of the proposition. \end{proof}

\begin{lemma} Consider two $ \lambda $-Casimir families $ C_{\lambda} $, $ \lambda\in{\mathcal U} $, and $ C'_{\lambda} $, $ \lambda\in{\mathcal U}' $, on $ M $.
Then
\begin{equation}
\left\{C_{\lambda},C'_{\mu}\right\}_{1}=\left\{C_{\lambda},C'_{\mu}\right\}_{2}=0,\qquad \lambda\in{\mathcal U},\quad \mu\in{\mathcal U}'.
\notag\end{equation}
\end{lemma}

\begin{proof} To simplify notations assume $ \infty\notin{\mathcal U} $. Let $ \left\{,\right\}^{\lambda}\buildrel{\text{def}}\over{=}\lambda\left\{,\right\}_{1}+\left\{,\right\}_{2} $. Since
$ \left\{C_{\lambda},f\right\}^{\lambda}=\left\{f,C'_{\mu}\right\}^{\mu}=0 $ for any function $ f $, we see that $ \left\{C_{\lambda},C'_{\mu}\right\}^{\nu}=0 $ if $ \nu=\lambda $ or
$ \nu=\mu $. Since $ \left\{,\right\}^{\nu} $ is linear in $ \nu $, $ \left\{C_{\lambda},C'_{\mu}\right\}^{\nu}=0 $ for any $ \nu $ as far as $ \lambda\not=\mu $. On
the other hand, the same identity is true for $ \lambda=\mu $ by continuity in $ \lambda $. \end{proof}

Proposition~\ref{prop48.10} provides a way to obtain a giant collection
of functions which commute with each other w.r.t.~{\em both\/} the brackets. Out
of this huge collection of functions on $ M $ only a finite number of
functions are independent (since this number is bounded by the dimension
of the manifold). One needs a way to extract a finite subset out of this
continuum. The anchored Lenard scheme provides such a way, moreover, it
allows one to find this small collection without actually finding the whole
continuum of Casimir functions.

The idea of the anchored Lenard scheme is to put $ \lambda^{0}=\infty $ and write a
formal series in\footnote{To follow the standard description of Lenard scheme, we use expansion in
formal variable $ \lambda^{-1} $, though by exchanging $ \left\{,\right\}_{1} $ and $ \left\{,\right\}_{2} $ one might use
more natural expansion in $ \lambda $.} $ \lambda^{-1} $ for a $ \lambda $-Casimir family $ C_{\lambda} $ defined near $ \lambda_{0} $:
\begin{equation}
C_{\lambda}=H_{0}+\lambda^{-1}H_{1}+\lambda^{-2}H_{2}+\dots .
\notag\end{equation}
Obviously, commutativity of Casimir functions implies $ \left\{H_{i},H_{j}\right\}_{1}=\left\{H_{i},H_{j}\right\}_{2}=0 $
for any $ i,j $. On the other hand, the condition
\begin{equation}
\left\{H_{0}+\lambda^{-1}H_{1}+\lambda^{-2}H_{2}+\dots ,f\right\}_{1}+\lambda^{-1}\left\{H_{0}+\lambda^{-1}H_{1}+\lambda^{-2}H_{2}+\dots ,f\right\}_{2}=0,
\notag\end{equation}
which describes the formal-variables analogue of the condition on a
$ \lambda $-Casimir family, can be written as
\begin{enumerate}
\item
function $ H_{0} $ is a Casimir function for $ \left\{,\right\}_{1} $;
\item
for any function $ f $ on $ M $
\begin{equation}
\left\{H_{i},f\right\}_{2}=-\left\{H_{i+1},f\right\}_{1}.
\label{equ48.20}\end{equation}\myLabel{equ48.20,}\relax 

\end{enumerate}
\begin{remark} It is easy to see that given $ H_{i} $, the relation~\eqref{equ48.20} is
equivalent for a system of equations of the form $ dH_{i+1}|_{{\mathcal L}}=\omega_{{\mathcal L}} $, here $ {\mathcal L} $ runs
over symplectic leaves of $ \left\{,\right\}_{1} $, and $ \omega_{{\mathcal L}} $ is a $ 1 $-form on $ {\mathcal L} $ which is
determined by $ H_{i} $. In particular, if $ \left\{,\right\}_{1} $ has a constant rank, then
~\eqref{equ48.20} has a local solution iff all the forms $ \omega_{{\mathcal L}} $ are closed.

If so, one can find $ H_{i+1} $ by integrating $ \omega_{{\mathcal L}} $. Thus if a solution to
~\eqref{equ48.20} exists, it is easy to find. One can also see why in Lenard
scheme one takes $ \lambda^{0}=\infty $: in applications $ \left\{,\right\}_{1} $ is much simpler than $ \left\{,\right\}_{2} $ or
any other combination $ \lambda_{1}\left\{,\right\}_{1}+\lambda_{2}\left\{,\right\}_{2} $, thus taking $ \lambda^{0}=\infty $ simplifies the
integration of relations~\eqref{equ48.20}. \end{remark}

\begin{definition} \label{def48.25}\myLabel{def48.25}\relax  Consider a formal series $ H_{0}+\lambda^{-1}H_{1}+\lambda^{-2}H_{2}+\dots $ in $ \lambda^{-1} $
with $ H_{i} $ being functions on $ M $. Call it a {\em formal\/} $ \lambda $-{\em family\/} on $ M $ if the
sequence $ H_{k} $ satisfies the recurrence relation~\eqref{equ48.20}. Call this formal
$ \lambda $-family {\em anchored\/} if $ H_{0} $ is a Casimir function for $ \left\{,\right\}_{1} $. \end{definition}

\begin{proposition} Given two anchored formal $ \lambda $-families $ H_{0}+\lambda^{-1}H_{1}+\lambda^{-2}H_{2}+\dots $
and $ H'_{0}+\lambda^{-1}H'_{1}+\lambda^{-2}H'_{2}+\dots $ one has $ \left\{H_{i},H'_{j}\right\}_{1}=\left\{H_{i},H'_{j}\right\}_{2}=0 $ for any $ i $ and $ j $.
\end{proposition}

\begin{proof} Put $ H_{i}=H'_{i}=0 $ for $ i<0 $. This makes~\eqref{equ48.20} applicable for $ i<0 $
too. For any $ i $ and $ j $
\begin{equation}
\left\{H_{i},H'_{j}\right\}_{1}=-\left\{H_{i-1},H'_{j}\right\}_{2}=\left\{H'_{j},H_{i-1}\right\}_{2}=-\left\{H'_{j+1},H_{i-1}\right\}_{1}=\left\{H_{i-1},H'_{j+1}\right\}_{1}.
\notag\end{equation}
Repeating this process $ i+1 $ times, one gets $ \left\{H_{i},H'_{j}\right\}_{1}=\left\{H_{-1},H'_{i+j+1}\right\}_{1}=0 $.
Moreover, $ \left\{H_{i},H'_{j}\right\}_{2}=-\left\{H_{i+1},H'_{j}\right\}_{1}=0 $. \end{proof}

If one considers one chain of solutions to~\eqref{equ48.20}, then the
anchoring condition may be dropped:

\begin{amplification} Given a formal $ \lambda $-family $ H_{0}+\lambda^{-1}H_{1}+\lambda^{-2}H_{2}+\dots $ one has
$ \left\{H_{i},H_{j}\right\}_{1}=\left\{H_{i},H_{j}\right\}_{2}=0 $ for any $ i $ and $ j $. \end{amplification}

\begin{proof} For any $ i\geq1 $ and $ j\geq0 $ one gets $ \left\{H_{i},H_{j}\right\}_{1}=\left\{H_{i-1},H_{j+1}\right\}_{1} $ again.
Repeating this several times, one can decrease $ |i-j| $ until it becomes 0
or 1 (depending on $ i-j $ being even or odd). If $ i-j $ is even, use
$ \left\{H_{k-l},H_{k+l}\right\}_{1}=\left\{H_{k},H_{k}\right\}_{1}=0 $, if $ i-j $ is odd, use
$ \left\{H_{k-l+1},H_{k+l}\right\}_{1}=\left\{H_{k+1},H_{k}\right\}_{1}=-\left\{H_{k},H_{k}\right\}_{2}=0 $. Thus $ \left\{H_{i},H_{j}\right\}_{1} $ is always 0,
moreover, $ \left\{H_{i},H_{j}\right\}_{2}=-\left\{H_{i+1},H_{j}\right\}_{1}=0 $. \end{proof}

\begin{lemma} \label{lm48.30}\myLabel{lm48.30}\relax  Suppose that conditions of Proposition~\ref{prop48.10} are
satisfied. Fix $ n\geq0 $. Given $ m_{0}\in M $ and any sequence of functions $ H_{0},\dots ,H_{n} $ on
$ M $ such that $ H_{0} $ is a Casimir function for $ \left\{,\right\}_{1} $, and $ H_{k} $ satisfy equations
~\eqref{equ48.20} for $ i=0,\dots ,n-1 $, there exists a neighborhood $ U $ of $ \left(m_{0},\infty\right)\in M\times{\mathbb R}{\mathbb P}^{1} $
and a $ \lambda $-Casimir family $ C_{\lambda}\left(m\right) $ defined for $ \left(m,\lambda\right)\in U $ such that
\begin{equation}
C_{\lambda}=H_{0}+\lambda^{-1}H_{1}+\lambda^{-2}H_{2}+\dots +\lambda^{-n}H_{n}+o\left(\lambda^{-n}\right).
\notag\end{equation}

In particular, there is a function $ H_{n+1} $ defined near $ m_{0} $ which solves
~\eqref{equ48.20} for $ i=n $. \end{lemma}

\begin{proof} To simplify notations, suppose $ r=1 $ (the case of general $ r $ is
absolutely parallel). Then $ H_{0} $ is defined uniquely up to a change
$ H_{0}'=\varphi_{0}\left(H_{0}\right) $. Additionally, given $ H_{i} $, Equation~\eqref{equ48.20} determines $ H_{i+1} $ up
to a change $ H_{i+1}'=H_{i+1}+\varphi_{i+1}\left(H_{0}\right) $.

Moreover, the Taylor series for $ C_{1,\lambda} $ provides {\em one\/} solution to the
recursion relations~\eqref{equ48.20}. Since the change of the form $ H_{0}'=\varphi_{0}\left(H_{0}\right) $
corresponds to a change of the form $ c_{1}'=\varphi_{0}\left(c_{1}\right) $ in the proof of
Proposition~\ref{prop48.10}, we conclude that there is a locally defined
solution to the recursion relations~\eqref{equ48.20} for any initial data $ H_{0} $
which is a Casimir function for $ \left\{,\right\}_{1} $.

Next, proceed by induction in $ n $. To do the step of induction, it is
enough to prove the following statement: given a $ \lambda $-Casimir family $ C_{\lambda} $ near
$ \lambda=\infty $ such that $ C_{\infty}=H_{0} $, and given any function $ \varphi_{n}\left(h\right) $ of one variable defined
in a neighborhood of $ h=H_{0}\left(m_{0}\right) $ one can find another $ \lambda $-Casimir family $ C'_{\lambda} $
such that $ C'_{\lambda}-C_{\lambda} =\varphi_{n}\left(H_{0}\right)\lambda^{-n}+o\left(\lambda^{-n}\right) $. Putting $ C'_{\lambda}=C_{\lambda}+\varphi\left(C_{\lambda}\right)\lambda^{-n} $ finishes the
proof in the case $ r=1 $. \end{proof}

The following statement is obvious:

\begin{lemma} Suppose that $ r=1 $ and a sequence $ \left(H_{i}\right) $ satisfies conditions of
Lemma~\ref{lm48.30}. If $ H_{k} $ depends functionally on $ H_{0},\dots ,H_{k-1} $, then $ H_{l} $
depends functionally on $ H_{0},\dots ,H_{k-1} $ for any $ l $ such that $ k\leq l\leq n $. \end{lemma}

This shows that a maximal independent subset of the sequence $ \left(H_{l}\right) $
can be chosen to be the starting subsequence. The situation in the case
$ r>1 $ is slightly more complicated, however, it is easy to show that

\begin{proposition} \label{prop48.50}\myLabel{prop48.50}\relax  Consider a maximal collection $ H_{0}^{\left(1\right)},\dots ,H_{0}^{\left(r\right)} $ of
independent Casimir functions for $ \left\{,\right\}_{1} $ near $ m_{0}\in M $. Let $ H_{i}^{\left(t\right)} $, $ t=1,\dots ,r $,
$ i\geq0 $, be solutions to recursion relations~\eqref{equ48.20} with $ H_{0}^{\left(t\right)} $ as the
initial data. Then there are numbers $ k_{1},\dots ,k_{r} $ such that the collection
$ \left\{H_{i}^{\left(t\right)}\right\} $, $ 1\leq t\leq r $, $ 0\leq i\leq k_{t} $, is functionally independent, and all functions
$ H_{i}^{\left(t\right)} $, $ 1\leq t\leq r $, $ i\geq0 $, depend functionally on this collection. \end{proposition}

\begin{definition} {\em Anchored Lenard scheme\/} of finding a large family of
functions on a bihamiltonian structure which mutually commute w.r.t.~both
brackets consists of two steps: first one finds a maximal independent
collection of Casimir functions for the bracket $ \left\{,\right\}_{1} $, then one solves
recurrence relations~\eqref{equ48.20} with these functions as initial data until
new functions start depend on the old ones. \end{definition}

In fact it is not necessary to consider many chains of solutions of
recurrence relations:

\begin{amplification} In conditions of Proposition~\ref{prop48.10} suppose that the
bihamiltonian structure on $ M $ is analytic. Then there is a sequence of
functions $ H_{0},\dots ,H_{n} $ defined near a given point $ m_{0}\in M $ such that
\begin{enumerate}
\item
function $ H_{0} $ is a Casimir function for $ \left\{,\right\}_{1} $;
\item
functions $ \left(H_{i}\right) $ satisfy the recurrence relation~\eqref{equ48.20};
\item
functions $ \left(H_{i}\right) $ are independent, and for any $ 1\leq t\leq r $ and $ \lambda $ near $ \infty $ the
function $ C_{t,\lambda} $ of Proposition~\ref{prop48.50} depends on $ \left(H_{0},\dots ,H_{n}\right) $.
\end{enumerate}
\end{amplification}

\begin{proof} Since the Taylor series for $ C_{t,\lambda} $ in $ \lambda^{-1} $ converge, it is enough
to show that the Taylor coefficients for $ C_{t,\lambda} $ depend on $ \left(H_{0},\dots ,H_{n}\right) $.
Fix numbers $ \alpha_{t} $, $ 2\leq t\leq r $. Let $ n=\sum_{t=1}^{r}k_{t}+r-1 $, and
\begin{equation}
C_{\lambda}= C_{1,\lambda}+\alpha_{2}\lambda^{-k_{1}+1}C_{2,\lambda}+\dots +\alpha_{r}\lambda^{-k_{1}-\dots -k_{r-1}+1}C_{r,\lambda}.
\notag\end{equation}
Obviously, this is a $ \lambda $-Casimir family.

It is easy to show that for generic values of $ \alpha_{2},\dots ,\alpha_{r} $ the first
$ n+1 $ Taylor coefficients $ H_{0},\dots ,H_{n} $ of $ C_{\lambda} $ are independent, which finishes
the proof. \end{proof}

\begin{remark} Since the functions $ H_{i}^{\left(t\right)} $ of the anchored Lenard scheme are
obtained by doing manipulations (taking Taylor coefficients) with Casimir
functions, they can be pushed down to the web $ {\mathcal B}_{M} $ of $ M $. Thus they should
be considered as action functions on $ M $ (see Section~\ref{h02}).

In interesting cases (see Section~\ref{h55} and \cite{Pan99Ver}) the functions
$ H_{i}^{\left(t\right)} $ provide a local coordinate system on $ {\mathcal B}_{M} $. (This shows that in fact
$ {\mathcal B}_{U} $ is a smooth manifold if $ U $ is a small subset of $ M $.) In these cases the
submanifolds $ \left\{H_{i}^{\left(t\right)}=\operatorname{const}_{i}^{\left(t\right)} \mid i\geq0\right\} $ carry a natural local affine
structure, thus one can find a complementary set of {\em angle variables\/} $ \varphi_{j} $
such that functions $ \left\{H_{i},\varphi_{j}\right\}_{k}\buildrel{\text{def}}\over{=}c_{\text{ijk}} $ depends on $ H_{l} $ only\footnote{Another problem is to find such change-of-variables in action variables
$ \overset{\,\,{}_\circ}{H}_{i}=\overset{\,\,{}_\circ}{H}_{i}\left(H_{0},\dots \right) $ that the corresponding functions $ \overset{\,\,{}_\circ}{c}_{\text{ijk}} $ become as simple as
possible. As Theorem~\ref{th47.13} shows, in general it is {\em not possible\/} to
make all $ \overset{\,\,{}_\circ}{c}_{\text{ijk}} $ into constants. However, it is obviously possible for
bihamiltonian structures with constant coefficients, thus for flat
structures.}.

\end{remark}

\begin{example} \label{ex48.80}\myLabel{ex48.80}\relax  Consider the bihamiltonian structure defined by
~\eqref{equ47.10}. In this case $ r=1 $, and Casimir functions are functions of $ x $ and
$ y $ only. Thus $ H_{i} $ are functions of $ x $ and $ y $ too. Moreover, one can write an
explicit formula for $ H_{i} $.

Indeed, let $ \Phi\left(x,y\right)=\frac{\partial f/\partial x}{\partial f/\partial y} $. Obviously, the symplectic leaves
for $ \left\{,\right\}_{1}+\lambda^{-1}\left\{,\right\}_{2} $ can be described as surfaces $ \left\{\left(x,y,z\right) \mid y=\Psi\left(x\right)\right\} $, here $ \Psi $
is a solution of the ODE
\begin{equation}
\frac{d\Psi}{dx}=-\lambda^{-1}\Phi\left(x,\Psi\right)
\label{equ48.81}\end{equation}\myLabel{equ48.81,}\relax 
Given $ \left(x_{0},y_{0}\right) $ which is close to (0,0), let $ \Psi_{\lambda,x_{0},y_{0}}\left(x\right) $ be the solution of
~\eqref{equ48.81} which passes through the point $ \left(x_{0},y_{0}\right) $. Let
$ F_{\lambda}\left(x_{0},y_{0}\right)\buildrel{\text{def}}\over{=}\Psi_{\lambda,x_{0},y_{0}}\left(0\right) $. Obviously, $ F_{\lambda}\left(x,y\right) $ is well-defined for large
$ |\lambda| $ and small $ \left(x,y\right) $. Moreover, $ F_{\lambda} $ is a Casimir function for $ \left\{,\right\}_{1}+\lambda^{-1}\left\{,\right\}_{2} $.

Taking Laurent coefficients of $ F_{\lambda} $ near $ \lambda=\infty $, one obtains functions $ H_{i} $
from the anchored Lenard scheme. Obviously,
\begin{equation}
H_{0}\left(x,y\right) = y,\qquad H_{1}\left(x,y\right)=\int_{0}^{x}\frac{\partial f\left(t,y\right)/\partial t}{\partial f\left(t,y\right)/\partial y}dt=x+o\left(x\right).
\notag\end{equation}
This implies that all other functions $ H_{i} $ depend on $ H_{0} $ and $ H_{1} $. One can see
that $ z $ provides an example of an angle variable, and any other angle
variable can be written as $ a\left(x,y\right)z+b\left(x,y\right) $ with arbitrary $ a\left(x,y\right) $ and
$ b\left(x,y\right) $.

\end{example}

Summing up, we obtain

\begin{proposition} The bihamiltonian structure $ M_{f} $ given by~\eqref{equ47.10} is
completely integrable by the anchored Lenard scheme. If $ \varphi $ satisfies
~\eqref{equ47.20} and $ \varphi\left(x,y\right)\not\equiv x+y $, then $ M_{\varphi} $ is not flat. \end{proposition}

\begin{remark} It is possible to provide similar examples of homogeneous but
not flat bihamiltonian structures of any given type. In Section~\ref{h55} we
will see that all these structures are completely integrable by the
anchored Lenard scheme. In the case of type $ \left(2k-1\right) $, $ k\in{\mathbb N} $, one can write
such a bihamiltonian structure\footnote{In a slightly different language such bihamiltonian structures were
described in \cite{GelZakhWeb} and \cite{GelZakh93}.} based on $ k-1 $ functions
$ \varphi_{1}\left(x,y\right),\dots ,\varphi_{k-1}\left(x,y\right) $ of two complex variables (though one cannot do it
as explicitly as in~\eqref{equ47.10}). Any two of these bihamiltonian structures
are not locally isomorphic, thus only one of them (for any given $ k\in{\mathbb N} $) is
flat. What is very surprising is that (apparently) they did not appear in
examples of integrable systems arising in problems of mathematical
physics. \end{remark}

\begin{remark} \label{rem48.91}\myLabel{rem48.91}\relax  Let us point out the relation of the anchored Lenard
scheme with the {\em algebraic Zakharov\/}--{\em Shabat scheme\/} of \cite{DriSok84Alg}.
Recall how the latter scheme works. Given a Poisson structure $ \left\{,\right\} $ on $ M $ and a
function $ H $ on $ M $, define the {\em Hamiltonian vector field\/} $ {\mathcal V}_{H} $ of $ H $ by the
identity $ {\mathcal V}_{H}\cdot f=\left\{H,f\right\} $ for any function $ f $ on $ M $. Given two Poisson
structures $ \left\{,\right\}_{1} $, $ \left\{,\right\}_{2} $, one obtains two Hamiltonian vector fields $ {\mathcal V}_{H}^{\left(1\right)} $,
$ {\mathcal V}_{H}^{\left(2\right)} $. Note that the Hamiltonian vector field of $ H $ for the bracket
$ \lambda\left\{,\right\}_{1}+\left\{,\right\}_{2} $ is $ \lambda{\mathcal V}_{H}^{\left(1\right)}+{\mathcal V}_{H}^{\left(2\right)} $.

Consider a family $ {\mathcal H}_{\lambda} $ of function on $ M $ which depends polynomially on
$ \lambda $. Say that a vector field $ V $ on $ M $ is {\em associated\/} with $ {\mathcal H}_{\lambda} $ if
$ V=\lambda{\mathcal V}_{{\mathcal H}_{\lambda}}^{\left(1\right)}+{\mathcal V}_{{\mathcal H}_{\lambda}}^{\left(2\right)} $ for any $ \lambda $, in particular, for the association to hold,
the right-hand side should not depend on $ \lambda $. The associated vector fields
are the central tool of the algebraic Zakharov--Shabat scheme. In many
examples such vector fields commute, and are plentiful enough to
completely integrate the bihamiltonian structure.

To explain this phenomenon write $ {\mathcal H}_{\lambda}=\sum_{k=0}^{K}H_{k}\lambda^{K-k} $. Clearly, the finite
sequence $ \left(H_{k}\right) $ satisfies the same conditions as an anchored $ \lambda $-series:
function $ H_{0} $ is a Casimir function for $ \left\{,\right\}_{1} $, and the relation~\eqref{equ48.20}
holds. Moreover, any vector field in the span of Hamiltonian vector
fields of $ \left(H_{k}\right) $ w.r.t.~any Poisson structure $ \lambda_{1}\left\{,\right\}_{1}+\lambda_{2}\left\{,\right\}_{2} $ can be written
as an associated vector field of an appropriate family.

Thus one can consider the algebraic Zakharov--Shabat scheme as
a different formulation of the anchored Lenard scheme. \end{remark}

\section{Lenard-integrable structures }\label{h55}\myLabel{h55}\relax 

Here we show that the class of bihamiltonian structures for which
the anchored Lenard scheme gives ``many'' functions in involution coincides
with the class of homogeneous structures. In fact, since our approach to
Lenard scheme is based on a formal analogue of $ \lambda $-Casimir families, the
result of this section are closely related to ones in \cite{Bol91Com} (compare
with discussion of ``completeness'' in \cite{Pan99Ver}).

\begin{definition} The {\em action dimension\/} of a Poisson structure $ \left(M,\left\{,\right\}_{1}\right) $ of
constant corank $ r $ is $ \frac{\dim  M+r}{2} $. The action dimension of an arbitrary
Poisson structure on $ M $ at $ m_{0}\in M $ is the minimum action dimension of open
subsets $ U\subset M $ which contain $ m_{0} $ in its closure, and such that the Poisson
structure is of constant corank on $ U $. \end{definition}

This definition gives a {\em lower\/} bound on the number of functions in
involution which are enough to completely integrate the dynamical system
on $ M $ given by some Hamiltonian $ H $. Indeed, in the case of constant corank
$ r $ one needs $ r $ functions to disambiguate symplectic leaves, and $ \frac{\dim 
M-r}{2} $ functions to provide action variables inside the leaves.

To do the same in the case of a bihamiltonian structure, introduce

\begin{definition} The {\em action dimension\/} of a complex vector space $ V $ with two
skewsymmetric bilinear pairings is $ \frac{\dim  V +r}{2} $, here $ r $ is the number of
Kronecker blocks of $ V $. \end{definition}

\begin{definition} The {\em action dimension\/} at $ m_{0}\in M $ of a bihamiltonian structure on
$ M $ is the lower limit of action dimensions of\footnote{If $ M $ is analytic, one should consider $ {\mathcal T}_{m}^{*}M $ instead of $ {\mathcal T}_{m}^{*}M\otimes{\mathbb C} $.} $ {\mathcal T}_{m}^{*}M\otimes{\mathbb C} $ for $ m \to m_{0} $. \end{definition}

Note that the number of Kronecker blocks of a pair of skewsymmetric
pairings $ \left(,\right)_{1} $, $ \left(,\right)_{2} $ is equal to $ \min _{\lambda_{1},\lambda_{2}}\dim  \operatorname{Ker} \left(\lambda_{1}\left(,\right)_{1}+\lambda_{2}\left(,\right)_{2}\right) $, here $ \operatorname{Ker} $
denotes null-space of the pairing. Thus the action dimension of a
bihamiltonian structure provides a {\em lower\/} bound on the number of functions
in involution necessary to completely integrate the structure w.r.t.~{\em at
least one particular\/} Poisson structure of the form $ \lambda_{1}\left\{,\right\}_{1}+\lambda_{2}\left\{,\right\}_{2} $ on an
open subset of $ M $ near $ m_{0} $.

\begin{definition} Call a bihamiltonian structure on $ M $ {\em Lenard-integrable\/} at
$ m_{0}\in M $ if the number of independent functions provided by the anchored
Lenard scheme in an appropriate neighborhood of $ m_{0} $ coincides with the
action dimension of $ M $ at $ m_{0} $.

Call a bihamiltonian structure on $ M $ {\em strictly Lenard-integrable\/} at $ m_{0} $
if it is Lenard-integrable at $ m_{0} $ and the sequences $ H_{i}^{\left(t\right)} $ of the anchored
Lenard scheme can be continued for $ i>k_{t} $ as well. \end{definition}

\begin{remark} Recall that Section~\ref{h48} describes the anchored Lenard scheme
as a formal-series counterpart of $ \lambda $-Casimir families. For this
description to work one needs to assume some constant rank conditions, as
in Proposition~\ref{prop48.10}. The condition of Proposition~\ref{prop48.10} was
not very restrictive, since one could achieve it by a small deformation
of $ \left(m_{0},\lambda^{0}\right) $. However, in the anchored Lenard scheme $ \lambda^{0} $ is fixed to be $ \infty $,
thus the restriction of Proposition~\ref{prop48.10} is in fact non void. Thus
Lemma~\ref{lm48.30} {\em does not\/} imply that any Lenard-integrable structure is
strictly Lenard-integrable. \end{remark}

\begin{theorem} \label{th55.50}\myLabel{th55.50}\relax  If a bihamiltonian system on $ M $ is strictly
Lenard-integrable at $ m_{0}\in M $, then it is homogeneous on an open subset $ U $ of
$ M $ with $ m_{0} $ being in the closure of $ U $. \end{theorem}

\begin{proof} Indeed, if a structure is Lenard-integrable at $ m_{0} $, then it is
also Lenard-integrable at $ m $ for $ m $ in an appropriate open subset of $ M $. It
is easy to show that by decreasing this subset $ U $ one may assume that at
any point $ m\in U $ the sizes of Kronecker blocks of the pair of pairings on
$ {\mathcal T}_{m}^{*}M\otimes{\mathbb C} $ are the same.

Functions $ H_{i}^{\left(t\right)} $, $ 1\leq t\leq r $, $ 0\leq i\leq k_{t} $, given by the anchored Lenard scheme
provide a mapping $ {\mathbit H}\colon U \to {\mathbb R}^{K} $, $ K=\sum_{t}\left(k_{t}+1\right) $. Decreasing $ U $ yet more, we may
assume that the differential of this mapping is of constant rank $ K $
(recall that components of $ {\mathbit H} $ are independent). Fix a point $ m\in U $ and $ t $,
$ 1\leq t\leq r $. Let $ \beta_{i}=dH_{i}^{\left(t\right)}|_{m}\in{\mathcal T}_{m}^{*}M $. Let $ W_{{\mathbb R}}={\mathcal T}_{m}^{*}M $, $ W=W_{{\mathbb R}}\otimes{\mathbb C} $. By Equation~\eqref{equ48.20},
$ \beta_{0} $ is in the null space of pairing $ \left(,\right)_{1} $ on $ W $, and $ \left(\beta_{i},w\right)_{2}=\left(\beta_{i+1},w\right)_{1} $ for
any $ w\in W $.

An immediate check shows that if $ W\simeq{\mathcal J}_{2k,\mu} $, $ \mu\in{\mathbb C}{\mathbb P}^{1} $, then $ \beta_{i}=0 $, $ i=0,\dots . $
Similarly, if $ W\simeq{\mathcal K}_{2k-1} $, then all vectors $ \beta_{i} $ are in the subspace
$ W_{1}=\left<{\mathbit w}_{0},{\mathbit w}_{2},\dots ,{\mathbit w}_{2k-2}\right> $ of $ {\mathcal K}_{2k-1} $. The dimension of this subspace is $ k $ (it
is the same subspace which appears in a similar context in Lemma
~\ref{lm25.20}). In general case, taking a decomposition of $ W $ into a sum of
indecomposable components, one can see that all vectors $ \beta_{i} $ are in the sum
of Kronecker blocks of $ W $, moreover, they are in a direct sum of subspaces
$ W_{1} $ for these blocks.

This shows that $ N\leq\sum_{t=1}^{r}k_{r} $, here $ {\mathcal K}_{2k_{t}-1} $, $ 1\leq t\leq r $, are Kronecker blocks
of $ W $. The restriction on the action dimension shows that $ \dim 
W\leq\sum_{t=1}^{r}\left(2k_{r}-1\right) $, thus $ W $ has no Jordan blocks, which finishes the proof of
the theorem. \end{proof}

\begin{proposition} \label{prop55.60}\myLabel{prop55.60}\relax  Any homogeneous bihamiltonian structure is
strictly Lenard-integrable on small open subsets. \end{proposition}

\begin{proof} A tiny modification of the above proof together with
Proposition~\ref{prop48.10} imply this statement immediately. \end{proof}

This shows that the ``strict'' anchored Lenard scheme integrates
homogeneous structures and only them. Note that a linear combination of
brackets of homogeneous structure is never symplectic.

\begin{remark} The Lenard schemes of \cite{Mag78Sim,GelDor79Ham,KosMag96Lax}
differ from what we describe here, the difference being that they
consider non-anchored formal $ \lambda $-families. Though our condition is more
restrictive, note that in applications the Lenard scheme usually provides
an anchored $ \lambda $-family. Moreover, in non-symplectic cases there is no
simple way to find a non-anchored family, thus it is not obvious whether
non-anchored Lenard scheme may be used to integrate a system (unless
applied to the traces of powers of recursion operator). \end{remark}

\begin{remark} \label{rem55.80}\myLabel{rem55.80}\relax  The amount of our knowledge about classification of
bihamiltonian structures is not enough to describe finite-dimensional
Lenard-integrable system which are not strict. The situation is slightly
more promising if one consider non-strict structures for which {\em one\/}
anchored Lenard chain provides enough functions in involution.

In this case slightly more elaborate arguments than those in the
proof of Theorem~\ref{th55.50} show that there is an open subset $ U\subset M $ such that
at a point $ m $ of $ U $ the pair of brackets in $ {\mathcal T}_{m}^{*}M $ has one Jordan block only,
and this block is of the form $ {\mathcal J}_{2k,\infty} $. (The remaining blocks are
Kronecker.)

In particular, if at least one linear combination $ \lambda_{1}\left\{,\right\}_{1}+\lambda_{2}\left\{,\right\}_{2} $ is
symplectic near $ m_{0} $, then there are no Kronecker blocks. Thus the pairings
on $ {\mathcal T}_{m}^{*}M $ are isomorphic to $ {\mathcal J}_{2k,\infty} $ for any $ m\in U $.

Such symplectic structures were classified in \cite{Tur89Cla}, they turn
out to be flat (thus isomorphic to the natural bihamiltonian structure on
the dual space to the vector space $ {\mathcal J}_{2k,\infty} $). These are exactly the
structures for which the arguments of \cite{Mag78Sim} and \cite{GelDor79Ham} are
actually applicable to {\em anchored\/} formal $ \lambda $-families. It is again an
interesting question to find physically interesting bihamiltonian
structures of this form. \end{remark}

\begin{remark} Note asymmetry between Theorem~\ref{th55.50} and Proposition
~\ref{prop55.60}: one of them is applicable on small open neighborhoods of any
point, another on small open neighborhood of a {\em dense\/} collection of
points. Note that \cite{Pan99Ver} introduces a more general notion than
homogeneity: bihamiltonian structure is {\em complete\/} if the pairs of pairings
at $ {\mathcal T}_{m}^{*}M $ for any $ m\in M $ do not contain a Jordan block (thus the condition of
Kronecker blocks having the same sizes for all the points of $ M $ is
dropped). For complete structures Proposition~\ref{prop48.10} is applicable
for any point of $ M $, and it is easy to see that the following statement
holds: \end{remark}

\begin{amplification} The class of bihamiltonian structures which are strictly
Lenard-integrable at any point of $ M $ coincides with the class of complete
bihamiltonian structures. \end{amplification}

\section{Bihamiltonian Toda lattices }\label{h0}\myLabel{h0}\relax 

\begin{definition} The {\em open Toda lattice\/} (\cite{FadTakh87Ham}) is the
$ \left(2k+1\right) $-dimensional vector space $ V_{2k+1} $ over $ {\mathbb C} $ with coordinates $ v_{0},\dots ,v_{2k} $
and the two compatible Poisson brackets defined as follows. The bracket
$ \left\{,\right\}_{1} $ is defined by the condition $ \left\{v_{i},v_{j}\right\}=0 $ for $ |i-j|>1 $, and
\begin{equation}
\left\{v_{2l},v_{2l\pm1}\right\}_{1} = \mp v_{2l\pm1}.
\label{equ0.20}\end{equation}\myLabel{equ0.20,}\relax 
The bracket $ \left\{,\right\}_{2} $ is defined by the condition $ \left\{v_{i},v_{j}\right\}=0 $ for $ |i-j|>2 $, and
\begin{equation}
\begin{aligned}
\left\{v_{2l},v_{2l\pm1}\right\}_{2} & = \mp v_{2l}v_{2l\pm1},
\\
\left\{v_{2l},v_{2l+2}\right\}_{2} & = -2v_{2l+1}^{2},
\\
\left\{v_{2l-1},v_{2l+1}\right\}_{2} & = -\frac{1}{2}v_{2l-1}v_{2l+1},
\end{aligned}
\label{equ0.10}\end{equation}\myLabel{equ0.10,}\relax 
for all $ l $ such that the the left-hand sides make sense. \end{definition}

We denote a point of $ V_{2k+1} $ by $ {\mathbit v} $. Define transformation $ {\mathfrak T}_{\lambda} $, $ \lambda\in{\mathbb C} $, by
\begin{equation}
{\mathfrak T}_{\lambda}\colon V_{2k+1} \to V_{2k+1}\colon {\mathbit v} \mapsto {\mathbit v}+\lambda{\mathbit v}^{0},\qquad {\mathbit v}^{0}=\left(1,0,1,0,1,\dots ,0,1\right).
\label{equ0.15}\end{equation}\myLabel{equ0.15,}\relax 
Translating bracket $ \left\{,\right\}_{2} $ by the transformation $ {\mathfrak T}_{-\lambda} $, one obtains a Poisson
bracket $ \left\{,\right\}^{\left(\lambda\right)} $ which depends on a parameter $ \lambda $.

\begin{remark} \label{rem01.57}\myLabel{rem01.57}\relax  Note that for any $ i $, $ j $ the bracket $ \left\{v_{i},v_{j}\right\}_{2} $ depends
linearly on $ v_{2l} $, $ l=0,\dots ,k $, thus $ \left\{,\right\}^{\left(\lambda\right)} $ depends linearly on $ \lambda $. In fact
$ \left\{,\right\}^{\left(\lambda\right)} $ may be written as
\begin{equation}
\left\{,\right\}^{\left(\lambda\right)} =\lambda\left\{,\right\}_{1} + \left\{,\right\}_{2}.
\notag\end{equation}

One can use this remark to simplify the proof of compatibility and
Poisson property of brackets $ \left\{,\right\}_{1} $ and $ \left\{,\right\}_{2} $. Indeed, if we know that $ \left\{,\right\}_{2} $
is Poisson, then $ \left\{,\right\}^{\lambda} $ is Poisson, thus is $ \left\{,\right\}_{1} $ as a limit of $ \left\{,\right\}^{\left(\lambda\right)}/\lambda $. To
check that $ \left\{,\right\}_{2} $ is Poisson, one can use the symmetry of~\eqref{equ0.10} of the
form $ l \mapsto 2m\pm l $, so it is enough to check Jacobi identity for $ v_{0},v_{1},v_{2} $,
for $ v_{1},v_{2},v_{3} $, for $ v_{0},v_{1},v_{3} $, for $ v_{0},v_{2},v_{3} $, for $ v_{0},v_{2},v_{4} $, and for $ v_{1},v_{3},v_{5} $.

\end{remark}

\begin{definition} The {\em infinite Toda lattice\/} is the manifold with coordinates
$ v_{l} $, $ l\in{\mathbb Z} $, and the Poisson brackets\footnote{These brackets are well-defined on functions which depend on finite
number of coordinates $ v_{l} $ only.}~\eqref{equ0.20},~\eqref{equ0.10}. Considering
sequences $ v_{l} $ with period $ 2k $, one obtains a pair of well-defined Poisson
brackets on a $ 2k $-dimensional subvariety. Denote this bihamiltonian
structure by $ V_{2k} $, call it the {\em periodic Toda lattice}. \end{definition}

In Sections~\ref{h4} and~\ref{h10} we prove that open dense subsets of the
bihamiltonian structures $ V_{2k+1} $ and $ V_{2k} $ are Kronecker bihamiltonian
structure. In other words, in these sections we prove the following
theorems:

\begin{theorem} \label{th01.60}\myLabel{th01.60}\relax  The open Toda lattice (of dimension $ 2k-1 $) is a
bihamiltonian structure which is generically Kronecker of type $ \left(2k-1\right) $. \end{theorem}

\begin{theorem} \label{th01.70}\myLabel{th01.70}\relax  The periodic Toda lattice (of dimension 2k) is a
bihamiltonian structure which is generically Kronecker of type $ \left(2k-1,1\right) $.
\end{theorem}

\begin{remark} \label{rem0.20}\myLabel{rem0.20}\relax  Note that one can also consider a manifold $ \widetilde{V}_{2k} $ with
coordinates $ v_{0},\dots ,v_{2k-1} $ and brackets~\eqref{equ0.20},~\eqref{equ0.10}. It is also
bihamiltonian, but it is not a Kronecker structure, so it cannot be
described by the methods of this paper. Say, at a generic point both the
Poisson structures are in fact symplectic, while all linear combinations
of Poisson structures of a Kronecker structure are degenerated. While
this structure may be described by the means of \cite{Tur89Cla,Mag88Geo,%
Mag95Geo,McKeanPC,GelZakh93}, note that the in applications $ \widetilde{V}_{2k} $
appears not by itself, but as a reduction of the structure $ V_{2k+1} $ w.r.t.~
forgetting the variable $ v_{2k} $.

This supports the point of view from Section~\ref{h005} that Kronecker
structures are more important in applications than structures which may
be described in symplectic terms. \end{remark}

\section{Casimir families on the open Toda lattice }\label{h4}\myLabel{h4}\relax 

Apply the description of Section~\ref{h2} to the bihamiltonian Toda
structure. First, construct a family of would-be semi-Casimir functions
$ F_{\lambda} $, $ \lambda\in{\mathbb C} $.

Consider the inclusion $ \iota $ of $ V_{2k+1} $ into $ \operatorname{Mat}\left(k+1,k+1\right) $ which sends
$ \left(v_{0},\dots ,v_{2k}\right) $ to a symmetric $ 3 $-diagonal matrix with diagonal elements
$ \left(v_{0},v_{2},\dots ,v_{2k}\right) $ and over-diagonal elements $ \left(v_{1},v_{3},\dots ,v_{2k-1}\right) $. Taking
determinant of the resulting matrix, one obtains a polynomial function $ F_{0} $
on $ V_{2k+1} $.

Any proof of integrability of Toda lattice is based on the following
statement:

\begin{lemma} \label{lm4.05}\myLabel{lm4.05}\relax  The function $ F_{0} $ is Casimir, in other words, for any
function $ f $ on $ V_{2k+1} $ the Poisson bracket $ \left\{F_{0},f\right\}_{2} $ is identically 0. \end{lemma}

\begin{proof} Let $ d_{2m} $ be the determinant of the upper-left minor of
$ \iota\left({\mathbit v}\right) $ of size $ \left(m+1\right)\times\left(m+1\right) $. We need to show that $ \left\{v_{l},d_{2k}\right\}_{2}=0 $, $ 0\leq l\leq2k $. Let
us show that $ \left\{v_{l},d_{2m}\right\}_{2}=0 $, $ 0\leq l\leq2m $, $ m\leq k $.

Use induction in $ m $. Plugging the identity
\begin{equation}
d_{2m}=v_{2m}d_{2m-2}-v_{2m-1}^{2}d_{2m-4}
\notag\end{equation}
into $ \left\{v_{l},d_{2m}\right\}_{2} $ shows that the step of induction will work as far as
$ l\leq2m-4 $. On the other hand, due to obvious symmetry $ v_{t}\iff v_{2m-t} $ of brackets
~\eqref{equ0.10} and the determinant $ d_{2m} $, it is enough to check $ \left\{v_{l},d_{2m}\right\}_{2}=0 $ for
$ 0\leq l\leq m $. Moreover, if we know $ \left\{v_{l},d_{2m}\right\}_{2}=0 $ for $ 0\leq l\leq m-1 $, then we know it for
$ m+1\leq l\leq2m $, thus $ \left\{d_{2m},d_{2m}\right\}=\frac{\partial d_{2m}}{\partial v_{m}}\left\{v_{m},d_{2m}\right\} $. Since all these expressions
are polynomials in $ v_{i} $, and $ \frac{\partial d_{2m}}{\partial v_{m}}\not\equiv 0 $, one would be able to conclude
that $ \left\{v_{m},d_{2m}\right\}=0 $.

Thus the only relations to check are $ \left\{v_{l},d_{2m}\right\}_{2}=0 $ for $ 0\leq l\leq m-1 $ such
that $ 2m-l\leq3 $. This leaves only $ \left\{v_{0},d_{2}\right\} $ and $ \left\{v_{1},d_{4}\right\} $, which are easy to
check (using one step of induction for the latter one). \end{proof}

\begin{remark} In Remark~\ref{rem01.57} we used the fact that the right-hand sides
of~\eqref{equ0.10} are linear in variables $ v_{2l} $. The last sentence of the above
proof is the only other place were we use the particular form of
right-hand sides of~\eqref{equ0.10}. \end{remark}

Consider the translation $ {\mathfrak T} $ defined in~\eqref{equ0.15}. Motivated by the
above lemma, define $ F_{\lambda}\buildrel{\text{def}}\over{=}{\mathfrak T}_{\lambda}^{*}F $, $ \lambda\in{\mathbb C} $. By definition of $ \left\{,\right\}^{\left(\lambda\right)} $, the bracket
$ \left\{F_{\lambda},f\right\}^{\left(\lambda\right)} $ is identically 0 for any function $ f $. On the other hand, for any
given $ {\mathbit v}\in V_{2k+1} $ the function $ F_{\lambda}\left({\mathbit v}\right) $ of $ \lambda $ is the characteristic polynomial of
$ \iota\left({\mathbit v}\right) $. Thus the degree of $ F_{\lambda}\left({\mathbit v}\right)+\left(-1\right)^{k}\lambda^{k+1} $ in $ \lambda $ is $ k $. We obtain

\begin{proposition} The family $ \overset{\,\,{}_\circ}{F}_{\lambda}\left({\mathbit v}\right)\buildrel{\text{def}}\over{=}F_{\lambda}\left({\mathbit v}\right)+\left(-1\right)^{k}\lambda^{k+1} $ of functions on $ V_{2k+1} $
depends polynomially on $ \lambda $ with the degree being $ k $. For each $ \lambda $ the
function $ \overset{\,\,{}_\circ}{F}_{\lambda} $ is a Casimir function for the bracket $ \lambda\left\{,\right\}_{1}+\left\{,\right\}_{2} $. \end{proposition}

However, this proposition is not yet enough to put us in the context
of Theorem~\ref{th1.10}, since we do not know the dimension of the span of
$ d\overset{\,\,{}_\circ}{F}_{\lambda}|_{{\mathbit v}} $ for any given $ {\mathbit v} $ and variable $ \lambda $. To find this dimension, we need to
investigate the functions $ F_{\lambda} $ in more details.

Denote the set of polynomials of degree $ d $ in $ \lambda $ with the leading
coefficient $ \left(-1\right)^{d} $ by $ {\mathfrak P}_{d} $. Functions $ F_{\lambda} $ (considered as polynomials in $ \lambda $)
define a mapping $ F_{\bullet}\colon V_{2k+1} \to {\mathfrak P}_{k+1} $, $ {\mathbit v} \mapsto F_{\bullet}\left({\mathbit v}\right) $.

To describe the geometry of this mapping, associate with each
$ {\mathbit v}=\left(v_{i}\right)\in V_{2k+1} $ a finite sequence of polynomials $ C_{I_{p}} $ in $ \lambda $. First, construct
a partition of the set of even numbers $ \left\{0,2,\dots ,2k\right\} $: consider numbers
$ 2l+1 $ such that $ v_{2l+1}=0 $ as walls, they separate $ \left\{0,2,\dots ,2k\right\} $ into
continuous intervals $ I_{1},\dots ,I_{q} $, which we call {\em runs}. To each run
$ I_{p}=\left\{2l_{p},2l_{p}+2,\dots ,2l_{p+1}-2\right\} $ associate the characteristic polynomial $ C_{I_{p}} $ of
the corresponding principal minor (with columns and rows $ l_{p}+1,\dots ,l_{p+1} $)
of the matrix $ \iota\left({\mathbit v}\right) $. Obviously, $ \det \left(\iota\left({\mathbit v}\right)-\lambda\right) $ coincides with the product of
polynomials $ C_{I_{p}} $.

Call $ {\mathbit v}\in V_{2k+1} S $-{\em generic\/} if any two of polynomials $ C_{I_{p}} $ are mutually
prime. Non-$ S $-generic points form a submanifold of codimension 2: one of
$ v_{2l+1} $ should vanish, and two polynomials should have a common zero.

\begin{proposition} At an $ S $-generic point $ {\mathbit v}\in V_{2k+1} $ the mapping $ F_{\bullet}\colon V_{2k+1} \to {\mathfrak P}_{k+1} $
is a submersion\footnote{I.e., its derivative is an epimorphism.}. At non-$ S $-generic points it is not a submersion. \end{proposition}

\begin{proof} It is enough to consider the case when no $ v_{2l+1} $ vanishes.
Indeed, if we leave all the variables $ v_{m} $ except $ v_{2l+1} $ fixed, then $ \det \iota\left({\mathbit v}\right) $
is quadratic in $ v_{2l+1} $ without the linear term. Thus $ v_{2l+1}=0 $ implies
$ \frac{\partial\det }{\partial v_{2l+1}}=0 $. On the other hand, if $ v_{2l+1}=0 $, the matrix breaks into
two blocks, and the derivatives w.r.t.~other variables can be calculated
when we consider two blocks separately. Now the case when some $ v_{2l+1} $
vanish can be proved by induction using the following obvious

\begin{lemma} \label{lm4.12}\myLabel{lm4.12}\relax  The multiplication mapping $ {\mathfrak P}_{a}\times{\mathfrak P}_{b} \to {\mathfrak P}_{a+b} $ is a submersion
at $ \left(P_{1},P_{2}\right) $ iff $ P_{1} $ and $ P_{2} $ are mutually prime. \end{lemma}

In the case when all $ v_{2l+1}\not=0 $ the matrix $ \iota\left({\mathbit v}\right) $ is similar to a
$ 3 $-diagonal matrix with diagonal entries $ v_{2l} $, above-diagonal entries 1,
and below-diagonal entries $ v_{2l+1}^{2} $. Denote by $ Q_{k+1} $ the set of $ 3 $-diagonal
$ \left(k+1\right)\times\left(k+1\right) $ matrices with the above-diagonal entries being 1. Denote by
$ \widetilde{F}_{\bullet} $ the mapping $ Q_{k+1} \to {\mathfrak P}_{k+1} $ of taking the characteristic polynomial.
Denote the diagonal entries of $ q\in Q $ by $ a_{l} $, $ l=0,\dots ,k $, the below-diagonal
entries by $ b_{l} $, $ l=1,\dots ,k $. Now the proposition is an immediate corollary
of the following

\begin{lemma} The mapping $ \widetilde{F}_{\bullet} $ restricted on the subset $ b_{l}\not=0 $, $ l=1,\dots ,k $, is a
submersion. \end{lemma}

To prove this lemma, denote the characteristic polynomial of the
upper-left principal $ l\times l $ minor by $ d_{l} $. The lemma is an immediate corollary
of

\begin{lemma} The mapping $ \left(d_{k},d_{k+1}\right)\colon Q_{k+1} \to {\mathfrak P}_{k}\times{\mathfrak P}_{k+1} $ restricted on the subset
$ b_{l}\not=0 $, $ l=1,\dots ,k $, is a bijection onto the subset of mutually prime
polynomials $ \left(P_{1},P_{2}\right)\in{\mathfrak P}_{k}\times{\mathfrak P}_{k+1} $. \end{lemma}

This lemma is a direct discrete analogue of the inverse problem for
Sturm--Liouville equation by the spectrum with fixed ends and normalizing
numbers (compare \cite{Lev87Inv}). In fact zeros of $ d_{k+1} $ determine the
spectrum, and
values of $ d_{k} $ at these points determine the normalizing numbers.

\begin{proof} Indeed, extending the sequence $ d_{l} $ by $ d_{0}=1 $, $ d_{-1}=0 $, one can see
that this sequence is uniquely determined by the recurrence relation
\begin{equation}
d_{l}=\left(a_{l-1}-\lambda\right)d_{l-1}-b_{l-1}d_{l-2}.
\notag\end{equation}
From this relation one can immediately see that if $ b_{m} $, $ m<l $, do not
vanish, then $ d_{l} $ and $ d_{l-1} $ are mutually prime. On the other hand, given
mutually prime $ d_{l}\in{\mathfrak P}_{l} $ and $ d_{l-1}\in{\mathfrak P}_{l-1} $, one can uniquely determine $ d_{l-2}\in{\mathfrak P}_{l-2} $
and two numbers $ a_{l-1} $ and $ b_{l-1} $ from the above relation, and $ b_{l-1}\not=0 $. \end{proof}

This finishes the proof of the proposition. \end{proof}

We conclude that at an $ S $-generic point $ {\mathbit v} $ the derivatives $ d\overset{\,\,{}_\circ}{F}_{\lambda}|_{{\mathbit v}} $ span
$ k+1 $-dimensional space (since $ \dim {\mathfrak P}_{k+1}=k+1 $). Now the only condition of
Theorem~\ref{th1.10} (in fact, of Amplification~\ref{amp1.12}) which is missing is
the calculation of the rank of $ \lambda_{1}\left\{,\right\}_{1}+\lambda_{2}\left\{,\right\}_{2} $ for an appropriate $ \lambda_{1} $ and
$ \lambda_{2} $. One can easily see that

\begin{lemma} \label{lm6.50}\myLabel{lm6.50}\relax  The rank of the bracket $ \left\{,\right\}_{1} $ at the point $ {\mathbit v} $ is $ 2k-2d $, here
$ d $ is the number of indices $ l=0,\dots ,k-1 $, such that $ v_{2l+1}=0 $. \end{lemma}

This shows that on the subset $ v_{2l+1}\not=0 $, $ l=0,\dots ,k-1 $, the bihamiltonian
structure satisfies conditions of Theorem~\ref{th1.10} and Amplification
~\ref{amp1.12}, thus is flat indecomposable. This finishes the proof of Theorem
~\ref{th01.60}.

Moreover, since for a flat indecomposable structure both brackets
have corank 1 everywhere, Lemma~\ref{lm6.50} implies that in a neighborhood of
a point $ {\mathbit v} $ with $ v_{2l+1}=0 $ for some $ l=0,\dots ,k-1 $ the bihamiltonian open Toda
structure is {\em not\/} flat indecomposable.

\section{Periodic Toda lattice }\label{h10}\myLabel{h10}\relax 

Recall that $ V_{2k} $ denotes the periodic Toda lattice.

\begin{lemma} The function $ N=v_{1}v_{3}\dots v_{2k-1} $ on $ V_{2k} $ is Casimir w.r.t.~both Poisson
brackets $ \left\{,\right\}_{1} $ and $ \left\{,\right\}_{2} $. \end{lemma}

\begin{proof} Since this function is invariant w.r.t.~translation $ {\mathfrak T}_{\lambda} $, it is
enough to show this for the bracket $ \left\{,\right\}_{2} $. When one calculates $ \left\{N,v_{2l}\right\} $,
only the factor $ v_{2l-1}v_{2l+1} $ of $ N $ matters, and by~\eqref{equ0.10}
$ \left\{v_{2l-1}v_{2l+1},v_{2l}\right\} $ vanishes. Similarly, for $ \left\{N,v_{2l-1}\right\} $ only
$ \left\{v_{2l-3}v_{2l+1},v_{2l-1}\right\} $ matters, and it also vanishes. \end{proof}

Since dimension of $ V_{2k} $ is even, this shows that symplectic leaves of
$ \lambda_{1}\left\{,\right\}_{1}+\lambda_{2}\left\{,\right\}_{2} $ have codimension at least 2. Any hypersurface $ N=\operatorname{const} $ is
decomposed into a union of such leaves for any $ \left(\lambda_{1},\lambda_{2}\right)\not=\left(0,0\right) $. In
particular, each hypersurface $ N=\operatorname{const} $ carries an odd-dimensional
bihamiltonian structure.

\begin{theorem} \label{th10.10}\myLabel{th10.10}\relax  For any $ c\not=0 $ the bihamiltonian structure on the
hypersurface $ N=c $ is generically flat indecomposable. \end{theorem}

Note that this theorem implies Theorem~\ref{th01.70}, since one can
easily modify Theorem~\ref{th2.07} to cover families of bihamiltonian
structures as well:

\begin{amplification} Consider a family of bihamiltonian structures
$ \left(\left\{,\right\}_{1}^{\left(\mu\right)},\left\{\right\}_{2}^{\left(\mu\right)}\right) $ on a manifold $ M $ which depends smoothly on a parameter
$ \mu\in{\mathcal M} $. Suppose that for any $ \mu $ the bihamiltonian structure is flat
indecomposable. Then for any $ m_{0}\in M $ and $ \mu_{0}\in{\mathcal M} $ there is a neighborhood $ U $ of
$ m $, a neighborhood $ U' $ of $ \mu_{0} $ and a family of coordinate system $ \left(x_{i}^{\left(\mu\right)}\right) $ on $ U $
depending smoothly on a parameter $ \mu\in U' $ such that the bihamiltonian
structure $ \left(\left\{,\right\}_{1}^{\left(\mu\right)},\left\{\right\}_{2}^{\left(\mu\right)}\right) $ in the coordinate system $ \left(x_{i}^{\left(\mu\right)}\right) $ is given by
~\eqref{equ45.20} for any $ \mu\in U' $. \end{amplification}

Since the bihamiltonian structure corresponding to $ {\mathcal K}_{1} $ has both
bracket being 0, this amplification implies Theorem~\ref{th01.70}.

\begin{proof}[Proof of Theorem~\ref{th10.10} ] Associate to a point $ {\mathbit v} $ of the infinite Toda
lattice an infinite $ 3 $-diagonal matrix $ \iota\left({\mathbit v}\right) $ in the same way we did it in
Section~\ref{h4}. Consider a matrix equation $ \iota\left({\mathbit v}\right){\mathbit x}=0 $, here $ {\mathbit x}\in{\mathbb C}^{\infty} $ is a
two-side-infinite vector. Since this equation may be written as the
recursion relation
\begin{equation}
v_{2l-1}x_{l-1}+v_{2l}x_{l}+v_{2l+1}x_{l+1}=0,\qquad l\in{\mathbb Z},
\label{equ10.10}\end{equation}\myLabel{equ10.10,}\relax 
this matrix equation has a two-dimensional space of solutions if $ v_{2l-1}\not=0 $
for any $ l\in{\mathbb Z} $.

If $ {\mathbit v} $ is in the periodic Toda lattice, then the equation $ \iota\left({\mathbit v}\right){\mathbit x}=0 $ is
invariant with respect to the shift $ x_{l} \mapsto x_{l+k} $ of coordinates of $ {\mathbit x} $. This
shift induces a linear transformation $ {\mathcal M}={\mathcal M}\left({\mathbit v}\right) $ of monodromy in the
$ 2 $-dimensional vector space of solutions. As in Section~\ref{h0}, denote by $ {\mathbit v}^{0} $
an element of $ {\mathbb C}^{\infty} $ with 1 on even positions, 0 on odd positions.

\begin{lemma} If $ v_{2l-1}\not=0 $ for any $ l\in{\mathbb Z} $, then $ \det  {\mathcal M}=1 $, and $ \operatorname{Tr} {\mathcal M}\left({\mathbit v}-\lambda{\mathbit v}^{0}\right) $ is a
polynomial of degree $ k $ in $ \lambda $ with the leading coefficient $ N^{-1} $. \end{lemma}

\begin{proof} Indeed, the recursion~\eqref{equ10.10} induces a linear transformation
$ \left(x_{l},x_{l+1}\right)=m_{l}\left(x_{l-1},x_{l}\right)/v_{2l+1} $, $ m_{l}=\left(
\begin{matrix}
0 & v_{2l+1} \\ -v_{2l-1} & -v_{2l}
\end{matrix}
\right) $. In an appropriate
basis $ N\cdot{\mathcal M} $ can be written as $ m_{k}m_{k-1}\dots m_{1} $, and each matrix $ m_{l}=m_{l}\left({\mathbit v}\right) $ has
determinant $ v_{2l-1}v_{2l+1} $. Moreover, $ m_{l}\left({\mathbit v}-\lambda{\mathbit v}^{0}\right) $ is of degree 1 in $ \lambda $ with the
leading term $ \left(
\begin{matrix}
0 & 0 \\ 0 & \lambda
\end{matrix}
\right) $.

Thus $ N\cdot{\mathcal M}\left({\mathbit v}-\lambda{\mathbit v}^{0}\right) $ is a polynomial in $ \lambda $ of degree $ k $ with the leading term
being $ \left(
\begin{matrix}
0 & 0 \\ 0 & \lambda^{k}
\end{matrix}
\right) $, which finishes the proof. \end{proof}

\begin{lemma} \label{lm10.50}\myLabel{lm10.50}\relax  The function $ \operatorname{Tr} {\mathcal M}\left({\mathbit v}\right) $ defined on the open subset $ v_{2l-1}\not=0 $,
$ l=1,\dots ,k $, of $ V_{2k} $ is a Casimir function for the Poisson bracket $ \left\{,\right\}_{2} $. \end{lemma}

We do not prove this standard statement about the periodic Toda
lattice. As in the case of Lemma~\ref{lm4.05}, the proof is reduced to a check
of a finite number of identities.

The following lemma is obvious:

\begin{lemma} On the open subset $ v_{2l-1}\not=0 $, $ l=1,\dots ,k $, of $ V_{2k} $ the Poisson
bracket~\eqref{equ0.20} has symplectic leaves of codimension 2 given by
the equations $ v_{0}+v_{2}+\dots +v_{2k-2}=C_{0} $, $ v_{1}v_{3}\dots v_{2k-1}=C_{1} $. \end{lemma}

This shows that $ r=2 $ in Proposition~\ref{prop6.15}.

To demonstrate Theorem~\ref{th10.10} the only thing which remains to be
proved is that at a generic point $ {\mathbit v}\in V_{2k} $ the differentials $ d \operatorname{Tr} {\mathcal M}\left({\mathbit v}-\lambda{\mathbit v}^{0}\right)|_{{\mathbit v}} $
for different $ \lambda\in{\mathbb C} $ and the differential of $ N\equiv v_{1}v_{3}\dots v_{2k-1} $ span a
$ k+1 $-dimensional vector subspace of $ {\mathcal T}_{{\mathbit v}}^{*}V_{2k} $. It is enough to show that for
a generic $ {\mathbit v} $ the differentials of $ N\cdot{\mathcal M}\left({\mathbit v}-\lambda{\mathbit v}^{0}\right) $ for different $ \lambda\in{\mathbb C} $
span a $ k $-dimensional vector subspace of the hyperplane $ d\left(v_{1}v_{3}\dots v_{2k-1}\right)=0 $
in $ {\mathcal T}_{{\mathbit v}}^{*}V_{2k} $.

The leading coefficient in $ \lambda $ of $ N\cdot\operatorname{Tr} {\mathcal M}\left({\mathbit v}-\lambda{\mathbit v}^{0}\right) $ is 1, thus the
function $ N\cdot\operatorname{Tr} {\mathcal M}\left({\mathbit v}-\lambda{\mathbit v}^{0}\right)-\lambda^{k} $ defines a mapping $ {\mathfrak M}\colon V_{2k} \to {\mathcal P}_{k-1} $. Again, it is
enough to show that the restriction of this polynomial mapping to
$ H_{c}=\left\{v_{1}v_{3}\dots v_{2k-1}=c\right\} $ is a submersion for a generic $ {\mathbit v} $ and $ c\not=0 $. On the other
hand, multiplication of $ v_{i} $ by the same non-zero constant does not change
$ {\mathcal M}\left({\mathbit v}\right) $, thus if we prove this statement for one $ c\not=0 $, is it true for any
$ c\not=0 $. Thus it is enough to demonstrate this statement for $ c\approx0 $, $ c\not=0 $. Again,
it is enough to show that the restriction of $ {\mathfrak M} $ to an open subset of $ c=0 $
is a submersion.

However, if $ v_{1}=v_{2}=\dots =v_{2k-1}=0 $, then
\begin{equation}
\lambda^{k}+{\mathfrak M}\left({\mathbit v}\right)=\left(\lambda-v_{2}\right)\left(\lambda-v_{4}\right)\dots \left(\lambda-v_{2k}\right),
\notag\end{equation}
thus the restriction of $ {\mathfrak M} $ to $ \left\{v_{1}=v_{2}=\dots =v_{2k-1}=0\right\} $ is a surjection, thus is
a submersion in a generic point. This shows that Theorem~\ref{th1.10} is
applicable, thus the bihamiltonian structure is indeed flat
indecomposable at a generic point. \end{proof}

\section{Lax structures }\label{h60}\myLabel{h60}\relax 

The following definition is inspired by \cite{KosMag96Lax}. In this paper
a notion of a Lax operator is introduced, this is a matrix-valued
function on a bihamiltonian structure which satisfies some compatibility
relations. However, since these relations are expressed in terms of the
characteristic polynomial of the matrix, it is more convenient to work
directly with the mapping into polynomials.

Recall that $ {\mathcal P}_{n} $ was defined in Section~\ref{h45}. Denote the value at $ \lambda $ of
a polynomial $ p\in{\mathcal P}_{n} $ by $ p|_{\lambda} $.

\begin{definition} Consider a bihamiltonian structure $ \left(M,\left\{,\right\}_{1},\left\{,\right\}_{2}\right) $. Consider a
mapping $ {\mathbit L} $ from $ M $ to the set $ {\mathcal P}_{n-1} $ of polynomials of degree $ n-1 $. This
mapping is a {\em weak Lax structure\/} on $ M $ of rank $ n $ if for any $ \lambda\in{\mathbb R} $ the
function $ C_{\lambda} $ on $ M $ defined by $ m \mapsto {\mathbit L}\left(m\right)|_{\lambda} $ is a Casimir function for
$ \lambda\left\{,\right\}_{1}+\left\{,\right\}_{2} $.

Consider a point $ m_{0}\in M $. Suppose that the action dimension of $ M $ at
$ m_{0}\in M $ is $ n $. A {\em Lax structure\/} on $ M $ near $ m_{0} $ is a weak Lax structure $ {\mathbit L} $ of rank
$ n $ such that the mapping $ {\mathbit L} $ is a submersion. \end{definition}

Note that if the bihamiltonian structure is in fact analytic, then
$ C_{\lambda} $ is Casimir for complex $ \lambda $ too (since the conditions of being a
$ \lambda $-Casimir family are polynomial in $ \lambda $).

\begin{theorem} \label{th60.30}\myLabel{th60.30}\relax  If an analytic bihamiltonian structure on $ M $ admits a
Lax structure near $ m_{0}\in M $, and for one particular $ \left(\lambda_{1},\lambda_{2}\right)\in{\mathbb C}^{2} $ the Poisson
structure $ \lambda_{1}\left\{,\right\}_{1}+\lambda_{2}\left\{,\right\}_{2} $ has a constant corank 1, then the bihamiltonian
structure is a Kronecker structure of type $ \left(\dim  M\right) $ near $ m_{0} $. \end{theorem}

In other words, the manifold $ M $ is odd-dimensional and one can find a
local coordinate system where both brackets have constant coefficients
and are given by~\eqref{equ45.20}. In particular, all such bihamiltonian
structures of the same dimension are locally isomorphic.

\begin{proof} Reduce this statement to one of Amplification~\ref{amp1.12}.

In our case $ d=n-1 $, and, by submersion condition, $ \dim  W_{1}=n $. Thus the
only thing one needs to show is that $ \dim  M=2n-1 $. This momentarily follows
from the definition of the action dimension. \end{proof}

\begin{remark} In applications the Poisson bracket $ \left\{,\right\}_{1} $ usually has a much
simpler form than $ \left\{,\right\}_{2} $, thus most of the time one would check the rank
condition for the bracket $ \left\{,\right\}_{1} $. (Recall that for Kronecker structures {\em all\/}
the nonzero linear combinations of brackets have the same rank.) \end{remark}

Let us spell out the relation of our definition with one of
\cite{KosMag96Lax}. Consider the Newton symmetric functions $ s_{k}=\sum_{i}\lambda_{i}^{k} $ of roots
$ \left\{\lambda_{i}\right\} $ of polynomial $ \lambda^{n}+p\left(\lambda\right) $, $ p\in{\mathcal P}_{n-1} $ as functions on $ {\mathcal P}_{n-1} $, let
$ H_{k-1}\buildrel{\text{def}}\over{=}s_{k}\circ{\mathbit L}/k $. Then the conditions of \cite{KosMag96Lax} are that $ H_{k} $, $ k\geq0 $,
satisfy Lenard recursion relations~\eqref{equ48.20}. As in Section~\ref{h48},
consider a formal power series $ c\left(t\right)=\sum_{k\geq1}s_{k}t^{1-k}/k $ in $ t^{-1} $ with coefficients
in functions on $ {\mathcal P}_{n-1} $. Then $ e^{-c\left(t\right)/t}=\Pi_{i}\left(1-\lambda_{i}/t\right)=t^{-n}\left(t^{n}+p\left(t\right)\right) $, $ p\in{\mathcal P}_{n-1} $. Let
$ C\left(t\right)=\sum_{k\geq0}H_{k}t^{-k} $, then $ e^{-C\left(t\right)/t}|_{m}=t^{-n}{\mathbit L}\left(m\right)|_{t} $ for any $ m\in M $.

Since the latter expression is a formal series in $ t^{-1} $ with a finite
number of non-zero coefficients, it is an anchored formal $ \lambda $-family iff it
is a $ \lambda $-Casimir family. Since for any function $ \alpha $ of one variable $ \alpha\left(C\right) $ is a
Casimir function if $ C $ is such, we conclude that $ {\mathbit L}|_{t} $ is a $ \lambda $-Casimir family
iff $ C\left(t\right) $ is an anchored formal $ \lambda $-family. Thus the condition that $ {\mathbit L} $ is a
weak Lax structure is equivalent to the pair of conditions: of $ H_{k} $
satisfying Lenard recursion relations~\eqref{equ48.20}, {\em and additionally\/} of $ H_{0} $
being a Casimir function for $ \left\{,\right\}_{1} $. This shows

\begin{proposition} Suppose that $ L\colon M \to \operatorname{Mat}\left(n\right) $ is a Lax operator in the
sense of \cite{KosMag96Lax}. Let $ {\mathbit L} $ be the mapping
\begin{equation}
M \to {\mathcal P}_{n-1}\colon m \mapsto \det \left(t\boldsymbol1-L\left(m\right)\right)-t^{n}.
\notag\end{equation}
Then $ {\mathbit L} $ is a weak Lax structure iff $ \operatorname{Tr} L $ is a Casimir function for $ \left\{,\right\}_{1} $.

\end{proposition}

As in Section~\ref{h55}, note that in applications the Lenard scheme is
most frequently used when $ \operatorname{Tr} L $ is a Casimir function for $ \left\{,\right\}_{1} $. Note also
that one can consider a weak Lax structure as an ``anchored'' variant of a
Lax operator of \cite{KosMag96Lax} (compare with Remark~\ref{rem48.02} and
Definition~\ref{def48.25}).

\begin{remark} By Theorem~\ref{th01.60}, in conditions of Theorem~\ref{th60.30} an open
subset of the bihamiltonian structure is locally isomorphic to the
structure of Toda lattice. This isomorphism provides the subset $ U $ with a
Lax operator in the most usual sense of this word, i.e., with a mapping
$ L\colon U \to \operatorname{Mat}\left(n\right) $ such that for any action function\footnote{See Section~\ref{h02}.} $ H $ on $ U $ there is a
mapping $ A_{H}\colon U \to \operatorname{Mat}\left(n\right) $ such that $ H $-Hamiltonian flow on $ U $ corresponds to
$ \frac{dL}{dt}=\left[A_{H},L\right] $.

In other words, Theorem~\ref{th60.30} provides a partial explanation for
the relation between Lax operator and Lax--Nijenhuis operators discovered
in \cite{KosMag96Lax}. \end{remark}

\begin{remark} Note that the conditions of Theorem~\ref{th60.30} break into four
separate parts: the condition of being a weak Lax structure, the
condition that coefficients of $ {\mathbit L} $ provide enough functions to completely
integrate $ M $, the submersion condition, and the condition of having small
corank. Note that the corank of the structure cannot be less than 1,
since we {\em require\/} existence of Casimir function for any $ \lambda $. Thus two last
conditions taken together may be interpreted as conditions of
non-degeneracy of the Lax structure. \end{remark}

\begin{nwthrmi} Which conditions on a weak Lax family imply that the
bihamiltonian structure is Kronecker at generic points? \end{nwthrmi}

Conjecture~\ref{con01.100} claims that many bihamiltonian structures
which admit a Lax structure are in fact Kronecker at generic points. An
answer on the above question might have provided a better understanding
for the statement of Conjecture~\ref{con01.100}.

\section{Geometric conjectures }\label{h005}\myLabel{h005}\relax 

Note that the Theorems~\ref{th01.60},~\ref{th01.70}, and~\ref{th60.30} run against
the common intuition, which says that integrable systems should be
expressed as direct products of two-dimensional blocks. However, this
point of view comes from the symplectic approach to integrable systems,
where everything is {\em forced\/} to be even-dimensional.

The above theorems show that this common intuition has historical
roots only, and some new type of intuition for geometric approach to
integrable systems may be needed.

Our meta-conjecture is that the mindset of ``everything is a product
of odd-dimensional components (given by~\eqref{equ45.20})'' is much more
appropriate for the geometric study of bihamiltonian structures, compare
with Remark~\ref{rem01.95} and Conjecture~\ref{con01.100}.

Again, if one believes in the above meta-conjecture, one can see
that the Procrustean approach of symplectic geometry forces a reduction
of dimension (as in Remark~\ref{rem0.20}, which gives an analogue of
restriction to a hypersurface), which reduces a feature-rich
bihamiltonian structure to a non-rigid symplectic structure.

\begin{remark} \label{rem01.95}\myLabel{rem01.95}\relax  Definition~\ref{def01.105} provides an example of micro-local
approach to bihamiltonian systems. By Theorem~\ref{th6.10}, in each tangent
space any bihamiltonian structure decomposes into a direct sum of Jordan
blocks and Kronecker blocks. Thus a natural question arises: given a
bihamiltonian structure $ M $, which indecomposable pairs $ {\mathcal J}_{2k,\lambda} $ and $ {\mathcal K}_{2k-1} $
appear at which points of $ M? $

Theorems~\ref{th01.60} and~\ref{th01.70} answer this question for generic
points of the open and the periodic {\em Toda lattice}. We think we can answer
this question\footnote{After the initial release of this paper M.~Gekhtman explained us that
the result on the open Toda lattice implies the statements about
the open odd-dimensional Kac--van~Moerbeke--Volterra lattice, as well as a
similar statement about the open relativistic Toda lattice \cite{Ruj90Rel}.
\endgraf
This is an immediate corollary of the existence of local
isomorphisms of these bihamiltonian systems similar to those constructed in
\cite{DeiLi91Poi,Dam94Mul}, see \cite{GekhShap99Non} and \cite{FayGekh99Ele}.} for generic points of the odd-dimensional open or
even-dimensional periodic {\em Kac\/}--{\em van\/}~{\em Moerbeke\/}--{\em Volterra system\/}
\cite{Kacvan75Exp,FerSan97Int}, of the {\em full Toda lattice\/} \cite{Kos79Sol},
and of the multidimensional {\em Euler top\/} \cite{MorPiz96Eul}. In tangent spaces at
generic points the open Toda lattice is an indecomposable Kroneker block,
the periodic Toda lattice is a direct product of indecomposable
$ 1 $-dimensional and $ 2k-1 $-dimensional Kroneker blocks. The complete Toda
lattice and the multidimensional Euler top are products of Kroneker
blocks with the dimensions of components being $ \left(2k-1,2k-3,2k-5,\dots \right) $ and
$ \left(2k-1,2k-5,2k-9,\dots \right) $ correspondingly.

Additionally, results of \cite{Pan99Ver} show that a similar
decomposition exists for the regular case of Example~\ref{ex002.45}. In this
case the dimensions of components have the form $ 2e_{1}-1,\dots ,2e_{r}-1 $, $ e_{i} $ being
the exponents of the Weyl group of $ {\mathfrak g} $, $ r $ being the rank of $ {\mathfrak g} $. \end{remark}

The above descriptions of tangent spaces together with Theorems
~\ref{th01.60} and~\ref{th01.70} suggest the following

\begin{conjecture} \label{con01.100}\myLabel{con01.100}\relax  The odd-dimensional open Volterra system, the
even-dimensional periodic Volterra system, the full Toda lattice, the
multidimensional Euler top, and the regular case of Example~\ref{ex002.45}
are\footnote{Paper \cite{Zakh99Kro} contains a proof of the part of the conjecture related
to Example~\ref{ex002.45}, see the previous footnote for some other cases.} generically Kronecker bihamiltonian structures. \end{conjecture}

As shown in this paper, the powerful methods of \cite{GelZakhWeb,%
GelZakh93} are enough to translate some simple properties\footnote{The existence of Casimir functions given by Lemmas~\ref{lm4.05} and~\ref{lm10.50}.} of the open
and the periodic Toda lattices into description of the {\em local\/} geometry of
these structures. One may hope that it is possible to generalize the
results of \cite{GelZakhWeb,GelZakh93} so that they cover structures with
geometry of tangent spaces as in Remark~\ref{rem01.95}. This would allow one
to prove Conjecture~\ref{con01.100} using some simple results about these
integrable systems\footnote{Again, since the geometry of these system is very well investigated, it
may be possible to prove this conjecture directly using appropriate
systems of action-angle variables for these manifolds.
\endgraf
However, an approach based on Conjecture~\ref{con01.110} would allow one to
prove Conjecture~\ref{con01.100} using only simple-to-obtain action
variables, i.e., families of Hamiltonians for the above manifolds.}.

Using language of Section~\ref{h02}, one can state such conjectures in the
following form.

\begin{conjecture} \label{con01.110}\myLabel{con01.110}\relax  Suppose that two bihamiltonian structures
$ \left(M,\left\{\right\}_{1},\left\{\right\}_{2}\right) $ and $ \left(M',\left\{\right\}'_{1},\left\{\right\}'_{2}\right) $ are both homogeneous. Consider webs\footnote{See Section~\ref{h02}.} $ {\mathcal B}_{U} $
and $ {\mathcal B}_{U'} $ which correspond to small open subsets $ U\subset M $, $ U'\subset M' $. If webs $ {\mathcal B}_{U} $ and
$ {\mathcal B}_{U'} $ are locally isomorphic, then the bihamiltonian structures on $ M $ and $ M' $
are locally isomorphic. In particular, the types of $ M $ and $ M' $ coincide. \end{conjecture}

This conjecture may be augmented by the following description
of webs for homogeneous structures \cite{Pan99Ver}:

\begin{proposition} The web $ {\mathcal B}_{U} $ corresponding to a small open subset $ U $ of
homogeneous bihamiltonian structure of type $ \left(2k_{1}-1,2k_{2}-1,\dots ,2k_{l}-1\right) $ is a
manifold of dimension $ k_{1}+k_{2}+\dots +k_{l} $, and the subspace $ {\mathfrak C}_{\lambda} $ of the space of
functions on $ {\mathcal B}_{U} $ consists of local equations of a foliation $ {\mathcal F}_{\lambda} $ on $ {\mathcal B}_{U} $ of
codimension $ l $. \end{proposition}

Conjecture~\ref{con01.110}, together with Amplification~\ref{amp1.07}, lead to
the following

\begin{conjecture} \label{con01.120}\myLabel{con01.120}\relax  Consider a manifold $ M $ with two compatible Poisson
structures $ \left\{,\right\}_{1} $ and $ \left\{,\right\}_{2} $. Consider a finite set $ L $ with $ r $ elements.
Consider families of smooth functions $ F_{l,\lambda} $, $ l\in L $, $ \lambda\in{\mathbb C} $, on $ M $ such that for
any $ l\in L $ and any $ \lambda\in{\mathbb C} $ the function $ F_{l,\lambda} $ is Casimir w.r.t.~the Poisson
bracket $ \lambda\left\{,\right\}_{1}+\left\{,\right\}_{2} $. Suppose that $ F_{l,\lambda} $ depends polynomially on $ \lambda $
\begin{equation}
F_{l,\lambda}\left(m\right)= \sum_{k=0}^{d_{l}}f_{l,k}\left(m\right)\lambda^{k},
\notag\end{equation}
with smooth coefficients $ f_{l,k}\left(m\right) $. For $ m\in M $ denote by $ W_{1}\left(m\right)\subset{\mathcal T}_{m}^{*}M $ the vector
subspace spanned by the the differentials $ df_{l,k}|_{m} $ for all possible $ l $ and
$ 0\leq k\leq d_{l} $. If
\begin{enumerate}
\item
for one particular value $ m_{0}\in M $ one has $ \dim  W_{1}\left(m_{0}\right)\geq\frac{\dim  M+r}{2} $;
\item
for one particular value of $ \lambda_{1},\lambda_{2}\in{\mathbb C}^{2} $ the Poisson structure
$ \lambda_{1}\left\{,\right\}_{1}+\lambda_{2}\left\{,\right\}_{2} $ has at most $ r $ independent Casimir functions on any open
subset of $ M $ near $ m_{0} $;
\item
the degrees $ d_{l} $ satisfy $ \sum_{L}\left(2d_{l}+1\right)\leq\dim  M $;
\end{enumerate}
then $ \dim  M-r $ is even, $ \dim  W_{1}\left(m_{0}\right)=\frac{\dim  M+r}{2} $, the degrees $ d_{l} $ satisfy
$ 2\sum_{L}d_{l}+r=\dim  M $, and the bihamiltonian structure on $ M $ is Kronecker of
type $ \left(2d_{1}+1,\dots ,2d_{r}+1\right) $ on an open subset $ U\subset M $ such that $ m_{0} $ is in the
closure of $ U $. \end{conjecture}

Conjecture~\ref{con01.120} immediately implies Conjecture~\ref{con01.100},
since the explicit formulae for Hamiltonians for the dynamic systems of
Conjecture~\ref{con01.100} are well-known and may be included into families as
in Conjecture~\ref{con01.120}.

To understand the significance of Conjecture~\ref{con01.120}, note that
by Remark~\ref{rem6.13} all the Kronecker structures of the given type are
locally isomorphic, and obviously satisfy the conditions of the
conjecture. Thus this conjecture provides a criterion of being a
Kronecker structure in terms of the mutual position of Casimir functions
for the combinations of brackets of bihamiltonian structure.

\begin{conjecture} In the settings of Conjecture~\ref{con01.120} if one supposes
that the Poisson structure $ \lambda_{1}\left\{,\right\}_{1}+\lambda_{2}\left\{,\right\}_{2} $ has constant corank $ r $, then one
may weaken the condition on $ \dim  W_{1} $ to become $ \dim  W_{1}\left(m_{0}\right)\geq\frac{\dim  M+r-1}{2} $,
and amplify the conclusion to so that the open subset $ U $ contains $ m_{0} $. \end{conjecture}

The above theorems and conjectures lead one to the following

\begin{nwthrmii} Why each ``classical'' finite-dimensional bihamiltonian
structure has an open subset which is Kronecker, or may be ``naturally''
considered as a reduction of dimension starting from a larger
bihamiltonian structure which is Kronecker? \end{nwthrmii}

This question is amplified by the fact that in \cite{GelZakhWeb,%
GelZakh93} we constructed a huge family of non-Kronecker integrable
bihamiltonian structures (see also examples in Section~\ref{h47} for the
dimension being 3). Such integrable systems are {\em actually nonlinear}, as
opposed to {\em manifestly nonlinear\/} systems, which may become linear after an
appropriate coordinate change (compare with Definition~\ref{def002.43}). One
would see that an answer to the above question would unravel some
mechanism by which the actually nonlinear integrable systems avoid
attention of mathematical physicists.

Note that Theorem~\ref{th01.60} allows one to restate the above question
using direct products of open Toda lattices instead of Kronecker
structures:

{\em Why many ``classical'' bihamiltonian structures are (in generic points)
locally isomorphic to direct products of open Toda lattices?\/}

While Section~\ref{h60} singles out flat indecomposable structures as
those which admit non-degenerate Lax structures, we do not consider this
as a legitimate explanation to the above selection principle. Lax
representation is only one of multiple approaches to integration of
dynamical systems, so explaining the above selection principle by using
Theorem~\ref{th60.30} just substitutes one question (why all the classical
systems are flat) by another one (why all the classical systems admit Lax
representation).

\begin{remark} Note that a flat bihamiltonian structure of dimension $ d $ may be
extended locally to a $ d\left(d-1\right)/2 $-parametric linear family of Poisson
structures: those which have constant coefficients in the above
coordinate system. Our meta-conjecture about the r\^ole of Kronecker
structures may explain an abundance of multi-hamiltonian structures in
mathematical physics (for example, see \cite{DamPasSok95Tri,OlvRos96Tri,%
TsuTakKaj97Tri,BlaFerGom98Som}).\footnote{However, note that what is commonly called a ``multi-hamiltonian''
structure is frequently just a figure of speech: the additional
``brackets'' which augment the bihamiltonian structure are not only not
Poisson (thus do not satisfy Jacobi condition), but not even brackets
(thus $ \left\{f,g\right\}_{3} $ would be defined for {\em some\/} $ f $ and $ g $ only).} \end{remark}

\bibliography{ref,outref,mathsci}
\end{document}